\documentclass[11pt,a4paper]{article}
\usepackage[T1]{fontenc}
\usepackage{lmodern}
\usepackage[margin=26mm]{geometry}
\usepackage{amsmath,amssymb,amsthm,mathtools}
\usepackage[numbers,sort]{natbib}
\usepackage[hidelinks,unicode]{hyperref}
\allowdisplaybreaks[2]
\newtheorem{theorem}{Theorem}[section]
\newtheorem{lemma}[theorem]{Lemma}
\newtheorem{proposition}[theorem]{Proposition}
\newtheorem{corollary}[theorem]{Corollary}
\theoremstyle{definition}
\newtheorem{assumption}[theorem]{Assumption}
\newtheorem{definition}[theorem]{Definition}
\theoremstyle{remark}
\newtheorem{remark}[theorem]{Remark}
\numberwithin{equation}{section}
\newcommand{\R}{\mathbb R}
\newcommand{\E}{\mathbb E}
\newcommand{\esssup}{\operatorname*{ess\,sup}}

\hypersetup{pdftitle={Uniform Spatial Schauder Estimates for General Backward Stochastic Partial Differential Equations with Random Coefficients},pdfauthor={Yunzhang Li, Shanjian Tang, Yuyang Ye}}
\title{Uniform Spatial Schauder Estimates for General Backward Stochastic Partial Differential Equations with Random Coefficients\footnote{This work is partially supported by National Natural Science Foundation of China (Nos. 12631016 and 11631004), and by Key Laboratory of Mathematics for Nonlinear Sciences (Fudan University), Ministry of Education, Handan Road 220, Shanghai 200433, China.}}
\author{
Yunzhang Li\thanks{Research Institute of Intelligent Complex Systems, Fudan University, Shanghai 200433, China. \textit{E-mail:} \texttt{li\_yunzhang@fudan.edu.cn}.}
\and
Shanjian Tang\thanks{Department of Finance and Control Sciences, School of Mathematical Sciences, Fudan University, Shanghai 200433, China. \textit{E-mail:} \texttt{sjtang@fudan.edu.cn}.}
\and
Yuyang Ye\thanks{Corresponding author. Department of Finance and Control Sciences, School of Mathematical Sciences, Fudan University, Shanghai 200433, China. \textit{E-mail:} \texttt{yyye22@m.fudan.edu.cn}.}}
\date{}
\begin{document}
\maketitle
\begin{abstract}
We prove spatial Schauder estimates for super-parabolic backward stochastic partial differential equations. All coefficients can be random and depend on time and space, and the estimates are uniform in  $(\omega, t)$ . The key ingredient is a weighted conditional energy norm adapted to the parabolic scale of each frequency block. Together with heat-kernel bounds and coefficient freezing, it controls the gradient-martingale coupling without assuming a bounded martingale integrand. As an application, we establish Schauder solvability for semilinear stochastic Hamilton--Jacobi--Bellman equations.
\end{abstract}
\noindent\textbf{Keywords:} backward stochastic partial differential equation, Schauder estimate, local Sobolev estimates, coefficient freezing, stochastic Hamilton--Jacobi--Bellman equation.

\section{Introduction}\label{sec:intro}
We consider backward stochastic partial differential equations (BSPDEs) of the general linear form
\begin{equation}\label{eq:full}
\left\{
\begin{aligned}
&-du=\bigl[a^{ij}\partial_{ij}u+\sigma^{ik}\partial_i v^k
             +b^i\partial_i u+cu+\gamma^kv^k+f\bigr]dt\\
&\qquad-v^k\,dW_t^k,\qquad (t,x)\in[0,T]\times\mathbb T^n,\\
&u(T,x)=g(x),\qquad x\in\mathbb T^n.
\end{aligned}
\right.
\end{equation}
Here $\mathbb T^n:=(\R/2\pi\mathbb Z)^n$ is the $n$-dimensional torus, and $W$ is an $m$-dimensional Brownian motion. All coefficients are random and depend on time and space. We write $\partial_i:=\partial_{x_i}$ and $\partial_{ij}:=\partial_{x_i}\partial_{x_j}$, and use $D:=(\partial_1,\ldots,\partial_n)$ and $D^2:=(\partial_{ij})_{i,j=1}^n$ for the gradient and Hessian. We use the  convention of Einstein summation for repeated spatial indexes $i,j=1,\ldots,n$ and noise indexes $k=1,\ldots,m$.

In the theory of stochastic control, BSPDEs naturally arise as adjoint equations in the maximum principle of optimal control,  and as dual equations for controlled stochastic partial differential equations. Bensoussan~\cite{Bensoussan1983PartialObservation,Bensoussan1983Distributed} used them in the maximum principles for partially observed optimal control of diffusions,  and for optimal control of stochastic distributed parameter systems. Zhou~\cite{Zhou1993} extended this adjoint approach to controlled SPDEs with gradient noise. Fractional adjoint BSPDEs also appear in the control of partially observed systems driven by stable L\'evy processes~\cite{Ye2025}. In controllability problems, observability estimates for backward adjoint equations yield null controllability results for stochastic parabolic equations through duality~\cite{TangZhang2009}.

BSPDEs also characterize adapted random fields used in probabilistic representations and stochastic transformations. In mathematical finance, Ma and Yong~\cite[Section~8]{Ma1997} use BSPDE solutions as decoupling fields for pricing and hedging under random coefficients. Further applications include swing option pricing~\cite{BenderDokuchaev2016Swing} and option pricing under rough volatility~\cite{Bayer2022}. Reflected BSPDEs and related variational inequalities describe optimal stopping~\cite{QiuWei2014} and Dynkin games~\cite{YangTang2013}. BSPDEs also provide Zvonkin transformations for proving well-posedness of SDEs with random, irregular coefficients~\cite{Zhao2022}.

Stochastic Hamilton--Jacobi--Bellman (HJB) equations provide a direct motivation for the random principal coefficients and derivative coupling in~\eqref{eq:full}. Peng's pioneering work~\cite{Peng1992} introduced these equations as dynamic programming equations for control problems with random coefficients. When the state diffusion is random and depends on the state, the principal coefficients $a,\sigma$ are random and spatially dependent. The drift term $\sigma^{ik}\partial_i v^k$, involving spatial derivatives of the unknown martingale integrand $v$, arises from the quadratic covariation in the It\^o--Wentzell formula when the state process and the value function share the same Brownian noise. For control-independent diffusion, the HJB equation is semilinear, and fixing a feedback control yields a linear BSPDE with the same diffusion and derivative coupling~\cite{MollaQiu2026}. The same coupling appears in the common-noise mean field game system of Cardaliaguet, Delarue, Lasry, and Lions~\cite[Chapter~4]{CDLL2019}, where a stochastic HJB equation for the value function is coupled with a forward Fokker--Planck equation for the conditional population distribution.

When all coefficients and data are deterministic, the martingale integrand vanishes ($v=0$), and~\eqref{eq:full} reduces to the parabolic terminal-value problem
\[
-\partial_tu=a^{ij}\partial_{ij}u+b^i\partial_i u+cu+f,
\qquad u(T)=g.
\]
For uniformly parabolic equations with coefficients measurable in time and H\"older continuous in space, Schauder theory gives (see Lorenzi~\cite{Lorenzi2000} for the whole-space case)
\[
\|u\|_{L^\infty_tC_x^{2+\alpha}}
+\|\partial_tu\|_{L^\infty_tC_x^\alpha}
\le C\bigl(\|f\|_{L^\infty_tC_x^\alpha}+\|g\|_{C_x^{2+\alpha}}\bigr).
\]

With random coefficients or data, an adapted solution of~\eqref{eq:full} generally requires a martingale term, as in a backward stochastic differential equation (BSDE). When the coefficients and data are independent of space, equation~\eqref{eq:full} reduces to a linear BSDE. The martingale term generally prevents classical differentiability in time. We ask whether spatial H\"older bounds uniform in  $(\omega, t)$  hold with random principal coefficients and gradient-martingale coupling.

The Sobolev theory of BSPDEs is well developed; see, for example,~\cite{Ma1999,Hu2002,Du2010,DuQiuTang2012,Du2012Dirichlet}. Schauder estimates are comparatively less developed. Finite-moment Schauder estimates were established by Du and Liu~\cite{Du2016,Du2019} for forward SPDEs and by Tang and Wei~\cite{tang2016cauchy} for BSPDEs. Uniform spatial H\"older estimates of the type sought here are available in several mean field game settings. Cardaliaguet, Delarue, Lasry, and Lions~\cite[Chapter~4]{CDLL2019} and Cardaliaguet, Maillet, and Yan~\cite[Proposition~3.1]{CardaliaguetMailletYan2025} obtained such estimates for equations with constant principal coefficients in the common-noise setting. In major--minor mean field games, Delarue and Mou~\cite[Lemma~3.1 and Proposition~3.6]{DelarueMou2025} obtained uniform spatial H\"older bounds together with BMO estimates for the spatially centered martingale integrand.

We give spatial Schauder estimates uniform in  $(\omega, t)$, much stronger than those  of previous BSPDE theory~\cite{tang2016cauchy} for essentially bounded data. We extend the existing uniform estimates for the case of constant principal coefficients in the common-noise mean field game setting~\cite{CDLL2019,CardaliaguetMailletYan2025} to the case of  random, spatially varying principal coefficients, with gradient-martingale coupling and general lower-order terms. Under uniform super-parabolicity, we obtain the full gain of two spatial derivatives from $C_x^\alpha$ forcing to $C_x^{2+\alpha}$ solutions for $C_x^{2+\alpha}$ terminal data, together with a conditional energy bound for the martingale integrand. We also prove that the $C_x^{1+\alpha}$ regularity assumption on both $\sigma$ and $\gamma$ in equation~\eqref{eq:full} is sharp, by counterexamples with fixed deterministic coefficients and uniformly bounded $C_x^{2+\alpha}$ terminal data.

The above uniform bounds fail to hold in general for forward SPDEs: the exogenous  stochastic forcing can produce solutions with an infinite essential supremum over $\omega$, even for bounded smooth initial data and noise coefficients.  In contrast to that of general forward SPDEs, the spatial regularity of the BSPDE~\eqref{eq:full} is more like that of deterministic parabolic PDEs.

A challenging difficulty is that the unknown martingale integrand $v$ might be unbounded even for bounded and spatially invariant terminal data. Sobolev or finite-moment estimates alone do not provide the required essential bound over  $(\omega, t)$. Related conditional energy estimates with Sobolev weights on Fourier modes are used by Delarue and Mou~\cite[proof of Proposition~3.6]{DelarueMou2025} to obtain BMO bounds. Here we introduce a weighted conditional energy norm for the frequency blocks of $v$, localized at spatial scale $r$ and time scale $r^2$. This norm links weighted It\^o estimates for $v$ with heat-kernel and perturbation estimates for $u$, so that each unknown can be controlled in terms of the other and the data. This norm also controls the heat-kernel terms involving $v$ under random spatial shifts. Combining these bounds yields the uniform Schauder estimate with a gain of two spatial derivatives from the forcing.

As an application, we establish Schauder solvability for semilinear stochastic HJB equations with random state-varying diffusion and a compact control set. The drift and cost coefficients need only be H\"older continuous in space. Starting from Sobolev solvability, our linear estimates yield a spatial $C_x^{2+\alpha}$ bound uniform in  $(\omega, t)$ . The same framework may also provide a regularity tool for the backward component of mean field games with common noise, where the unknown martingale integrand is coupled to the state dynamics through the shared noise.

The paper is organized as follows. In Section~\ref{sec:framework}, we state the main theorems and outline the proof. In Section~\ref{sec:windows}, we  document properties of the conditional energy norms, spatial interpolation inequalities, and global Sobolev results. In Section~\ref{sec:energy}, we prove the energy estimates under random spatial shifts and combine a local Sobolev estimate with heat-kernel bounds for the frozen model. In Section~\ref{sec:kernel}, we derive perturbation estimates for variable coefficients and lower-order terms. In Section~\ref{sec:closure}, we complete the proof of the main theorems. In Section~\ref{sec:coefficient-regularity}, we show sharpness of the regularity assumption on the martingale coefficients, and in  Section~\ref{sec:hjb}, we apply these estimates to semilinear stochastic HJB equations. The appendices contain auxiliary proofs.

\section{Main results and proof strategy}\label{sec:framework}
\subsection{Functional setting and main theorem}\label{subsec:notation}
Fix $0<\alpha<1$. Let $(\Omega,\mathcal F,\mathbb P)$ be a probability space supporting an $m$-dimensional Brownian motion $W:=(W^1,\ldots,W^m)$. Let $(\mathcal F_t)_{0\le t\le T}$ be its usual augmented natural filtration, with predictable $\sigma$-algebra $\mathcal P$ on $\Omega\times[0,T]$. For a stopping time $\tau\le T$, write $\E_\tau[\cdot]:=\E[\cdot\mid\mathcal F_\tau]$. The coefficients and forcing in~\eqref{eq:full},
\[
(a,b,c,\sigma,\gamma,f):\Omega\times[0,T]\times\mathbb T^n\to\R^{n\times n}\times\R^n\times\R\times\R^{n\times m}\times\R^m\times\R,
\]
are $\mathcal P\otimes\mathcal B(\mathbb T^n)$-measurable, with $a=a^\top$. The terminal datum $g:\Omega\times\mathbb T^n\to\R$ is $\mathcal F_T\otimes\mathcal B(\mathbb T^n)$-measurable.

For two vectors $x$ and $y$ in a Euclidean space, denote by $x\cdot y$ the Euclidean inner product. For matrices $M$ and $N$ of identical sizes,  denote by $M:N$  the Frobenius inner product: 
\[
M:N:=\operatorname{tr}(M^\top N)=\sum_{i,j}M_{ij}N_{ij}.
\]
With $(Dv)^{ik}:=\partial_i v^k$, we have
\[
a:D^2u=a^{ij}\partial_{ij}u,
\qquad
\sigma:Dv=\sigma^{ik}\partial_i v^k.
\]

For simplicity, we write $H_x^s:=H^s(\mathbb T^n)$ for the Sobolev spaces, with $s\in\R$, and $C_x^\beta:=C^\beta(\mathbb T^n)$ for the H\"older spaces, with noninteger $\beta>0$. We use the notation
\[
\begin{gathered}
L^2_{\omega,t}H_x^s:=L^2\bigl(\Omega;L^2(0,T;H_x^s)\bigr),
\qquad
L^2_\omega C_tH_x^s:=L^2\bigl(\Omega;C([0,T];H_x^s)\bigr),\\
L^\infty_\omega L_t^2L_x^2:=L^\infty\bigl(\Omega;L^2(0,T;L^2(\mathbb T^n))\bigr).
\end{gathered}
\]
In particular, when $s=0$, we write $L^2_{\omega,t,x}:=L^2_{\omega,t}H_x^0=L^2(\Omega\times(0,T)\times\mathbb T^n)$.

Uniform H\"older bounds for a time-dependent function $h:\Omega\times[0,T]\times\mathbb T^n\to\R$ and a terminal datum $g=g(\omega,x)$ are measured by
\[
\begin{aligned}
\|h\|_{L^\infty_{\omega,t}C_x^\beta}
&:=\esssup_{\omega,t}\|h(\omega,t,\cdot)\|_{C^\beta(\mathbb T^n)},\\
\|g\|_{L^\infty_\omega C_x^\beta}
&:=\esssup_\omega\|g(\omega,\cdot)\|_{C^\beta(\mathbb T^n)}.
\end{aligned}
\]

The Sobolev strong formulation below is the periodic version of the standard notion used by Du and Tang~\cite{Du2012Dirichlet} and Du, Qiu, and Tang~\cite{DuQiuTang2012}.

\begin{definition}[Sobolev strong and Schauder solutions]\label{def:energy}
An $\R\times\R^m$-valued pair $(u,v)$, with $u$ adapted and $v$ $\mathcal P\otimes\mathcal B(\mathbb T^n)$-measurable, is a Sobolev strong solution of~\eqref{eq:full} if
\[
u\in L^2_\omega C_tH_x^1\cap L^2_{\omega,t}H_x^2,
\qquad
v\in L^2_{\omega,t}H_x^1,
\]
and, almost surely for all $t\in[0,T]$, the following identity holds in $L^2(\mathbb T^n)$:
\begin{equation}\label{eq:energy-weak-form}
\begin{aligned}
u_t
={}&g+\int_t^T\bigl[a_s:D^2u_s+\sigma_s:Dv_s
+b_s\cdot Du_s+c_su_s+\gamma_s\cdot v_s+f_s\bigr]ds\\
&-\int_t^T v_s\cdot dW_s.
\end{aligned}
\end{equation}
All spatial derivatives are Sobolev weak derivatives. In~\eqref{eq:energy-weak-form}, the time integral is a Bochner integral in $L^2(\mathbb T^n)$, and the stochastic integral is an $L^2(\mathbb T^n)$-valued It\^o integral. We call a Sobolev strong solution a \emph{Schauder solution} if $u\in L^\infty_{\omega,t}C_x^{2+\alpha}$.
\end{definition}

For bounded coefficients, $f\in L^2_{\omega,t,x}$, and $g\in L^2_\omega H_x^1$, the regularity in Definition~\ref{def:energy} makes all terms in~\eqref{eq:energy-weak-form} well defined in $L^2(\mathbb T^n)$.

We assume $C^\alpha$ spatial regularity for $a,b,c$ and $C^{1+\alpha}$ spatial regularity for the martingale coefficients $\sigma,\gamma$, together with uniform super-parabolicity.
\begin{assumption}\label{ass:coeff}
There exist deterministic constants $K<\infty$ and $\lambda>0$ such that
\begin{equation}\label{eq:coeff-bounds}
\begin{aligned}
&\|a\|_{L^\infty_{\omega,t}C_x^\alpha}
+\|b\|_{L^\infty_{\omega,t}C_x^\alpha}
+\|c\|_{L^\infty_{\omega,t}C_x^\alpha}\\
&\qquad
+\|\sigma\|_{L^\infty_{\omega,t}C_x^{1+\alpha}}
+\|\gamma\|_{L^\infty_{\omega,t}C_x^{1+\alpha}}
\le K,
\end{aligned}
\end{equation}
and the effective diffusion satisfies
\begin{equation}\label{eq:ellipticity}
A:=a-\frac12\sigma\sigma^\top\ge\lambda I.
\end{equation}
\end{assumption}

Constants may depend on the fixed parameters $n,m,\alpha,T$ without explicit mention. Unless stated otherwise, $C=C(\lambda,K)>0$ is deterministic and may change from line to line. Subscripts on constants indicate additional dependence.

\begin{theorem}[Schauder solvability and estimate]\label{thm:linear}
Under Assumption~\ref{ass:coeff}, for $f\in L^\infty_{\omega,t}C_x^\alpha$ and $g\in L^\infty_\omega C_x^{2+\alpha}$, equation~\eqref{eq:full} has a unique Sobolev strong solution $(u,v)$. Moreover, $u$ belongs to $L^\infty_{\omega,t}C_x^{2+\alpha}$ and satisfies
\begin{equation}\label{eq:linear-estimate}
\|u\|_{L^\infty_{\omega,t}C_x^{2+\alpha}}
\le C\left(
\|f\|_{L^\infty_{\omega,t}C_x^\alpha}
+\|g\|_{L^\infty_\omega C_x^{2+\alpha}}
\right),
\end{equation}
where $C=C(\lambda,K)$. In particular, $(u,v)$ is a Schauder solution in the sense of Definition~\ref{def:energy}.
\end{theorem}

Estimate~\eqref{eq:linear-estimate} controls spatial H\"older regularity uniformly in  $(\omega, t)$ , as in the common-noise mean field game theory of Cardaliaguet, Delarue, Lasry, and Lions~\cite[Chapter~4]{CDLL2019}. The norm used here is stronger than that of Tang and Wei~\cite[Definition~2.1]{tang2016cauchy}, who measure $C^{2+\alpha}$ regularity using the $L^2(\Omega\times(0,T))$ norm before taking spatial suprema. Theorem~\ref{thm:joint-estimate} complements this estimate with a bound for $v$ in the conditional energy norm introduced in Section~\ref{subsec:proof-strategy}. The counterexamples in Section~\ref{sec:coefficient-regularity} show that the spatial H\"older exponent $1+\alpha$ required of both $\sigma$ and $\gamma$ is optimal for~\eqref{eq:linear-estimate}.

\begin{remark}[Time and spatial regularity]\label{rem:time-spatial-regularity}
By Definition~\ref{def:energy}, $u$ has continuous $H_x^1$-valued paths. The uniform spatial bound in Theorem~\ref{thm:linear} does not imply continuity in the $C_x^{2+\alpha}$ norm. Such continuity can fail even for the deterministic heat equation. In the stochastic setting, the martingale term also generally prevents classical differentiability in time.
\end{remark}

\begin{remark}[An unbounded martingale integrand]\label{rem:unbounded-integrand}
The assumptions of Theorem~\ref{thm:linear} do not imply that $v$ is essentially bounded, even when the forcing and terminal data are bounded and spatially constant.

For example, take $n=m=1$, $a=1$, $\sigma=b=c=\gamma=f=0$, and $g=\mathbf1_{\{W_T\ge0\}}$, the indicator of $\{W_T\ge0\}$. The martingale representation in Arai and Suzuki~\cite[Example~2.14]{AraiSuzuki2021}, specialized to standard Brownian motion, gives $u_t=\E_tg$ and
\[
v_t=\frac{1}{\sqrt{2\pi(T-t)}}
\exp\!\left(-\frac{W_t^2}{2(T-t)}\right),\qquad 0\le t<T.
\]
For $0<t<T$, this formula gives
\[
\esssup_\omega|v_t|=[2\pi(T-t)]^{-1/2},
\qquad\text{hence}\qquad v\notin L^\infty_{\omega,t,x}.
\]
However, we can still bound the conditional energy norm $\|v\|_{\mathcal E}$ defined in~\eqref{eq:old-2-7} in terms of the data; see Theorem~\ref{thm:joint-estimate}.
\end{remark}

\subsection{Littlewood--Paley decomposition and H\"older norms}\label{subsec:frequency-strategy}
We first introduce the periodic Fourier coefficients and the Littlewood--Paley blocks $\Delta_j$ used to decompose $u$ and $v$; see Littlewood and Paley~\cite{LittlewoodPaley1931} for the original theory and~\cite[Section~3]{GubinelliPerkowski2015} for the smooth periodic construction.

For $h\in\mathcal D'(\mathbb T^n)$ and $\phi\in C^\infty(\mathbb T^n)$, write $\langle h,\phi\rangle:=h(\phi)$. If $h\in L^1(\mathbb T^n)$, then
\[
\langle h,\phi\rangle=\int_{\mathbb T^n}h(x)\phi(x)\,dx.
\]
We define the Fourier coefficients by
\[
\widehat h(\ell):=\frac{1}{(2\pi)^n}
\langle h,e^{-i\ell\cdot x}\rangle,
\qquad \ell\in\mathbb Z^n,
\]
and define the Fourier support by
\[
\operatorname{supp}\widehat h:=\{\ell\in\mathbb Z^n:\widehat h(\ell)\ne0\}.
\]

Choose a smooth radial $\chi\in C_c^\infty(\R^n)$ equal to one on $\{|\xi|\le1\}$ and zero on $\{|\xi|\ge2\}$, and put
$r_j:=2^{-j}$, $\varphi(\xi):=\chi(\xi)-\chi(2\xi)$, $m_0:=\chi$, and $m_j(\xi):=\varphi(r_j\xi)$ for $j\ge1$. Define
\[
\Delta_jh(x):=\sum_{\ell\in\mathbb Z^n}
m_j(\ell)\widehat h(\ell)e^{i\ell\cdot x}.
\]
Each block is a finite Fourier sum, and $h=\sum_{j\ge0}\Delta_jh$ in $\mathcal D'(\mathbb T^n)$. The Fourier support of $\Delta_0h$ lies in $\{\ell\in\mathbb Z^n:|\ell|\le2\}$. For $j\ge1$, define the frequency annulus
\[
\mathcal C_j:=\{\xi\in\R^n:2^{j-1}\le|\xi|\le2^{j+1}\}.
\]
Then $m_j$ is supported in $\mathcal C_j$, and the Fourier support of $\Delta_jh$ lies in $\mathcal C_j\cap\mathbb Z^n$.

For each $j\ge1$, we fix an auxiliary multiplier that equals one on the annulus $\mathcal C_j$:
\begin{equation}\label{eq:reproducing-cutoff}
\widetilde m_j(\xi):=\chi(r_j\xi/2)-\chi(4r_j\xi),\qquad j\ge1.
\end{equation}
It is supported in $\{2^{j-2}\le|\xi|\le2^{j+2}\}$. In particular, $\widetilde m_jm_j=m_j$, so applying this multiplier leaves the block $\Delta_jh$ unchanged. These frequency cutoffs are fixed throughout the paper, and constants may depend on their choice without explicit mention.

Spatial translation and differentiation act on Fourier coefficients by
\[
\widehat{h(\cdot+z)}(\ell)=e^{i\ell\cdot z}\widehat h(\ell),
\qquad
\widehat{\partial_i h}(\ell)=i\ell_i\widehat h(\ell).
\]
Both operations therefore preserve the support bounds $|\ell|\le2$ for $\Delta_0h$ and $\ell\in\mathcal C_j$ for $\Delta_jh$, $j\ge1$.

For $h\in L^2(\mathbb T^n)$ and each fixed $j$, Cauchy--Schwarz and Parseval give
\begin{equation}\label{eq:finite-fourier-block}
\|\Delta_jh\|_\infty
\le\sum_{\ell\in\mathbb Z^n}|m_j(\ell)\widehat h(\ell)|
\le C_j\|h\|_{L_x^2}.
\end{equation}

For a random field $h:\Omega\times[0,T]\times\mathbb T^n\to\R$, Fourier support conditions are understood to hold for almost every $(\omega,t)$. Write $u_j:=\Delta_ju$ and $v_j:=\Delta_jv$.
For the components of $v_j$, we write
\[
v_j^k=\Delta_j(v^k),\qquad j\ge0,\quad 1\le k\le m,
\]
where $j$ is the frequency index and $k$ indexes the Brownian components.

We recall the standard H\"older--Besov norm equivalence; see~\cite[Lemma~9]{GubinelliPerkowski2015}, with the higher-order cases obtained by applying the characterization to spatial derivatives.
\begin{lemma}[H\"older--Besov norm equivalence]\label{lem:frequency}
For every noninteger $\beta>0$, $C^\beta(\mathbb T^n)=B^\beta_{\infty,\infty}(\mathbb T^n)$ with equivalent norms:
\begin{equation}\label{eq:old-2-15}
\|h\|_{C^\beta}
\asymp\|h\|_{B^\beta_{\infty,\infty}}
:=\sup_{j\ge0}r_j^{-\beta}\|\Delta_jh\|_\infty,
\end{equation}
where $F\asymp G$ means $C^{-1}G\le F\le CG$ for some $C=C(\beta)\ge1$.
\end{lemma}

The convolution kernel of $\Delta_j$ is the trigonometric polynomial
\[
\psi_j(y):=\frac{1}{(2\pi)^n}
\sum_{\ell\in\mathbb Z^n}m_j(\ell)e^{i\ell\cdot y},
\qquad
\Delta_jh=\psi_j*h,
\]
where convolution uses Lebesgue measure on $\mathbb T^n$. For $y\in\mathbb T^n$, define
\begin{equation}\label{eq:old-2-2}
\rho(y):=\operatorname{dist}_{\mathbb T^n}(y,0)
=\min_{q\in\mathbb Z^n}|y+2\pi q|.
\end{equation}

The following lemma gives uniform $L^1$ bounds and spatial moment estimates for the block kernels $\psi_j$. A standard proof is included in Appendix~\ref{app:block-kernels} for completeness.
\begin{lemma}\label{lem:block-kernels}
For every $j\ge0$, the kernel $\psi_j$ satisfies
\[
\int_{\mathbb T^n}\psi_j(y)\,dy=
\begin{cases}
1,&j=0,\\
0,&j\ge1
\end{cases}
\]
and, for every $0\le\nu\le2$,
\begin{equation}\label{eq:block-kernel-bounds}
\int_{\mathbb T^n}\rho(y)^\nu|\psi_j(y)|\,dy\le C_\nu r_j^\nu,
\end{equation}
where $C_\nu$ is independent of $j$.
\end{lemma}

\subsection{Conditional local energy and proof strategy}\label{subsec:proof-strategy}
We now define the conditional energy norm used to control the martingale integrand at each parabolic scale.

For spatial localization, introduce the weight $w_r$ and its periodization $\kappa_r$, for $0<r\le1$:
\begin{equation}\label{eq:old-2-3}
w_r(z):=c_nr^{-n}\left(1+\frac{|z|^2}{r^2}\right)^{-(n+1)/2},
\qquad z\in\R^n,
\end{equation}
and
\begin{equation}\label{eq:periodic-window-weight}
\kappa_r(y):=\sum_{\ell\in\mathbb Z^n}w_r(y+2\pi\ell),
\qquad y\in\mathbb T^n.
\end{equation}
The constant $c_n$ normalizes both weights to have integral one. Write $\kappa_{r,x}(y):=\kappa_r(y-x)$ for the window centered at $x$.

For $\kappa=\kappa_r$ or $\kappa_{r,x}$, set
\begin{equation}\label{eq:old-2-4}
\|h\|_\kappa^2:=\int_{\mathbb T^n}\kappa(y)|h(y)|^2\,dy,
\end{equation}
and
\begin{equation}\label{eq:weighted-inner-product}
\langle h,q\rangle_\kappa:=\int_{\mathbb T^n}\kappa(y)\langle h(y),q(y)\rangle\,dy.
\end{equation}
The pointwise product is Euclidean or Frobenius, with scalar--vector pairings interpreted componentwise. In particular, $\|h\|_{\kappa_{r,x}}=\|h(x+\cdot)\|_{\kappa_r}$.

The following lemma gives weighted pointwise bounds for functions with bounded Fourier support and their first two spatial derivatives, together with the standard $L^\infty$ Bernstein inequality~\cite[Lemma~7]{GubinelliPerkowski2015}. Its proof is given in Appendix~\ref{app:frequency}.
\begin{lemma}[Bernstein inequalities]\label{lem:bernstein}
Let $k\ge0$ and let $h_k\in L^2(\mathbb T^n)$ satisfy
\[
\operatorname{supp}\widehat h_k\subset\{\xi\in\mathbb Z^n:|\xi|\le C_*r_k^{-1}\}
\]
for a fixed $C_*\ge2$. Then
\begin{enumerate}
\renewcommand{\labelenumi}{(\roman{enumi})}
\item For $\ell=0,1,2$,
\begin{equation}\label{eq:old-2-13}
|D^\ell h_k|^2\le Cr_k^{-2\ell}\kappa_{r_k}*|h_k|^2.
\end{equation}
\item For $\ell=1,2$,
\begin{equation}\label{eq:bernstein-infty}
\|D^\ell h_k\|_\infty\le Cr_k^{-\ell}\|h_k\|_\infty.
\end{equation}
\end{enumerate}
Here $\kappa_{r_k}$ is the periodic weight defined in~\eqref{eq:periodic-window-weight}, and $C=C(C_*)$ is independent of $k$.
\end{lemma}

Throughout the conditional energy estimates, $h$ is a $\mathcal P\otimes\mathcal B(\mathbb T^n)$-measurable random field whose weighted square integrals over finite time windows are almost surely finite. Estimates are understood whenever their right-hand sides are finite. Let $\tau$ range over stopping times taking values in $[0,T]$, and keep the spatial center $x$ fixed during the time integral. Using $\kappa_{r,x}$ on the interval $[\tau,(\tau+r^2)\wedge T]$, we define the conditional local energy by
\begin{equation}\label{eq:old-2-6}
\mathcal E_r(h):=\sup_{\tau\le T}
\left\|\left[
\E_\tau\int_\tau^{(\tau+r^2)\wedge T}
\|h_s\|_{\kappa_{r,x}}^2\,ds
\right]^{1/2}\right\|_{L^\infty_{\omega,x}}.
\end{equation}
The energy is \emph{local} in both space and time: the weight $\kappa_{r,x}$ emphasizes points within distance of order $r$ from $x$, with polynomially decaying tails, and the time integral is taken over an interval of length at most $r^2$. The local Sobolev estimate in Section~\ref{subsec:weighted-energy} uses $\mathcal E_r$ to control conditional weighted $L^2$ norms of spatial derivatives and martingale integrands.

The weighted energy identity couples the size of $u_j$ to the conditional energy of $v_j$. To match the frequency bound $\|u_j\|_{L^\infty_{\omega,t,x}}\le Cr_j^{2+\alpha}\|u\|_{L^\infty_{\omega,t}C_x^{2+\alpha}}$, we define
\begin{equation}\label{eq:old-2-7}
\|v\|_{\mathcal E}:=
\sup_{j\ge0}r_j^{-2-\alpha}\mathcal E_{r_j}(v_j),
\end{equation}
and abbreviate
\begin{equation}\label{eq:old-2-8}
e_j:=\mathcal E_{r_j}(v_j).
\end{equation}

Since $\mathcal E_{r_j}(v_j)\le r_j\|v_j\|_{L^\infty_{\omega,t,x}}$, Lemma~\ref{lem:frequency} gives
\[
\|v\|_{\mathcal E}
\le C\|v\|_{L^\infty_{\omega,t}C_x^{1+\alpha}}.
\]
The following theorem complements Theorem~\ref{thm:linear} by bounding $\|v\|_{\mathcal E}$ in terms of the data, also covering the unbounded integrand in Remark~\ref{rem:unbounded-integrand}.
\begin{theorem}\label{thm:joint-estimate}
Under the hypotheses of Theorem~\ref{thm:linear}, the unique Sobolev strong solution $(u,v)$ satisfies
\begin{equation}\label{eq:joint-estimate}
\|v\|_{\mathcal E}
\le C\left(\|f\|_{L^\infty_{\omega,t}C_x^\alpha}
+\|g\|_{L^\infty_\omega C_x^{2+\alpha}}\right),
\end{equation}
where $C=C(\lambda,K)$.
\end{theorem}

The conditional energy estimate also yields spatial H\"older regularity for the martingale integrand $v$ in a BMO norm. For a predictable $\R^m$-valued process $q$ and $\delta>0$, define
\[
\|q\|_{\mathbb H^2_{\mathrm{BMO}}(\delta)}
:=
\sup_{\tau\le T}
\left\|
\left(
\E_\tau\int_\tau^{(\tau+\delta)\wedge T}|q_s|^2\,ds
\right)^{1/2}
\right\|_{L^\infty_\omega}.
\]
Write $\mathbb H^2_{\mathrm{BMO}}:=\mathbb H^2_{\mathrm{BMO}}(T)$, with norm given by the $\mathrm{BMO}_2$ norm of $\int_0^\cdot q_s\cdot dW_s$, and denote the corresponding spatial H\"older space by $C_x^{1+\alpha}(\mathbb H^2_{\mathrm{BMO}})$.

\begin{corollary}\label{cor:bmo}
Under the hypotheses of Theorem~\ref{thm:linear}, the martingale integrand $v$ admits a version in $C_x^{1+\alpha}(\mathbb H^2_{\mathrm{BMO}})$ satisfying
\begin{equation}\label{eq:bmo-estimate}
\|v\|_{C_x^{1+\alpha}(\mathbb H^2_{\mathrm{BMO}})}
\le C\left(\|f\|_{L^\infty_{\omega,t}C_x^\alpha}
+\|g\|_{L^\infty_\omega C_x^{2+\alpha}}\right),
\end{equation}
where $C=C(\lambda,K)$.
\end{corollary}

\begin{proof}
Fix $j\ge0$, put $r=r_j$, and set
\[
B_j:=\sup_{x\in\mathbb T^n}
\|v_j(\cdot,x)\|_{\mathbb H^2_{\mathrm{BMO}}(r^2)}.
\]
For a stopping time $\tau$, let $\zeta=(\tau+r^2)\wedge T$. Conditional Fubini gives
\[
\E_\tau\int_\tau^\zeta\|v_j(s)\|_{\kappa_{r,x}}^2\,ds
=\int_{\mathbb T^n}\kappa_r(y-x)
\E_\tau\int_\tau^\zeta|v_j(s,y)|^2\,ds\,dy
\le B_j^2.
\]
Thus $\mathcal E_r(v_j)\le B_j$. Conversely, the local Bernstein bound~\eqref{eq:old-2-13} with $\ell=0$ gives
\[
|v_j(s,x)|^2
\le C(\kappa_r*|v_j(s)|^2)(x)
=C\|v_j(s)\|_{\kappa_{r,x}}^2.
\]
Integrating over $[\tau,\zeta]$ and conditioning at $\tau$ yields $B_j\le C\mathcal E_r(v_j)$, uniformly in $j$. This proves the equivalent characterization
\[
\|v\|_{\mathcal E}
\asymp\sup_{j\ge0}r_j^{-2-\alpha}B_j.
\]

Applying the local bound on successive intervals starting at $(\tau+\nu r^2)\wedge T$, we obtain
\[
\sup_{x\in\mathbb T^n}
\|v_j(\cdot,x)\|_{\mathbb H^2_{\mathrm{BMO}}}
\le \left\lceil T/r^2\right\rceil^{1/2}B_j
\le C_T r^{-1}\mathcal E_r(v_j).
\]
Summing the frequency blocks in $\mathbb H^2_{\mathrm{BMO}}$ and applying the Banach-valued version of Lemma~\ref{lem:frequency} give
\[
\|v\|_{C_x^{1+\alpha}(\mathbb H^2_{\mathrm{BMO}})}
\le C_T\|v\|_{\mathcal E}.
\]
Theorem~\ref{thm:joint-estimate} now yields~\eqref{eq:bmo-estimate}.
\end{proof}

Theorems~\ref{thm:linear} and~\ref{thm:joint-estimate} are proved together in Section~\ref{sec:closure}. The argument combines two complementary bounds through $\|v\|_{\mathcal E}$. We first work with spatially smooth approximations, for which both $\|u\|_{L^\infty_{\omega,t}C_x^{2+\alpha}}$ and $\|v\|_{\mathcal E}$ are finite. A weighted It\^o identity for each pair $(u_j,v_j)$, together with super-parabolicity, gives
\[
\|v\|_{\mathcal E}\le C_E\bigl(\|u\|_{L^\infty_{\omega,t}C_x^{2+\alpha}}+\|f\|_{L^\infty_{\omega,t}C_x^\alpha}+\|g\|_{L^\infty_\omega C_x^{2+\alpha}}\bigr).
\]

To bound $\|u\|_{L^\infty_{\omega,t}C_x^{2+\alpha}}$ in terms of this energy, we freeze $a,\sigma$ at a spatial point and use a stochastic integrating factor in Fourier space to cancel the derivative coupling. The resulting representation controls the martingale terms by $\|v\|_{\mathcal E}$. Coefficient regularity makes the high-frequency perturbations small, while comparison and conditional Sobolev estimates bound the low frequencies in terms of the data. For every $\varepsilon>0$, this yields
\[
\|u\|_{L^\infty_{\omega,t}C_x^{2+\alpha}}\le C_\varepsilon\bigl(\|f\|_{L^\infty_{\omega,t}C_x^\alpha}+\|g\|_{L^\infty_\omega C_x^{2+\alpha}}\bigr)+\varepsilon\|v\|_{\mathcal E}.
\]
Choosing $\varepsilon C_E\le1/2$ and combining these inequalities bounds both norms by the data. Sobolev convergence and lower semicontinuity transfer these uniform bounds from the smooth approximations to the original solution.

\section{Auxiliary lemmas}\label{sec:windows}\label{subsec:auxiliary-lemmas}
We collect conditional energy estimates, spatial interpolation inequalities, and global Sobolev results needed for the argument.

\subsection{Conditional energy estimates and spatial interpolation}\label{subsubsec:local-energy-estimates}
The next lemma bounds the conditional energy of spatial derivatives and compares different localization scales for the same functions $h_k\in L^2_{\omega,t,x}$ with bounded Fourier support. The proof is given in Appendix~\ref{app:weights}.

\begin{lemma}\label{lem:window-scales}
Let $k\ge0$ and let $h_k\in L^2_{\omega,t,x}$ be a predictable random field satisfying
\[
\operatorname{supp}\widehat h_k\subset\{\xi\in\mathbb Z^n:|\xi|\le C_*r_k^{-1}\}
\]
for almost every $(\omega,t)$, where $C_*\ge2$ is fixed. Then, for every $j\ge0$,
\begin{equation}\label{eq:old-2-18}
\mathcal E_{r_j}(D^\ell h_k)
\le
Cr_k^{-\ell}\mathcal E_{r_j}(h_k),
\qquad \ell=1,2,
\end{equation}
and
\begin{equation}\label{eq:old-2-19}
\mathcal E_{r_j}(h_k)
\le
C\max\left\{1,\frac{r_j}{r_k}\right\}
\mathcal E_{r_k}(h_k).
\end{equation}
Here $C=C(C_*)$ is independent of $j,k$.
\end{lemma}

We use the following spatial interpolation inequalities; see Gilbarg and Trudinger~\cite{GilbargTrudinger2001} for classical H\"older interpolation. A dyadic proof of both forms used here is given in Appendix~\ref{app:frequency}.
\begin{lemma}[Spatial interpolation]\label{lem:holder-interpolation}
For every $\delta>0$, there exists $C_\delta$ such that, for all $h\in C^{2+\alpha}(\mathbb T^n)$,
\begin{equation}\label{eq:old-2-16}
\|D^2h\|_\infty+\|h\|_{C^{1+\alpha}}
\le\delta\|h\|_{C^{2+\alpha}}+C_\delta\|h\|_\infty.
\end{equation}
Similarly,
\begin{equation}\label{eq:holder-sobolev-interpolation}
\|h\|_{C^{1+\alpha}}
\le\delta\|h\|_{C^{2+\alpha}}+C_\delta\|h\|_{H_x^1}.
\end{equation}
\end{lemma}

\subsection{Basic results for BSPDEs}\label{subsec:sobolev}
Sobolev well-posedness, stability, and comparison for BSPDEs on $\R^n$ are established in~\cite[Theorems~5.5, 6.1, and~6.3]{DuQiuTang2012} for $p=2$, and the same arguments apply on $\mathbb T^n$. We state the periodic solvability, stability, and comparison results below. Appendix~\ref{app:sobolev} gives the proofs of solvability, the conditional estimates, and comparison.

Throughout this subsection, we impose the following weaker assumption for coefficient regularity.
\begin{assumption}\label{ass:sobolev-coeff}
The principal coefficients $a,\sigma$ satisfy the bounds and super-parabolicity condition in Assumption~\ref{ass:coeff}. The lower-order coefficients $b,c,\gamma$ satisfy
\[
\|b\|_\infty+\|c\|_\infty+\|\gamma\|_\infty\le K.
\]
\end{assumption}

The following lemma gives Sobolev solvability and conditional estimates in terms of the data, independently of Schauder regularity.
\begin{lemma}[Periodic Sobolev theory]\label{lem:periodic-existence}
Under Assumption~\ref{ass:sobolev-coeff}, for $f\in L^2_{\omega,t,x}$ and $g\in L^2_\omega H_x^1$, equation~\eqref{eq:full} has a unique Sobolev strong solution $(u,v)$. For every stopping time $\tau\le T$, this solution satisfies
\begin{equation}\label{eq:conditional-sobolev}
\E_\tau\int_\tau^T\bigl(\|u_s\|_{H_x^2}^2+\|v_s\|_{H_x^1}^2\bigr)\,ds
\le C\E_\tau\left[\|g\|_{H_x^1}^2+\int_\tau^T\|f_s\|_{L_x^2}^2\,ds\right].
\end{equation}
In particular, if $f\in L^\infty_\omega L_t^2L_x^2$ and $g\in L^\infty_\omega H_x^1$, then
\begin{equation}\label{eq:old-2-23}
\sup_\tau\left\|
\left(\E_\tau\int_\tau^T\|v_s\|_{H_x^1}^2\,ds\right)^{1/2}
\right\|_{L^\infty_\omega}
\le C\left(\|f\|_{L^\infty_\omega L_t^2L_x^2}
+\|g\|_{L^\infty_\omega H_x^1}\right),
\end{equation}
where the supremum is over all stopping times $\tau\le T$, and $C=C(\lambda,K)$.
\end{lemma}

Applying Lemma~\ref{lem:periodic-existence} to the difference of two equations yields stability with respect to the coefficients and data.
\begin{corollary}\label{cor:sobolev-stability}
Under the hypotheses of Lemma~\ref{lem:periodic-existence}, suppose that $(a^\nu,\sigma^\nu,b^\nu,c^\nu,\gamma^\nu)$ satisfy the same coefficient assumptions with common bounds and converge to $(a,\sigma,b,c,\gamma)$ in $L^\infty_{\omega,t,x}$. If $f^\nu\to f$ in $L^2_{\omega,t,x}$ and $g^\nu\to g$ in $L^2_\omega H_x^1$, then the corresponding Sobolev strong solutions satisfy
\begin{equation}\label{eq:sobolev-stability}
u^\nu\to u\quad\mathrm{in}\quad L^2_{\omega,t}H_x^2,
\qquad v^\nu\to v\quad\mathrm{in}\quad L^2_{\omega,t}H_x^1.
\end{equation}
Here $(u,v)$ is the solution for the limiting coefficients and data.
\end{corollary}

\begin{proof}
We use the difference-equation argument of~\cite[Theorem~6.1]{DuQiuTang2012}. The pair $(u^\nu-u,v^\nu-v)$ solves the equation with coefficients $(a^\nu,\sigma^\nu,b^\nu,c^\nu,\gamma^\nu)$, terminal datum $g^\nu-g$, and forcing
\[
\begin{aligned}
F^\nu:={}&(a^\nu-a):D^2u+(\sigma^\nu-\sigma):Dv\\
&+(b^\nu-b)\cdot Du+(c^\nu-c)u+(\gamma^\nu-\gamma)\cdot v+(f^\nu-f).
\end{aligned}
\]
Since the coefficient differences converge to zero in $L^\infty_{\omega,t,x}$ and $D^2u,Dv,Du,u,v\in L^2_{\omega,t,x}$, each product in $F^\nu$ converges to zero in $L^2_{\omega,t,x}$. Together with $f^\nu\to f$, this gives $F^\nu\to0$ in the same space. Applying~\eqref{eq:conditional-sobolev} at $\tau=0$ and taking expectation gives, with a constant independent of $\nu$,
\[
\|u^\nu-u\|_{L^2_{\omega,t}H_x^2}
+\|v^\nu-v\|_{L^2_{\omega,t}H_x^1}
\le C\bigl(\|F^\nu\|_{L^2_{\omega,t,x}}
+\|g^\nu-g\|_{L^2_\omega H_x^1}\bigr)\longrightarrow0.
\]
This proves~\eqref{eq:sobolev-stability}.
\end{proof}

We also use the following comparison principle.
\begin{lemma}\label{lem:periodic-sobolev}
Under Assumption~\ref{ass:sobolev-coeff}, let $(u,v)$ and $(\bar u,\bar v)$ be Sobolev strong solutions of~\eqref{eq:full} with the same coefficients and data $(f,g)$ and $(\bar f,\bar g)$, respectively. If
\[
f\le\bar f\quad(\mathbb P\otimes dt\otimes dx)\mbox{-a.e.},
\qquad g\le\bar g\quad(\mathbb P\otimes dx)\mbox{-a.e.}
\]
then, $\mathbb P$-a.s., for every $t\in[0,T]$,
\[
u_t(x)\le\bar u_t(x)\qquad\text{for a.e. }x\in\mathbb T^n.
\]
\end{lemma}

The comparison principle yields a uniform bound for $u$, and hence bounds for its low-frequency blocks.
\begin{corollary}\label{cor:uniform-u}
Under the hypotheses of Lemma~\ref{lem:periodic-existence}, suppose in addition that $f\in L^\infty_{\omega,t,x}$ and $g\in L^\infty_{\omega,x}$. Then the solution $(u,v)$ satisfies
\begin{equation}\label{eq:old-2-22}
\|u\|_{L^\infty_{\omega,t,x}}
\le C\left(\|f\|_{L^\infty_{\omega,t,x}}
+\|g\|_{L^\infty_{\omega,x}}\right).
\end{equation}
Moreover, for every integer $N\ge1$,
\begin{equation}\label{eq:old-2-25}
\sup_{0\le j<N}r_j^{-2-\alpha}\|\Delta_ju\|_{L^\infty_{\omega,t,x}}
\le C_N\left(\|f\|_{L^\infty_{\omega,t,x}}
+\|g\|_{L^\infty_{\omega,x}}\right).
\end{equation}
Here $C=C(\lambda,K)$ and $C_N=C_N(\lambda,K)$.
\end{corollary}

\begin{proof}[Proof of Corollary~\ref{cor:uniform-u}]
Put $F_0:=\|f\|_{L^\infty_{\omega,t,x}}$ and $G_0:=\|g\|_{L^\infty_{\omega,x}}$, and define
\[
q(t):=e^{K(T-t)}G_0+F_0\int_t^T e^{K(s-t)}\,ds.
\]
Then the spatially constant pair $(q,0)$ solves the BSPDE (degenerate to PDE)
\[
-q'(t)=Kq(t)+F_0,\qquad q(T)=G_0.
\]
Since $q\ge0$ and $c\le K$, this pair also solves~\eqref{eq:full} with forcing $(K-c)q+F_0\ge f$ and terminal datum $G_0\ge g$. The comparison principle in Lemma~\ref{lem:periodic-sobolev} gives $u_t\le q(t)\le q(0)$. Applying the same argument to $(-u,-v)$ with data $(-f,-g)$ gives $-u_t\le q(t)$. Hence
\[
|u_t|\le q(0)\le C(F_0+G_0),
\]
proving~\eqref{eq:old-2-22}.

The uniform $L^1$ bounds for $\psi_j$ in Lemma~\ref{lem:block-kernels} give
\[
\|u_j\|_{L^\infty_{\omega,t,x}}\le\|\psi_j\|_{L^1}\|u\|_{L^\infty_{\omega,t,x}}
\le C\|u\|_{L^\infty_{\omega,t,x}}.
\]
Multiplying by $r_j^{-2-\alpha}$ and taking the maximum over $j<N$ yields~\eqref{eq:old-2-25} by~\eqref{eq:old-2-22}.
\end{proof}

The conditional Sobolev estimate gives the corresponding low-frequency bound for $v$.
\begin{lemma}\label{lem:low}
Under Assumption~\ref{ass:sobolev-coeff}, let $f\in L^\infty_\omega L_t^2L_x^2$ and $g\in L^\infty_\omega H_x^1$, and let $(u,v)$ be the Sobolev strong solution given by Lemma~\ref{lem:periodic-existence}. Then, for every integer $N\ge0$,
\begin{equation}\label{eq:old-2-24}
\sum_{k=0}^{N}r_k^{-1}\mathcal E_{r_k}(\Delta_kv)
\le C_N\left(\|f\|_{L^\infty_\omega L_t^2L_x^2}
+\|g\|_{L^\infty_\omega H_x^1}\right),
\end{equation}
where $C_N=C_N(\lambda,K)$.
\end{lemma}

\begin{proof}[Proof of Lemma~\ref{lem:low}]
For each fixed $k\ge0$, estimate~\eqref{eq:finite-fourier-block} gives
\[
\|v_k(s)\|_{\kappa_{r_k,x}}^2
\le\|v_k(s)\|_\infty^2
\le C_k\|v_s\|_{L_x^2}^2.
\]
Thus, by~\eqref{eq:conditional-sobolev}, for every stopping time $\tau\le T$ and center $x\in\mathbb T^n$,
\[
\begin{aligned}
\E_\tau\int_\tau^{(\tau+r_k^2)\wedge T}
\|v_k(s)\|_{\kappa_{r_k,x}}^2\,ds
&\le C_k\E_\tau\int_\tau^T\|v_s\|_{H_x^1}^2\,ds\\
&\le C_k\left(\|f\|_{L^\infty_\omega L_t^2L_x^2}^2
+\|g\|_{L^\infty_\omega H_x^1}^2\right).
\end{aligned}
\]
The definition of $\mathcal E_{r_k}$ and summation over $k\le N$ give~\eqref{eq:old-2-24}.
\end{proof}

\section{Local Sobolev and uniform estimates for the frozen model}\label{sec:energy}
We first consider a BSPDE with coefficients independent of space:
\begin{equation}\label{eq:frozen-model}
-du_j=[a_0:D^2u_j+\sigma_0:Dv_j+F_j]dt-v_j\cdot dW_t,
\qquad u_j(T)=g_j,
\end{equation}
where $a_0,\sigma_0\in L^\infty_{\omega,t}$ are predictable, $a_0=a_0^\top$, and
\[
A_0:=a_0-\tfrac12\sigma_0\sigma_0^\top\ge\lambda I,
\qquad |a_0|+|\sigma_0|\le K.
\]
Fix $j\ge1$ and put $r=r_j$. The forcing $F_j\in L^2_{\omega,t,x}$ and terminal datum $g_j\in L^2_\omega H_x^1$ have Fourier support in $\mathcal C_j\cap\mathbb Z^n$. The same is true of the Sobolev strong solution $(u_j,v_j)$, since the coefficients commute with spatial Fourier projections. The symbol $F_j$ denotes the full forcing of this model.

We first bound $\mathcal E_{r_j}(v_j)$ in terms of $\|u_j\|_{L^\infty_{\omega,t,x}}$ and $r_j\mathcal E_{r_j}(F_j)$ by a weighted It\^o identity. We then construct a conditional representation for $u_j$ and bound its forcing term by $r_j\mathcal E_{r_j}(F_j)$.

\subsection{Weighted local Sobolev estimate}\label{subsec:weighted-energy}
The coupling $\sigma_0:Dv_j$ is handled by combining its contribution with the quadratic variation in It\^o's formula. Super-parabolicity then gives a conditional local Sobolev estimate controlling the weighted square integrals of $Du_j$ and $v_j$.
\begin{lemma}\label{lem:weighted-window-energy}
For a Sobolev strong solution of~\eqref{eq:frozen-model} with $u_j\in L^\infty_{\omega,t,x}$,
\begin{equation}\label{eq:model-energy}
\mathcal E_r(Du_j)^2+\mathcal E_r(v_j)^2
\le C\left(\|u_j\|_{L^\infty_{\omega,t,x}}^2
+r\|u_j\|_{L^\infty_{\omega,t,x}}\mathcal E_r(F_j)\right).
\end{equation}
\end{lemma}
\begin{proof}
Fix a stopping time $\tau$ and let $\zeta=(\tau+r^2)\wedge T$. Spatial translation allows us to center the weight at $0$. Write $\kappa=\kappa_r$, with $r=r_j$, for the weight localizing the energy identity at the spatial scale of $u_j$ and $v_j$. Define $S_s=\|u_j(s)\|_\kappa^2$ and $R_s=v_j(s)+\sigma_0(s)^\top Du_j(s)$. It\^o's formula gives
\[
dS_s=\left[-2\langle u_j,a_0:D^2u_j+\sigma_0:Dv_j+F_j\rangle_\kappa+\|v_j\|_\kappa^2\right]ds+2\langle u_j,v_j\rangle_\kappa\cdot dW_s.
\]
Spatial integration by parts and the identity
\[
2\int\kappa a_0Du_j\cdot Du_j+2\int\kappa\sigma_0^\top Du_j\cdot v_j+\|v_j\|_\kappa^2
=2\int\kappa A_0Du_j\cdot Du_j+\|R_s\|_\kappa^2
\]
yield
\begin{equation}\label{eq:old-3-9}
dS_s=\left[2\int\kappa A_0Du_j\cdot Du_j+\|R_s\|_\kappa^2
-2\langle u_j,F_j\rangle_\kappa+\mathcal R_s\right]ds
+2\langle u_j,v_j\rangle_\kappa\cdot dW_s,
\end{equation}
where
\[
\mathcal R_s=2\int u_j[(a_0-\sigma_0\sigma_0^\top)D\kappa]\cdot Du_j
+2\int u_j(\sigma_0^\top D\kappa)\cdot R_s.
\]
With weight $1$, the remainder $\mathcal R_s$ vanishes and one recovers the usual spatial $L^2$ energy identity.

Differentiating the Euclidean weight in~\eqref{eq:old-2-3} gives
\[
Dw_r(z)=-(n+1)\frac{z}{r^2+|z|^2}w_r(z),
\qquad |Dw_r(z)|\le Cr^{-1}w_r(z).
\]
Summing over its lattice translates yields $|D\kappa|\le Cr^{-1}\kappa$. Young's inequality therefore gives
\[
|\mathcal R_s|\le\delta(\|Du_j(s)\|_\kappa^2+\|R_s\|_\kappa^2)
+C_\delta r^{-2}S_s.
\]
Integrate~\eqref{eq:old-3-9} from $\tau$ to $\zeta$, take conditional expectation, and choose $\delta$ small. A standard localization argument shows that the stochastic integral has zero conditional expectation. Using $A_0\ge\lambda I$ and $|v_j|^2\le2|R_s|^2+2K^2|Du_j|^2$, we obtain
\[
\begin{aligned}
&\E_\tau\int_\tau^\zeta(\|Du_j(s)\|_\kappa^2+\|v_j(s)\|_\kappa^2)ds\\
&\qquad\le C\E_\tau\left[S_\zeta+r^{-2}\int_\tau^\zeta S_sds
+\int_\tau^\zeta|\langle u_j,F_j\rangle_\kappa|ds\right].
\end{aligned}
\]
The finite Fourier sum for $u_j$ has continuous spatially uniform paths, so its essential bound also holds at stopping times. Hence the first two terms are bounded by $C\|u_j\|_{L^\infty_{\omega,t,x}}^2$. Conditional Cauchy--Schwarz gives
\[
\E_\tau\int_\tau^\zeta|\langle u_j,F_j\rangle_\kappa|ds
\le r\|u_j\|_{L^\infty_{\omega,t,x}}\mathcal E_r(F_j).
\]
The same bounds hold at every spatial center. Taking the norms and suprema in the definition of $\mathcal E_r$ proves~\eqref{eq:model-energy}.
\end{proof}

\subsection{Fourier representation and kernel estimates}\label{subsec:kernel-representation}
We now estimate $\|u_j\|_{L^\infty_{\omega,t,x}}$ in~\eqref{eq:model-energy}. A stochastic integrating factor in Fourier space cancels the contribution of $\sigma_0:Dv_j$ and gives a representation in terms of $g_j$ and $F_j$.

Taking spatial Fourier coefficients in~\eqref{eq:frozen-model} gives, for each $\xi\in\mathcal C_j\cap\mathbb Z^n$,
\[
\left\{
\begin{aligned}
-d\widehat u_j(s,\xi)
&=\bigl[-(\xi^\top a_0(s)\xi)\widehat u_j(s,\xi)
+i(\sigma_0(s)^\top\xi)\cdot\widehat v_j(s,\xi)\\
&\qquad+\widehat F_j(s,\xi)\bigr]ds-\widehat v_j(s,\xi)\cdot dW_s,\\
\widehat u_j(T,\xi)&=\widehat g_j(\xi).
\end{aligned}
\right.
\]

Fix $t<T$ and define
\[
Q_{t,s}:=\int_t^s A_0(\ell)\,d\ell,\qquad
I_{t,s}(\xi):=e^{-\xi^\top Q_{t,s}\xi},\qquad
M_{t,s}:=\int_t^s\sigma_0(\ell)\,dW_\ell,
\]
and let $\Gamma_{t,s}(\xi):=I_{t,s}(\xi)e^{i\xi\cdot M_{t,s}}$. It\^o's formula gives
\[
d\Gamma_{t,s}(\xi)
=-(\xi^\top a_0(s)\xi)\Gamma_{t,s}(\xi)\,ds
+i\Gamma_{t,s}(\xi)(\sigma_0(s)^\top\xi)\cdot dW_s.
\]
Applying the It\^o product formula to $\Gamma_{t,s}\widehat u_j$ yields
\begin{equation}\label{eq:old-4-7}
\begin{aligned}
d(\Gamma_{t,s}\widehat u_j)
={}&-\Gamma_{t,s}\widehat F_j\,ds
+\Gamma_{t,s}\bigl[\widehat v_j+i(\sigma_0^\top\xi)\widehat u_j\bigr]\cdot dW_s.
\end{aligned}
\end{equation}
For each $\xi\in\mathcal C_j\cap\mathbb Z^n$, $\Gamma_{t,s}(\xi)$ is adapted and $|\Gamma_{t,s}(\xi)|\le1$. The Sobolev integrability of $u_j,v_j$ therefore gives zero conditional mean for the stochastic integral in~\eqref{eq:old-4-7}.

The multiplier $\Gamma_{t,s}$ defines the frozen propagator $\mathcal P^0_{t,s}$. With the auxiliary cutoff $\widetilde m_j$ from~\eqref{eq:reproducing-cutoff}, write
\begin{equation}\label{eq:localized-heat-kernel-definition}
K_j(t,s,y):=\frac1{(2\pi)^n}\sum_{\xi\in\mathbb Z^n}
\widetilde m_j(\xi)e^{-\xi^\top Q_{t,s}\xi+i\xi\cdot y}.
\end{equation}
For $\operatorname{supp}\widehat h\subset\mathcal C_j\cap\mathbb Z^n$,
\[
\widehat{\mathcal P^0_{t,s}h}(\xi)=\Gamma_{t,s}(\xi)\widehat h(\xi)
\qquad\text{and}\qquad
(\mathcal P^0_{t,s}h)(z)=(K_j(t,s)*h)(z+M_{t,s}).
\]
Thus $\mathcal P^0_{t,s}$ combines the heat operator for $A_0$ with the random spatial shift $M_{t,s}$. This representation is the Fourier counterpart of the change of variables used for spatially constant coefficients; see~\cite[proof of Proposition~6.4]{DuQiuTang2012}.

Integrating~\eqref{eq:old-4-7}, conditioning on $\mathcal F_t$, and summing over $\xi\in\mathcal C_j\cap\mathbb Z^n$ gives
\begin{equation}\label{eq:old-4-8}
u_j(t,0)=\E_t\left[(K_j(t,T)*g_j)(M_{t,T})
+\int_t^T(K_j(t,s)*F_j(s))(M_{t,s})\,ds\right].
\end{equation}
The random kernels and spatial shifts remain inside the conditional expectation.

To estimate~\eqref{eq:old-4-8} in terms of $\mathcal E_{r_j}$, we use the following kernel bound, which combines spatial localization by $\kappa_{r_j}$ with exponential decay on the time scale $r_j^2$. The proof is given in Appendix~\ref{app:heat-kernel-proof}.
\begin{lemma}[Localized heat-kernel bounds]\label{lem:heat-kernel-decay}
Suppose that the matrix $Q_{t,s}$ in~\eqref{eq:localized-heat-kernel-definition} is symmetric and satisfies
\[
\lambda(s-t)I\le Q_{t,s}\le\Lambda(s-t)I,
\qquad 0\le t\le s\le T,
\]
for fixed $0<\lambda\le\Lambda$. Recall the localization weight $\kappa_r$ from~\eqref{eq:periodic-window-weight}. For $j\ge1$, $\ell=0,1$, and $0\le\nu\le2$, the periodic kernel $K_j$ satisfies
\begin{equation}\label{eq:old-4-9}
\rho(y)^\nu|D^\ell K_j(t,s,y)|
\le Cr_j^{\nu-\ell}e^{-c(s-t)/r_j^2}\kappa_{r_j}(y),
\qquad y\in\mathbb T^n,
\end{equation}
where $C=C(\lambda,\Lambda,\chi)$ and $c=c(\lambda)>0$.
\end{lemma}

\subsection{Conditional energy estimates under random spatial shifts}\label{subsec:displaced-energy}
In~\eqref{eq:old-4-8}, the convolution is evaluated at $M_{t,s}$, whereas $\mathcal E_r$ in~\eqref{eq:old-2-6} uses a center fixed during the time integral. The following lemma allows this fixed center to be chosen using the information available at the starting time; its proof is given in Appendix~\ref{app:window-calculus}.
\begin{lemma}\label{lem:random-center-energy}
Fix $0<r\le1$. Let $h$ be predictable as in~\eqref{eq:old-2-6}, with $\mathcal E_r(h)<\infty$. For every stopping time $\tau\le T$ and every $\mathcal F_\tau$-measurable $\mathbb T^n$-valued center $z$,
\begin{equation}\label{eq:random-center-energy}
\E_\tau\int_\tau^{(\tau+r^2)\wedge T}
\|h_s\|_{\kappa_{r,z}}^2\,ds
\le\mathcal E_r(h)^2
\qquad\text{a.s.}
\end{equation}
\end{lemma}

We now allow the center to move along a stochastic integral with a bounded predictable integrand.
\begin{lemma}\label{lem:displaced-energy}
Fix $0<r\le1$. Let $h$ be predictable as in~\eqref{eq:old-2-6}, with $\mathcal E_r(h)<\infty$, and let $\sigma$ be an $\R^{n\times m}$-valued predictable process with $|\sigma|\le K$. For a stopping time $\tau\le T$ and an $\mathcal F_\tau$-measurable $\mathbb T^n$-valued center $z$, set
\[
M_{\tau,s}:=\int_\tau^s\sigma_\ell\,dW_\ell,\qquad
X_s:=z+M_{\tau,s}\pmod{2\pi},\qquad \tau\le s\le T.
\]
Then
\begin{equation}\label{eq:displaced-local-energy}
\E_\tau\int_\tau^{(\tau+r^2)\wedge T}
\|h_s\|_{\kappa_{r,X_s}}^2\,ds
\le C\mathcal E_r(h)^2.
\end{equation}
Moreover, for every $c>0$,
\begin{equation}\label{eq:displaced-square-energy}
\E_\tau\int_\tau^T e^{-c(s-\tau)/r^2}
\|h_s\|_{\kappa_{r,X_s}}^2\,ds
\le C_c\mathcal E_r(h)^2,
\end{equation}
and
\begin{equation}\label{eq:displaced-linear-energy}
\E_\tau\int_\tau^T e^{-c(s-\tau)/r^2}(1+|M_{\tau,s}|/r)
\|h_s\|_{\kappa_{r,X_s}}\,ds
\le C_c r\mathcal E_r(h).
\end{equation}
The constants depend only on $n,K$ and, for $C_c$, on $c$. 
\end{lemma}

\begin{proof}
\emph{Step 1: Energy in a fixed window.} The explicit formula for $w_r$ implies $w_r(z-d)\le C(1+|d|/r)^{n+1}w_r(z)$ for $z,d\in\R^n$. Applying this inequality with $z=y+2\pi\ell$ and summing over $\ell\in\mathbb Z^n$ gives
\begin{equation}\label{eq:window-shift-comparison}
\kappa_r(y-d)\le C(1+|d|/r)^{n+1}\kappa_r(y),
\qquad d\in\R^n.
\end{equation}
For a stopping time $\tau$ and an $\mathcal F_\tau$-measurable center $z$, put $\zeta=(\tau+r^2)\wedge T$ and
\[
J_{\tau,z}:=\int_\tau^\zeta\|h_s\|_{\kappa_{r,z}}^2\,ds.
\]
By~\eqref{eq:random-center-energy}, $\E_\tau J_{\tau,z}\le\mathcal E_r(h)^2$.

Write $q_s=\|h_s\|_{\kappa_{r,z}}^2$ for $s\ge\tau$. On $\{\tau\le s\le\zeta\}$, both $q_s$ and $z$ are $\mathcal F_s$-measurable, and $\zeta\le(s+r^2)\wedge T$. Applying the same bound at time $s$ therefore gives $\E_s\int_s^\zeta q_\ell\,d\ell\le\mathcal E_r(h)^2$.

For bounded $q$, expanding the square of the time integral and conditioning first on $\mathcal F_s$ gives
\begin{equation}\label{eq:window-second-moment}
\E_\tau J_{\tau,z}^2
=2\E_\tau\int_\tau^\zeta q_s
\E_s\left[\int_s^\zeta q_\ell\,d\ell\right]ds
\le2\mathcal E_r(h)^2\E_\tau J_{\tau,z}
\le2\mathcal E_r(h)^4.
\end{equation}
For general $q$, set $q_s^{(N)}:=\min\{q_s,N\}$. These truncations are bounded and satisfy the same conditional energy bounds, since $0\le q_s^{(N)}\le q_s$. Applying the preceding calculation to $q^{(N)}$ and then letting $N\to\infty$, conditional monotone convergence yields~\eqref{eq:window-second-moment} for $q$.

\smallskip\noindent\emph{Step 2: Energy in a moving window.} Set
\[
M_{\tau,\zeta}^{*}:=\sup_{\tau\le s\le\zeta}|M_{\tau,s}|.
\]
The conditional Burkholder--Davis--Gundy inequality and $|\sigma|\le K$ give
\[
\E_\tau(1+M_{\tau,\zeta}^{*}/r)^{2(n+1)}\le C.
\]
By~\eqref{eq:window-shift-comparison}, conditional Cauchy--Schwarz, and~\eqref{eq:window-second-moment},
\[
\begin{aligned}
\E_\tau\int_\tau^\zeta\|h_s\|_{\kappa_{r,X_s}}^2\,ds
&\le C\E_\tau\bigl[(1+M_{\tau,\zeta}^{*}/r)^{n+1}J_{\tau,z}\bigr]\\
&\le C\bigl[\E_\tau(1+M_{\tau,\zeta}^{*}/r)^{2(n+1)}\bigr]^{1/2}
\bigl[\E_\tau J_{\tau,z}^2\bigr]^{1/2}\\
&\le C\mathcal E_r(h)^2.
\end{aligned}
\]
This proves~\eqref{eq:displaced-local-energy}.

\smallskip\noindent\emph{Step 3: Estimates on the full time interval.} Set $\tau_\nu=(\tau+\nu r^2)\wedge T$. Apply~\eqref{eq:displaced-local-energy} on each interval $[\tau_\nu,\tau_{\nu+1}]$ with initial center $X_{\tau_\nu}$, which is $\mathcal F_{\tau_\nu}$-measurable. On every nonempty such interval, $e^{-c(s-\tau)/r^2}\le e^{-c\nu}$. Conditioning at $\tau$ and summing the geometric series proves~\eqref{eq:displaced-square-energy}.

Finally, conditional It\^o isometry gives
\[
\E_\tau|M_{\tau,(\tau+q)\wedge T}|^2\le K^2q,\qquad q\ge0.
\]
Changing variables to $q=s-\tau$ and extending the nonnegative integral to $[0,\infty)$ yields
\[
\begin{aligned}
\E_\tau\int_\tau^T e^{-c(s-\tau)/r^2}
(1+|M_{\tau,s}|/r)^2\,ds
&\le2\int_0^\infty e^{-cq/r^2}(1+K^2q/r^2)\,dq\\
&\le C_c r^2.
\end{aligned}
\]
Conditional Cauchy--Schwarz with~\eqref{eq:displaced-square-energy} proves~\eqref{eq:displaced-linear-energy}.
\end{proof}

Taking $\tau=t$, $z=0$, and $\sigma=\sigma_0$ recovers the random spatial shift in the frozen representation. By~\eqref{eq:old-4-9}, Cauchy--Schwarz in space, and Lemma~\ref{lem:displaced-energy},
\begin{equation}\label{eq:kernel-window-bound}
\begin{aligned}
\left|\E_t\int_t^T(K_j(t,s)*h_s)(M_{t,s})ds\right|
&\le C\E_t\int_t^T e^{-c(s-t)/r_j^2}\|h_s\|_{\kappa_{r_j,M_{t,s}}}ds\\
&\le Cr_j\mathcal E_{r_j}(h).
\end{aligned}
\end{equation}
For every fixed $x\in\mathbb T^n$, estimate~\eqref{eq:kernel-window-bound} also holds with the convolution evaluated at $x+M_{t,s}$ and the window centered at $x+M_{t,s}$. The constant $C$ is independent of $x$, since $\mathcal E_{r_j}$ is invariant under spatial translations.

The representation and the local Sobolev estimate now give the following a priori bound for the frozen model.
\begin{lemma}\label{lem:frozen-model-apriori}
Let $(u_j,v_j)$ be the Sobolev strong solution of~\eqref{eq:frozen-model} under the assumptions of this section, with $j\ge1$. Suppose in addition that $g_j\in L^\infty_{\omega,x}$ and $\mathcal E_{r_j}(F_j)<\infty$. Then
\begin{equation}\label{eq:frozen-model-apriori}
\begin{aligned}
&\|u_j\|_{L^\infty_{\omega,t,x}}
+\mathcal E_{r_j}(Du_j)+\mathcal E_{r_j}(v_j)\\
&\qquad\le C\left(\|g_j\|_{L^\infty_{\omega,x}}
+r_j\mathcal E_{r_j}(F_j)\right),
\end{aligned}
\end{equation}
where $C=C(\lambda,K)$ is independent of $j$.
\end{lemma}
\begin{proof}
The unit mass of $\kappa_{r_j}$ and~\eqref{eq:old-4-9} bound the terminal term in~\eqref{eq:old-4-8} by $C\|g_j\|_{L^\infty_{\omega,x}}$. Combining the representation with~\eqref{eq:kernel-window-bound} at every spatial center therefore gives
\[
\|u_j\|_{L^\infty_{\omega,t,x}}
\le C\left(\|g_j\|_{L^\infty_{\omega,x}}
+r_j\mathcal E_{r_j}(F_j)\right).
\]
In particular, $u_j$ is essentially bounded, so Lemma~\ref{lem:weighted-window-energy} applies. Substituting the preceding bound into~\eqref{eq:model-energy} and taking square roots controls both conditional energies by the same right-hand side, proving~\eqref{eq:frozen-model-apriori}.
\end{proof}

For bounded $F_j$, the inequality $r_j\mathcal E_{r_j}(F_j)\le r_j^2\|F_j\|_{L^\infty_{\omega,t,x}}$ shows that~\eqref{eq:frozen-model-apriori} gives the usual gain of two spatial derivatives for the forcing term.

\section{Perturbation estimates for variable coefficients}\label{sec:kernel}
We apply the weighted local Sobolev estimate~\eqref{eq:model-energy} and the representation~\eqref{eq:old-4-8} to the full equation~\eqref{eq:full}, including all lower-order terms.

\begin{assumption}\label{ass:smooth-solution}
The pair $(u,v)$ is a Sobolev strong solution under the hypotheses of Theorem~\ref{thm:linear}. For almost every $(\omega,t)$, $u(\omega,t,\cdot)$ and each component of $v(\omega,t,\cdot)$ belong to $C^\infty(\mathbb T^n)$. Moreover,
\[
\|u\|_{L^\infty_{\omega,t}C_x^{2+\alpha}}+\|v\|_{\mathcal E}<\infty.
\]
\end{assumption}

\subsection{Coefficient freezing and commutator estimates}\label{subsec:perturbations}
Define
\[
a_0(s)=a(s,0),\qquad \sigma_0(s)=\sigma(s,0),\qquad
f_j=\Delta_jf,\quad g_j=\Delta_jg.
\]
For the commutator $[\Delta_j,h]q:=\Delta_j(hq)-h\Delta_jq$, write
\begin{equation}\label{eq:old-2-27}
R_j^a:=[\Delta_j,a]:D^2u,\qquad R_j^\sigma:=[\Delta_j,\sigma]:Dv.
\end{equation}
The block equation has exactly the model form~\eqref{eq:frozen-model}, with
\begin{equation}\label{eq:total-model-forcing}
F_j=f_j+\mathcal R_j^u+\mathcal R_j^v,
\end{equation}
where
\begin{equation}\label{eq:grouped-perturbations}
\left\{
\begin{aligned}
&\mathcal R_j^u=(a-a_0):D^2u_j+R_j^a+\Delta_j(b\cdot Du+cu),\\
&\mathcal R_j^v=(\sigma-\sigma_0):Dv_j+R_j^\sigma+\Delta_j(\gamma\cdot v).
\end{aligned}
\right.
\end{equation}
The representation~\eqref{eq:old-4-8} requires Fourier support in $\mathcal C_j$ only for the total forcing $F_j$, not for the individual terms in~\eqref{eq:grouped-perturbations}. Indeed,
\[
F_j=\Delta_j(a:D^2u+\sigma:Dv+b\cdot Du+cu+\gamma\cdot v+f)
-a_0:D^2u_j-\sigma_0:Dv_j
\]
has Fourier support in $\mathcal C_j\cap\mathbb Z^n$. Thus convolution with $K_j$ in~\eqref{eq:old-4-8} may be distributed over all terms of~\eqref{eq:total-model-forcing}, even though the individual terms need not have this support.

For the martingale perturbations, we introduce the global conditional energy
\begin{equation}\label{eq:global-conditional-energy}
\mathcal B_1(v):=\sup_\tau\left\|
\left(\E_\tau\int_\tau^T\|v_s\|_{H_x^1}^2\,ds\right)^{1/2}
\right\|_{L^\infty_\omega},
\end{equation}
where the supremum is over all stopping times taking values in $[0,T]$. Equivalently, $\mathcal B_1(v)$ is the $\mathrm{BMO}_2$ norm of the $H_x^1$-valued martingale $\int_0^\cdot v_s\cdot dW_s$.

We estimate $R_j^a+\Delta_j(b\cdot Du+cu)$ by $\delta\|u\|_{L^\infty_{\omega,t}C_x^{2+\alpha}}$ and $\|u\|_{L^\infty_{\omega,t,x}}$. For $R_j^\sigma$ and $\Delta_j(\gamma\cdot v)$, we use $\mathcal E_{r_j}$, separating a small multiple of $\|v\|_{\mathcal E}$ from a multiple of $\mathcal B_1(v)$. No smallness of the coefficients is assumed. The frequency decompositions proving these bounds are in Appendix~\ref{app:commutator-proofs}.
\begin{proposition}\label{prop:commutators}
Under Assumption~\ref{ass:smooth-solution}, there exists a deterministic constant $C=C(\lambda,K)>0$ with the following properties. For every $\delta>0$, there exists a deterministic constant $C_\delta=C_\delta(\lambda,K)>0$ such that, for all $j\ge0$,
\begin{equation}\label{eq:perturbation-u}
\|R_j^a+\Delta_j(b\cdot Du+cu)\|_{L^\infty_{\omega,t,x}}
\le Cr_j^\alpha\left(\delta\|u\|_{L^\infty_{\omega,t}C_x^{2+\alpha}}+C_\delta\|u\|_{L^\infty_{\omega,t,x}}\right).
\end{equation}
For all $j\ge0$, the martingale terms satisfy
\begin{equation}\label{eq:perturbation-v-energy}
\mathcal E_{r_j}(R_j^\sigma)+\mathcal E_{r_j}(\Delta_j(\gamma\cdot v))
\le Cr_j^{1+\alpha}\|v\|_{\mathcal E}.
\end{equation}
Moreover, for every $\varepsilon>0$, there exists a deterministic constant $C_\varepsilon=C_\varepsilon(\lambda,K)>0$ such that
\begin{equation}\label{eq:perturbation-v-small}
\sup_{j\ge0}r_j^{-1-\alpha}
\left[e_j+\mathcal E_{r_j}(R_j^\sigma)+\mathcal E_{r_j}(\Delta_j(\gamma\cdot v))\right]
\le\varepsilon\|v\|_{\mathcal E}+C_\varepsilon\mathcal B_1(v),
\end{equation}
where $e_j=\mathcal E_{r_j}(v_j)$ and $\mathcal B_1(v)$ is defined in~\eqref{eq:global-conditional-energy}.
\end{proposition}
In~\eqref{eq:perturbation-v-small}, the high-frequency blocks of $v$ contribute at most $\varepsilon\|v\|_{\mathcal E}$, while its low-frequency blocks contribute at most $C_\varepsilon\mathcal B_1(v)$. Estimate~\eqref{eq:old-2-23} controls $\mathcal B_1(v)$ independently of $\|u\|_{L^\infty_{\omega,t}C_x^{2+\alpha}}$.

\subsection{Conditional local Sobolev estimate}\label{subsec:dyadic-energy}
We apply~\eqref{eq:model-energy} with the forcing $F_j$ in~\eqref{eq:total-model-forcing}. Proposition~\ref{prop:commutators} then bounds $\|v\|_{\mathcal E}$ by $\|u\|_{L^\infty_{\omega,t}C_x^{2+\alpha}}$ and the norms of $f,g$.
\begin{proposition}\label{prop:energy}
Under Assumption~\ref{ass:smooth-solution}, there exists a deterministic constant $C_E=C_E(\lambda,K)>0$ such that
\begin{equation}\label{eq:old-3-1}
\|v\|_{\mathcal E}\le C_E\left(\|u\|_{L^\infty_{\omega,t}C_x^{2+\alpha}}
+\|f\|_{L^\infty_{\omega,t}C_x^\alpha}+\|g\|_{L^\infty_\omega C_x^{2+\alpha}}\right).
\end{equation}
\end{proposition}
\begin{proof}
For this and subsequent proofs, we use the abbreviation
\begin{equation}\label{eq:old-2-1}
\|h\|_{\mathcal S}:=\|h\|_{L^\infty_{\omega,t}C_x^{2+\alpha}}.
\end{equation}
Fix $j\ge1$ and put $r=r_j$. For $h\in L^\infty_{\omega,t,x}$, $\mathcal E_r(h)\le r\|h\|_{L^\infty_{\omega,t,x}}$. Boundedness of $a-a_0$ and Bernstein's inequality~\eqref{eq:bernstein-infty} therefore give
\[
\begin{aligned}
\mathcal E_r((a-a_0):D^2u_j)
&\le C\mathcal E_r(D^2u_j)
\le Cr\|D^2u_j\|_{L^\infty_{\omega,t,x}}\\
&\le Cr^{-1}\|u_j\|_{L^\infty_{\omega,t,x}}.
\end{aligned}
\]
Similarly, using boundedness of $\sigma-\sigma_0$ first and then applying Lemma~\ref{lem:window-scales} to $v_j$, we obtain
\[
\mathcal E_r((\sigma-\sigma_0):Dv_j)
\le C\mathcal E_r(Dv_j)
\le Cr^{-1}\mathcal E_r(v_j)=Cr^{-1}e_j.
\]
Apply the bound for bounded functions to $f_j$ and~\eqref{eq:perturbation-u} with $\delta=1$, and use~\eqref{eq:perturbation-v-energy}. This gives
\[
\begin{aligned}
\mathcal E_r(F_j)\le{}&Cr^{-1}(\|u_j\|_{L^\infty_{\omega,t,x}}+e_j)\\
&+Cr^{1+\alpha}\left(\|f\|_{L^\infty_{\omega,t}C_x^\alpha}+\|u\|_{\mathcal S}+\|v\|_{\mathcal E}\right).
\end{aligned}
\]
Substitution into the local Sobolev estimate for the model~\eqref{eq:model-energy}, followed by Young's inequality for $\|u_j\|_{L^\infty_{\omega,t,x}}e_j$, yields
\[
e_j^2\le C\|u_j\|_{L^\infty_{\omega,t,x}}^2+Cr^{2+\alpha}\|u_j\|_{L^\infty_{\omega,t,x}}\left(\|f\|_{L^\infty_{\omega,t}C_x^\alpha}+\|u\|_{\mathcal S}+\|v\|_{\mathcal E}\right).
\]
By Lemma~\ref{lem:frequency}, $\|u_j\|_{L^\infty_{\omega,t,x}}\le Cr_j^{2+\alpha}\|u\|_{\mathcal S}$. Divide by $r_j^{4+2\alpha}$ and take the supremum over $j\ge1$. The block $j=0$ is bounded independently by Lemma~\ref{lem:low}. Hence
\[
\|v\|_{\mathcal E}^2\le C\left[\|u\|_{\mathcal S}^2
+\left(\|f\|_{L^\infty_{\omega,t}C_x^\alpha}+\|g\|_{L^\infty_\omega C_x^{2+\alpha}}\right)^2
+\|u\|_{\mathcal S}\|v\|_{\mathcal E}\right].
\]
Applying Young's inequality to the mixed term gives~\eqref{eq:old-3-1}.
\end{proof}

\subsection{Perturbative Schauder estimate}
To estimate $\|u_j\|_{L^\infty_{\omega,t,x}}$, we retain $a-a_0$ and $\sigma-\sigma_0$ from~\eqref{eq:grouped-perturbations} inside the convolution in~\eqref{eq:old-4-8}. Lemma~\ref{lem:localized-heat-kernel} records the resulting factors $r_j^\alpha$ and $r_j$ for $K_j,M_{t,s}$ associated with $a_0(s)=a(s,0)$ and $\sigma_0(s)=\sigma(s,0)$.
\begin{lemma}\label{lem:localized-heat-kernel}
For $j\ge1$ and $t<T$,
\begin{equation}\label{eq:coefficient-variation-bounds}
\left\{
\begin{aligned}
&\left|\E_t\int_t^T\bigl(K_j(t,s)*[(a_s-a_0(s)):D^2u_j(s)]\bigr)(M_{t,s})ds\right|
\le Cr_j^\alpha\|u_j\|_{L^\infty_{\omega,t,x}},\\[0.6em]
&\left|\E_t\int_t^T\bigl(K_j(t,s)*[(\sigma_s-\sigma_0(s)):Dv_j(s)]\bigr)(M_{t,s})ds\right|
\le Cr_j e_j.
\end{aligned}
\right.
\end{equation}
\end{lemma}
\begin{proof}
Put $r=r_j$. Spatial H\"older continuity and the conditional moment bound give
\[
|a(s,M_{t,s}+y)-a(s,0)|\le C(\rho(y)^\alpha+|M_{t,s}|^\alpha),
\qquad \E_t|M_{t,s}|^\alpha\le C(s-t)^{\alpha/2}.
\]
By Bernstein's inequality~\eqref{eq:bernstein-infty} and~\eqref{eq:old-4-9}, the first integral is bounded by
\[
Cr^{-2}\|u_j\|_{L^\infty_{\omega,t,x}}
\int_t^T e^{-c(s-t)/r^2}(r^\alpha+(s-t)^{\alpha/2})ds
\le Cr^\alpha\|u_j\|_{L^\infty_{\omega,t,x}}.
\]
For the second integral, put $B_s(y)=\sigma(s,M_{t,s}+y)-\sigma(s,0)$. Then
\[
|B_s(y)|\le C(\rho(y)+|M_{t,s}|),\qquad |DB_s(y)|\le C.
\]
Evenness of $K_j$ and integration by parts give
\[
\bigl(K_j*[(\sigma_s-\sigma_0(s)):Dv_j(s)]\bigr)(M_{t,s})
=-\int\operatorname{div}(K_jB_s)(y)\cdot v_j(s,M_{t,s}+y)dy.
\]
The product rule and the bounds for $B_s$ give
\[
\begin{aligned}
|\operatorname{div}(K_jB_s)(y)|
&\le |DK_j(t,s,y)|\,|B_s(y)|+|K_j(t,s,y)|\,|\operatorname{div}B_s(y)|\\
&\le C\bigl[\rho(y)|DK_j(t,s,y)|+|M_{t,s}|\,|DK_j(t,s,y)|+|K_j(t,s,y)|\bigr]\\
&\le Ce^{-c(s-t)/r^2}(1+|M_{t,s}|/r)\kappa_r(y),
\end{aligned}
\]
where the last inequality applies~\eqref{eq:old-4-9} with $(\ell,\nu)=(1,1),(1,0),(0,0)$ to the three terms, respectively.
Spatial Cauchy--Schwarz and~\eqref{eq:displaced-linear-energy}, with $h=v_j$, now bound the time integral by $Cr\mathcal E_r(v_j)$.
\end{proof}

Combining~\eqref{eq:coefficient-variation-bounds} with Proposition~\ref{prop:commutators} gives the following perturbative Schauder bound, which will be used in Section~\ref{sec:closure}.
\begin{proposition}\label{prop:kernel}
Under Assumption~\ref{ass:smooth-solution}, for every $\varepsilon>0$, there exists a deterministic constant $C_\varepsilon=C_\varepsilon(\lambda,K)>0$ such that
\begin{equation}\label{eq:old-4-1}
\|u\|_{L^\infty_{\omega,t}C_x^{2+\alpha}}\le C_\varepsilon\left(\|f\|_{L^\infty_{\omega,t}C_x^\alpha}
+\|g\|_{L^\infty_\omega C_x^{2+\alpha}}\right)+\varepsilon\|v\|_{\mathcal E}.
\end{equation}
\end{proposition}
\begin{proof}
Apply the model representation~\eqref{eq:old-4-8} with the forcing~\eqref{eq:total-model-forcing}. The terminal term is bounded by $Cr_j^{2+\alpha}\|g\|_{L^\infty_\omega C_x^{2+\alpha}}$. For the forcing block, Lemma~\ref{lem:frequency} gives
\[
\|f_j\|_{L^\infty_{\omega,t,x}}\le Cr_j^\alpha\|f\|_{L^\infty_{\omega,t}C_x^\alpha}.
\]

The kernel has time-integrated $L^1$ norm at most $Cr_j^2$. Lemma~\ref{lem:localized-heat-kernel} controls the coefficient differences, while~\eqref{eq:kernel-window-bound} and~\eqref{eq:perturbation-u} control the remaining $v$ and $u$ terms, respectively. Together, these bounds give, uniformly in the freezing point,
\begin{equation}\label{eq:full-block-estimate}
\begin{aligned}
\|u_j\|_{L^\infty_{\omega,t,x}}\le{}&Cr_j^{2+\alpha}\left(
\|f\|_{L^\infty_{\omega,t}C_x^\alpha}+\|g\|_{L^\infty_\omega C_x^{2+\alpha}}
+\delta\|u\|_{\mathcal S}+C_\delta\|u\|_{L^\infty_{\omega,t,x}}\right)\\
&+Cr_j^\alpha\|u_j\|_{L^\infty_{\omega,t,x}}
+Cr_j\left[e_j+\mathcal E_{r_j}(R_j^\sigma)+\mathcal E_{r_j}(\Delta_j(\gamma\cdot v))\right].
\end{aligned}
\end{equation}
The conditional Sobolev bound~\eqref{eq:old-2-23} and the uniform bound~\eqref{eq:old-2-22} give
\begin{equation}\label{eq:independent-data-bound}
\|u\|_{L^\infty_{\omega,t,x}}+\mathcal B_1(v)
\le C\left(\|f\|_{L^\infty_{\omega,t}C_x^\alpha}+\|g\|_{L^\infty_\omega C_x^{2+\alpha}}\right).
\end{equation}
Fix $\varepsilon>0$. Normalize~\eqref{eq:full-block-estimate} by $r_j^{2+\alpha}$ and split the blocks at an integer $J\ge1$. For $j\ge J$, the coefficient-variation term contributes at most $Cr_J^\alpha\|u\|_{\mathcal S}$, while~\eqref{eq:perturbation-v-small}, with parameter $\eta>0$, controls the martingale terms by $C\eta\|v\|_{\mathcal E}+C_\eta\mathcal B_1(v)$. The blocks $j<J$ are bounded by the data using~\eqref{eq:old-2-25}. Combining these bounds with~\eqref{eq:independent-data-bound} and Lemma~\ref{lem:frequency} yields
\[
\begin{aligned}
\|u\|_{\mathcal S}\le{}&C_{\delta,J,\eta}\left(
\|f\|_{L^\infty_{\omega,t}C_x^\alpha}+\|g\|_{L^\infty_\omega C_x^{2+\alpha}}\right)\\
&+C_0(\delta+r_J^\alpha)\|u\|_{\mathcal S}+C_0\eta\|v\|_{\mathcal E},
\end{aligned}
\]
where $C_0\ge1$ depends only on the structural parameters. Choose $\delta>0$ small and $J$ large, and take $\eta>0$ so that
\[
C_0(\delta+r_J^\alpha)\le\tfrac12,\qquad
\eta=\frac{\varepsilon}{2C_0}.
\]
With this choice, the preceding inequality gives~\eqref{eq:old-4-1}.
\end{proof}

\section{Schauder estimates for Sobolev strong solutions}\label{sec:closure}
We prove Theorems~\ref{thm:linear} and~\ref{thm:joint-estimate} by applying Propositions~\ref{prop:energy} and~\ref{prop:kernel} to spatially mollified equations, then passing to the original Sobolev strong solution by Sobolev convergence and lower semicontinuity.

We approximate the coefficients and data of~\eqref{eq:full} by spatial convolution. Choose $\varrho\in C_c^\infty(\R^n)$ with $\varrho\ge0$, $\operatorname{supp}\varrho\subset B(0,1)$, and $\int_{\R^n}\varrho(z)\,dz=1$. For $\epsilon>0$, define its scaled periodization by
\[
\varrho_\epsilon(x):=\sum_{k\in\mathbb Z^n}\epsilon^{-n}
\varrho\!\left(\frac{x+2\pi k}{\epsilon}\right),
\qquad x\in\mathbb T^n.
\]
The sum is independent of the representative of $x$ and defines a smooth nonnegative periodic kernel. Unfolding the integral gives
\[
\int_{\mathbb T^n}\varrho_\epsilon(x)\,dx
=\int_{\R^n}\epsilon^{-n}\varrho(z/\epsilon)\,dz=1.
\]
We set
\begin{equation}\label{eq:old-5-4}
h^\epsilon:=\varrho_\epsilon*h,
\qquad h\in\{a,\sigma,b,c,\gamma,f,g\}.
\end{equation}
The approximating BSPDE is
\begin{equation}\label{eq:mollified-full}
\left\{
\begin{aligned}
&-du^\epsilon=\bigl[a^\epsilon:D^2u^\epsilon+\sigma^\epsilon:Dv^\epsilon
+b^\epsilon\cdot Du^\epsilon+c^\epsilon u^\epsilon\\
&\hphantom{-du^\epsilon=\bigl[}{}+\gamma^\epsilon\cdot v^\epsilon+f^\epsilon\bigr]dt-v^\epsilon\cdot dW_t,\\
&u^\epsilon(T,\cdot)=g^\epsilon.
\end{aligned}
\right.
\end{equation}
Spatial convolution preserves predictability and terminal measurability, as well as the coefficient and data bounds. Jensen's inequality gives $\sigma^\epsilon(\sigma^\epsilon)^\top\le(\sigma\sigma^\top)*\varrho_\epsilon$, so
\begin{equation}\label{eq:old-5-5}
a^\epsilon-\frac12\sigma^\epsilon(\sigma^\epsilon)^\top
\ge\left(a-\frac12\sigma\sigma^\top\right)*\varrho_\epsilon
\ge\lambda I.
\end{equation}
Thus Lemma~\ref{lem:periodic-existence} gives a unique Sobolev strong solution $(u^\epsilon,v^\epsilon)$ of~\eqref{eq:mollified-full}. For each fixed $\epsilon>0$, the higher-order version of Lemma~\ref{lem:periodic-existence}, applied to the smooth coefficients and data, gives, for every integer $q\ge1$,
\begin{equation}\label{eq:mollified-higher-order}
\begin{aligned}
&\|u^\epsilon\|_{L^\infty_{\omega,t}H_x^q}^2
+\sup_{\tau\le T}\left\|\E_\tau\int_\tau^T\|v^\epsilon_s\|_{H_x^q}^2\,ds\right\|_{L^\infty_\omega}\\
&\le C_{\epsilon,q}\left(\|f\|_{L^\infty_{\omega,t}C_x^\alpha}^2
+\|g\|_{L^\infty_\omega C_x^{2+\alpha}}^2\right).
\end{aligned}
\end{equation}
Since $q$ is arbitrary, these estimates give spatial smoothness of $(u^\epsilon,v^\epsilon)$.

Choose an integer $q>n/2+2+\alpha$. Sobolev embedding gives $\|u^\epsilon\|_{L^\infty_{\omega,t}C_x^{2+\alpha}}<\infty$. For $j\ge1$, the annulus $\mathcal C_j$ contains at most $Cr_j^{-n}$ lattice points, all with $|\xi|\asymp r_j^{-1}$. Cauchy--Schwarz and Parseval therefore give
\[
\|\Delta_jv^\epsilon_s\|_\infty
\le Cr_j^{-n/2}\|\Delta_jv^\epsilon_s\|_{L_x^2}
\le Cr_j^{q-n/2}\|v^\epsilon_s\|_{H_x^q}.
\]
The same bound holds for $j=0$, since $r_0=1$ and $\Delta_0$ has finite Fourier support. Together with~\eqref{eq:mollified-higher-order}, this implies, for every $j\ge0$,
\[
\mathcal E_{r_j}(\Delta_jv^\epsilon)
\le C_{\epsilon,q}r_j^{q-n/2}
\left(\|f\|_{L^\infty_{\omega,t}C_x^\alpha}
+\|g\|_{L^\infty_\omega C_x^{2+\alpha}}\right).
\]
Since $q-n/2>2+\alpha$, the definition of $\|\cdot\|_{\mathcal E}$ gives $\|v^\epsilon\|_{\mathcal E}<\infty$. Thus $(u^\epsilon,v^\epsilon)$ satisfies Assumption~\ref{ass:smooth-solution} for each fixed $\epsilon>0$.

As $\epsilon\downarrow0$, the mollified coefficients converge in $L^\infty_{\omega,t,x}$, while $f^\epsilon\to f$ in $L^2_{\omega,t,x}$ and $g^\epsilon\to g$ in $L^2_\omega H_x^1$. Corollary~\ref{cor:sobolev-stability} therefore gives
\begin{equation}\label{eq:old-5-7}
u^\epsilon\longrightarrow u\quad\mbox{in }L^2_{\omega,t}H_x^2,
\qquad
v^\epsilon\longrightarrow v\quad\mbox{in }L^2_{\omega,t}H_x^1,
\end{equation}
where $(u,v)$ is the unique Sobolev strong solution of~\eqref{eq:full}.

The convergence in~\eqref{eq:old-5-7} does not imply convergence in $\|\cdot\|_{L^\infty_{\omega,t}C_x^{2+\alpha}}$ or $\|\cdot\|_{\mathcal E}$. Lemma~\ref{lem:energy-limit} establishes the lower semicontinuity of these two norms under Sobolev convergence.
\begin{lemma}\label{lem:energy-limit}
With the weights $\kappa_r$ fixed, if $u^\nu\to u$ in $L^2_{\omega,t}H_x^2$ and $v^\nu\to v$ in $L^2_{\omega,t}H_x^1$, then
\begin{equation}\label{eq:energy-limit-bound}
\begin{aligned}
\|u\|_{L^\infty_{\omega,t}C_x^{2+\alpha}}+\|v\|_{\mathcal E}
\le\liminf_{\nu\to\infty}\bigl(
\|u^\nu\|_{L^\infty_{\omega,t}C_x^{2+\alpha}}+\|v^\nu\|_{\mathcal E}\bigr).
\end{aligned}
\end{equation}
\end{lemma}

\begin{proof}
Assume the right-hand side of~\eqref{eq:energy-limit-bound} is finite. Choose a subsequence with $\|u^\nu\|_{\mathcal S}\to U_*$ and $\|v^\nu\|_{\mathcal E}\to V_*$, where $U_*+V_*$ equals that lower limit.

After a further subsequence, $u^\nu(\omega,t)\to u(\omega,t)$ in $L_x^2$ for almost every $(\omega,t)$. Fix a pair $(\omega,t)$ for which this convergence and the bounds $\|u^\nu(\omega,t)\|_{C^{2+\alpha}}\le\|u^\nu\|_{\mathcal S}$ hold for every $\nu$. This sequence is bounded in $C^{2+\alpha}$, so compactness gives a subsequence $u^{\nu_k}(\omega,t)$ converging in $C^2$ to some $h$. Since this also implies convergence in $L_x^2$, uniqueness of the $L_x^2$ limit gives $h=u(\omega,t)$. Lower semicontinuity of the $\alpha$-H\"older seminorm of the second derivatives then yields
\[
\|u(\omega,t)\|_{C^{2+\alpha}}
\le\liminf_{k\to\infty}\|u^{\nu_k}(\omega,t)\|_{C^{2+\alpha}}
\le U_*.
\]
Taking the essential supremum over $(\omega,t)$ gives $\|u\|_{\mathcal S}\le U_*$.

Fix $j\ge0$, a stopping time $\tau\le T$, and $x\in\mathbb T^n$, and set
\[
I^\nu:=\int_\tau^{(\tau+r_j^2)\wedge T}
\|\Delta_jv^\nu(s)\|_{\kappa_{r_j,x}}^2\,ds,
\]
and define $I$ with $v$ in place of $v^\nu$. Sobolev convergence and boundedness of $\kappa_{r_j}$ give $I^\nu\to I$ in $L^1(\Omega)$, hence $\E_\tau I^\nu\to\E_\tau I$ in $L^1(\Omega)$. By the definition of $\|\cdot\|_{\mathcal E}$ and Lemma~\ref{lem:random-center-energy},
\[
\E_\tau I^\nu\le r_j^{4+2\alpha}\|v^\nu\|_{\mathcal E}^2.
\]
Passing to the limit gives $\E_\tau I\le r_j^{4+2\alpha}V_*^2$ almost surely. The suprema defining $\|\cdot\|_{\mathcal E}$ then yield $\|v\|_{\mathcal E}\le V_*$. Adding the two bounds proves~\eqref{eq:energy-limit-bound}.
\end{proof}

\begin{proof}[Proof of Theorems~\ref{thm:linear} and~\ref{thm:joint-estimate}]
Lemma~\ref{lem:periodic-existence} gives the unique Sobolev strong solution. Apply Propositions~\ref{prop:energy} and~\ref{prop:kernel} to the smooth approximations~\eqref{eq:mollified-full}. Their coefficient bounds are uniform in $\epsilon$. Since $\varrho_\epsilon$ is nonnegative and has mass one,
\[
\begin{aligned}
\|f^\epsilon\|_{L^\infty_{\omega,t}C_x^\alpha}
&\le\|f\|_{L^\infty_{\omega,t}C_x^\alpha},\\
\|g^\epsilon\|_{L^\infty_\omega C_x^{2+\alpha}}
&\le\|g\|_{L^\infty_\omega C_x^{2+\alpha}}.
\end{aligned}
\]
Thus, with constants independent of $\epsilon$,
\begin{equation}\label{eq:old-5-1}
\left\{
\begin{aligned}
&\|v^\epsilon\|_{\mathcal E}\le C_E\left(\|u^\epsilon\|_{\mathcal S}
+\|f\|_{L^\infty_{\omega,t}C_x^\alpha}+\|g\|_{L^\infty_\omega C_x^{2+\alpha}}\right),\\
&\|u^\epsilon\|_{\mathcal S}\le C_\varepsilon\left(
\|f\|_{L^\infty_{\omega,t}C_x^\alpha}+\|g\|_{L^\infty_\omega C_x^{2+\alpha}}\right)
+\varepsilon\|v^\epsilon\|_{\mathcal E},\qquad \varepsilon>0.
\end{aligned}
\right.
\end{equation}
The constant $C_E$ in the local Sobolev estimate is fixed before choosing $\varepsilon$. Substituting the first inequality into the second gives
\[
\|u^\epsilon\|_{\mathcal S}\le(C_\varepsilon+\varepsilon C_E)
\left(\|f\|_{L^\infty_{\omega,t}C_x^\alpha}+\|g\|_{L^\infty_\omega C_x^{2+\alpha}}\right)
+\varepsilon C_E\|u^\epsilon\|_{\mathcal S}.
\]
Choose $\varepsilon=[2(1+C_E)]^{-1}$. The preceding inequality and the first bound in~\eqref{eq:old-5-1} then yield
\begin{equation}\label{eq:old-5-6}
\|u^\epsilon\|_{\mathcal S}+\|v^\epsilon\|_{\mathcal E}
\le C\left(\|f\|_{L^\infty_{\omega,t}C_x^\alpha}+\|g\|_{L^\infty_\omega C_x^{2+\alpha}}\right),
\end{equation}
uniformly in the mollification scale.

By~\eqref{eq:old-5-7}, these approximations converge to the original Sobolev strong solution. Lemma~\ref{lem:energy-limit} therefore gives
\[
\begin{aligned}
\|u\|_{\mathcal S}+\|v\|_{\mathcal E}
&\le\liminf_{\epsilon\downarrow0}
(\|u^\epsilon\|_{\mathcal S}+\|v^\epsilon\|_{\mathcal E})\\
&\le C\left(\|f\|_{L^\infty_{\omega,t}C_x^\alpha}+\|g\|_{L^\infty_\omega C_x^{2+\alpha}}\right).
\end{aligned}
\]
This proves~\eqref{eq:linear-estimate} and~\eqref{eq:joint-estimate}.
\end{proof}

Applying Theorems~\ref{thm:linear} and~\ref{thm:joint-estimate} to the difference of two solutions with the same coefficients gives the following stability estimate.
\begin{corollary}\label{cor:stability}
Let $(u,v)$ and $(\bar u,\bar v)$ solve~\eqref{eq:full} with the same coefficients and data $(f,g)$ and $(\bar f,\bar g)$, respectively, under the hypotheses of Theorem~\ref{thm:linear}. Then
\begin{equation}\label{eq:old-5-19}
\|u-\bar u\|_{L^\infty_{\omega,t}C_x^{2+\alpha}}+\|v-\bar v\|_{\mathcal E}
\le C\Big(\|f-\bar f\|_{L^\infty_{\omega,t}C_x^\alpha}
+\|g-\bar g\|_{L^\infty_\omega C_x^{2+\alpha}}\Big).
\end{equation}
\end{corollary}

\section{Optimal regularity of the martingale coefficients}\label{sec:coefficient-regularity}
Fix $0<\alpha<1$. For every $1<\beta<1+\alpha$, we construct two real-valued scalar equations~\eqref{eq:full} on $[0,1]\times\mathbb T$, with $a=1$ and $b=c=f=0$. In each construction, the deterministic coefficients $\sigma,\gamma$ are fixed, and only the terminal datum $g_N$ varies with $N=2^J$, $J\ge1$. They satisfy
\begin{equation}\label{eq:sharpness-uniform-data}
\left\{
\begin{aligned}
&1-\tfrac12\sigma^2\ge\tfrac12,\\
&\|\sigma\|_{C^\beta}+\|\gamma\|_{C^\beta}
+\sup_{J\ge1}\|g_{2^J}\|_{L^\infty_\omega C_x^{2+\alpha}}<\infty.
\end{aligned}
\right.
\end{equation}
Each corresponding Sobolev strong solution belongs to $L^\infty_{\omega,t}C_x^{2+\alpha}$, but
\begin{equation}\label{eq:sharpness-divergence}
\|u_N\|_{L^\infty_{\omega,t}C_x^{2+\alpha}}\longrightarrow\infty
\qquad\text{as }N=2^J\to\infty.
\end{equation}
Thus the Schauder estimate~\eqref{eq:linear-estimate} fails in general. The constructions use a deterministic complex parabolic equation to represent real-valued BSPDE solutions.

\subsection{Parabolic representation and Fourier estimates}\label{subsec:sharpness-construction}
Let $T,a>0$ and $1<\theta<1+\alpha$ be fixed. Let $\sigma,\gamma\in C^\theta(\mathbb T)$ be fixed deterministic real-valued functions satisfying $a-\sigma^2/2\ge\lambda>0$, and let $\phi$ be a smooth deterministic complex-valued function. For each positive integer $N$, consider
\begin{equation}\label{eq:sharpness-bspde}
\left\{
\begin{aligned}
&-du_N=\bigl[a\partial_{xx}u_N+\sigma(x)\partial_xv_N+\gamma(x)v_N\bigr]dt-v_N\,dW_t,\\
&u_N(T,x)=g_N(x):=\operatorname{Re}\bigl(e^{iNW_T}\phi(x)\bigr).
\end{aligned}
\right.
\end{equation}
Let $z_N$ solve the deterministic complex equation
\begin{equation}\label{eq:sharpness-deterministic}
\left\{
\begin{aligned}
&\partial_s z_N=a\partial_{xx}z_N+iN\sigma\partial_xz_N+iN\gamma z_N,\\
&z_N(0,x)=\phi(x).
\end{aligned}
\right.
\end{equation}
For each fixed $N$, the real and imaginary parts form a uniformly parabolic system with $C^\theta$ lower-order coefficients. Since $\theta>\alpha$, deterministic parabolic regularity gives
\[
\sup_{0\le s\le T}\|z_N(s)\|_{C_x^{2+\alpha}}<\infty,
\]
with a bound that may depend on $N$. Writing $s=T-t$, It\^o's formula gives the real-valued Sobolev strong solution
\begin{equation}\label{eq:sharpness-solution}
\left\{
\begin{aligned}
&u_N(t,x)=e^{-N^2s/2}\operatorname{Re}\bigl(e^{iNW_t}z_N(s,x)\bigr),\\
&v_N(t,x)=-Ne^{-N^2s/2}\operatorname{Im}\bigl(e^{iNW_t}z_N(s,x)\bigr).
\end{aligned}
\right.
\end{equation}
Uniqueness follows from the $L^2$ energy estimate, using super-parabolicity and $\sigma\in C^1$.

For every real- or complex-valued $q\in C_x^{2+\alpha}(\mathbb T)$, we claim that
\begin{equation}\label{eq:sharpness-fourier}
\|q\|_{C_x^{2+\alpha}}\ge c|k|^{2+\alpha}|\widehat q(k)|,
\qquad k\in\mathbb Z\setminus\{0\}.
\end{equation}
Indeed, for $h=\pi/|k|$, the Fourier translation formula for the second derivative $q''$ gives
\[
2|k|^2|\widehat q(k)|
=\left|\widehat{q''(\cdot+h)}(k)-\widehat{q''}(k)\right|
\le[q'']_\alpha h^\alpha.
\]
The same argument, applied to the derivative of order $\lfloor\eta\rfloor$, gives, for every noninteger $\eta>0$,
\begin{equation}\label{eq:sharpness-fourier-decay}
|\widehat q(k)|\le C_\eta|k|^{-\eta}\|q\|_{C^\eta},
\qquad k\ne0.
\end{equation}

Taking real parts in~\eqref{eq:sharpness-solution} gives
\[
\widehat u_N(t,k)
=\frac{e^{-N^2s/2}}2
\left[e^{iNW_t}\widehat z_N(s,k)
+e^{-iNW_t}\overline{\widehat z_N(s,-k)}\right].
\]
For $t>0$, the Gaussian law of $W_t$ has full support on $\R$. Since
\[
\sup_{\varphi\in\R}|e^{i\varphi}A+e^{-i\varphi}B|=|A|+|B|,
\qquad A,B\in\mathbb C,
\]
we obtain
\begin{equation}\label{eq:sharpness-real-mode}
\begin{aligned}
\|\widehat u_N(t,k)\|_{L^\infty_\omega}
&=\frac{e^{-N^2s/2}}2
\left(|\widehat z_N(s,k)|+|\widehat z_N(s,-k)|\right)\\
&\ge\frac{e^{-N^2s/2}}2|\widehat z_N(s,k)|,
\qquad s=T-t.
\end{aligned}
\end{equation}
Thus~\eqref{eq:sharpness-fourier} reduces the desired divergence~\eqref{eq:sharpness-divergence} to finding frequencies $k_N\ne0$ such that
\begin{equation}\label{eq:sharpness-target}
|k_N|^{2+\alpha}
\|\widehat u_N(\cdot,k_N)\|_{L^\infty_{\omega,t}}
\longrightarrow\infty.
\end{equation}

\subsection{Counterexample for \texorpdfstring{$\sigma$}{sigma}}\label{subsec:sharpness-sigma}
Fix $1<\beta<1+\alpha$ and choose $\beta<\theta<1+\alpha$. For a sufficiently small fixed $\varepsilon_*>0$, take $T=a=1$ and
\begin{equation}\label{eq:sharpness-sigma-data}
\sigma(x)=\varepsilon_*\sum_{j=1}^\infty2^{-j\theta}\cos(2^jx),
\qquad \gamma=0,\qquad \phi(x)=e^{ix}.
\end{equation}
The H\"older--Besov characterization in Lemma~\ref{lem:frequency} gives $\|\sigma\|_{C_x^\theta}\le C_\theta\varepsilon_*$, since each frequency block contains only a bounded number of terms of this series. Hence $\sigma\in C_x^\theta\subset C_x^\beta$. Moreover,
\[
\widehat\sigma(2^j)=\frac{\varepsilon_*}{2}2^{-j\theta},
\qquad
(2^j)^{1+\alpha}|\widehat\sigma(2^j)|
=\frac{\varepsilon_*}{2}2^{j(1+\alpha-\theta)}\longrightarrow\infty,
\]
so~\eqref{eq:sharpness-fourier-decay} implies $\sigma\notin C_x^{1+\alpha}$. Choose $\varepsilon_*$ small enough that $\|\sigma\|_\infty\le1$, which gives $1-\sigma^2/2\ge1/2$.

For $N=2^J$, $J\ge1$, the BSPDE is
\begin{equation}\label{eq:sharpness-sigma-bspde}
\left\{
\begin{aligned}
&-du_N=\bigl[\partial_{xx}u_N+\sigma(x)\partial_xv_N\bigr]dt-v_N\,dW_t,\\
&u_N(1,x)=\cos(x+NW_1).
\end{aligned}
\right.
\end{equation}
All spatial derivatives of the terminal data are uniformly bounded. To prove~\eqref{eq:sharpness-target} at $k_N=N+1$, we use a short-time perturbation estimate to obtain a lower bound on $\widehat z_N(s,N+1)$.

Set $\tau_N=N^{-2}$ and let $P_s=e^{s\partial_{xx}}$ be the heat semigroup on $\mathbb T$. Its boundedness on $C_x^\theta$ and the smoothing estimate in~\cite[Lemma~A.7]{GubinelliImkellerPerkowski2015} give
\begin{equation}\label{eq:sharpness-heat-bounds}
\left\{
\begin{aligned}
&\|P_s q\|_{C_x^\theta}\le C\|q\|_{C_x^\theta},\\
&\|P_s q\|_{C_x^{1+\theta}}\le Cs^{-1/2}\|q\|_{C_x^\theta},
\end{aligned}
\right.
\qquad 0<s\le1.
\end{equation}
Subtract the free solution and set $r_N(s,x):=z_N(s,x)-e^{-s}e^{ix}$. Duhamel's formula gives
\[
r_N(s)=iN\int_0^sP_{s-\ell}
\bigl[\sigma\,\partial_x(e^{-\ell}e^{ix}+r_N(\ell))\bigr]\,d\ell.
\]
By~\eqref{eq:sharpness-heat-bounds} and the H\"older product estimate,
\[
\sup_{s\le\tau_N}\|r_N(s)\|_{C_x^{1+\theta}}
\le C\varepsilon_*N\sqrt{\tau_N}
\left(1+\sup_{s\le\tau_N}\|r_N(s)\|_{C_x^{1+\theta}}\right).
\]
The supremum is finite since $1+\theta<2+\alpha$. As $N\sqrt{\tau_N}=1$, choosing $\varepsilon_*$ sufficiently small and absorbing gives
\begin{equation}\label{eq:sharpness-sigma-perturbation}
\sup_{0\le s\le\tau_N}\|r_N(s)\|_{C_x^{1+\theta}}\le C\varepsilon_*.
\end{equation}
Here and below, constants may depend on $\alpha,\theta$ but are independent of $N$ and sufficiently small $\varepsilon_*$.

The decomposition $z_N(s,x)=e^{-s}e^{ix}+r_N(s,x)$ gives
$\partial_xz_N=ie^{-s}e^{ix}+\partial_xr_N$. Substituting this into~\eqref{eq:sharpness-deterministic}, with $a=1$ and $\gamma=0$, yields
\[
\partial_s z_N
=\partial_{xx}z_N-Ne^{-s}\sigma e^{ix}
+iN\sigma\,\partial_xr_N.
\]
Since $\widehat{\sigma e^{ix}}(N+1)=\frac{\varepsilon_*}{2}N^{-\theta}$, taking the Fourier coefficient at $N+1$ gives
\[
\frac{d}{ds}\widehat z_N(s,N+1)+(N+1)^2\widehat z_N(s,N+1)
=-\frac{\varepsilon_*}{2}N^{1-\theta}e^{-s}+\rho_N(s),
\]
where $\rho_N(s):=iN\widehat{\sigma\,\partial_xr_N(s)}(N+1)$. By~\eqref{eq:sharpness-fourier-decay} and~\eqref{eq:sharpness-sigma-perturbation},
\[
|\rho_N(s)|
\le CN^{1-\theta}\|\sigma\,\partial_xr_N(s)\|_{C_x^\theta}
\le C\varepsilon_*^2N^{1-\theta},
\qquad 0\le s\le\tau_N.
\]
Solving the scalar equation with zero initial datum therefore gives
\begin{equation}\label{eq:sharpness-sigma-remainder}
\begin{aligned}
&\left|\widehat z_N(s,N+1)
+\frac{\varepsilon_*N^{-\theta}}{2(N+2)}
\left(e^{-s}-e^{-(N+1)^2s}\right)\right|\\
&\qquad\le C\varepsilon_*^2N^{1-\theta}
\int_0^se^{-(N+1)^2(s-\ell)}\,d\ell
\le C\varepsilon_*^2N^{-1-\theta},
\qquad 0\le s\le\tau_N.
\end{aligned}
\end{equation}
For $s\in I_N:=[(2N^2)^{-1},N^{-2}]$, the difference of exponentials is at least $e^{-1/4}(1-e^{-1/2})$. Since $N+2\le2N$, fixing $\varepsilon_*$ sufficiently small yields
\begin{equation}\label{eq:sharpness-sigma-mode}
|\widehat z_N(s,N+1)|
\ge(c_0\varepsilon_*-C\varepsilon_*^2)N^{-1-\theta}
\ge c\varepsilon_*N^{-1-\theta},
\qquad s\in I_N.
\end{equation}
On the positive-length interval $t\in[1-N^{-2},1-(2N^2)^{-1}]$, we have $e^{-N^2(1-t)/2}\ge e^{-1/2}$. Combining~\eqref{eq:sharpness-fourier} and~\eqref{eq:sharpness-real-mode} therefore gives
\[
\|u_N\|_{L^\infty_{\omega,t}C_x^{2+\alpha}}
\ge c(N+1)^{2+\alpha}\varepsilon_*N^{-1-\theta}
\ge c\varepsilon_*N^{1+\alpha-\theta}\longrightarrow\infty.
\]
This proves sharpness for the fixed real-valued coefficient $\sigma$ in~\eqref{eq:sharpness-sigma-data}.

\subsection{Counterexample for \texorpdfstring{$\gamma$}{gamma}}\label{subsec:sharpness-gamma}
Choose $\beta,\theta,\varepsilon_*$ as in Section~\ref{subsec:sharpness-sigma}, and take $T=a=1$ with
\begin{equation}\label{eq:sharpness-gamma-data}
\sigma=0,\qquad
\gamma(x)=\varepsilon_*\sum_{j=1}^\infty2^{-j\theta}\cos(2^jx),
\qquad \phi(x)=1.
\end{equation}
The same frequency argument gives $\gamma\in C_x^\theta\subset C_x^\beta$ and $\gamma\notin C_x^{1+\alpha}$. For $N=2^J$, $J\ge1$, the BSPDE is
\begin{equation}\label{eq:sharpness-gamma-bspde}
\left\{
\begin{aligned}
&-du_N=\bigl[\partial_{xx}u_N+\gamma(x)v_N\bigr]dt-v_N\,dW_t,\\
&u_N(1,x)=\cos(NW_1).
\end{aligned}
\right.
\end{equation}
The terminal data are spatially constant and bounded by one, and the effective diffusion is one.

Set $r_N(s,x):=z_N(s,x)-1$ and $\tau_N=N^{-2}$. The same Duhamel and absorption argument, using the first bound in~\eqref{eq:sharpness-heat-bounds} and $N\tau_N=N^{-1}$, gives
\[
\sup_{0\le s\le\tau_N}\|r_N(s)\|_{C_x^\theta}\le C\varepsilon_*/N.
\]
At frequency $N$, equation~\eqref{eq:sharpness-deterministic} becomes
\[
\frac{d}{ds}\widehat z_N(s,N)+N^2\widehat z_N(s,N)
=\frac{i\varepsilon_*}{2}N^{1-\theta}+\rho_N(s),
\]
where, by~\eqref{eq:sharpness-fourier-decay},
\[
\rho_N(s):=iN\widehat{\gamma r_N(s)}(N),
\qquad
|\rho_N(s)|\le CN^{1-\theta}\|\gamma r_N(s)\|_{C_x^\theta}
\le C\varepsilon_*^2N^{-\theta}.
\]
Integrating this scalar equation with zero initial datum gives
\begin{equation}\label{eq:sharpness-gamma-mode-expansion}
\left|\widehat z_N(s,N)-\frac{i\varepsilon_*}{2}N^{-1-\theta}(1-e^{-N^2s})\right|
\le C\varepsilon_*^2N^{-2-\theta},
\qquad 0\le s\le\tau_N.
\end{equation}
For $s\in I_N=[(2N^2)^{-1},N^{-2}]$, the factor $1-e^{-N^2s}$ is at least $1-e^{-1/2}$. Hence, for sufficiently small fixed $\varepsilon_*$,
\[
|\widehat z_N(s,N)|
\ge(c_0\varepsilon_*-C\varepsilon_*^2/N)N^{-1-\theta}
\ge c\varepsilon_*N^{-1-\theta},
\qquad s\in I_N.
\]
As in Section~\ref{subsec:sharpness-sigma}, applying~\eqref{eq:sharpness-fourier} and~\eqref{eq:sharpness-real-mode} on the corresponding time interval yields
\[
\|u_N\|_{L^\infty_{\omega,t}C_x^{2+\alpha}}
\ge c\varepsilon_*N^{1+\alpha-\theta}\longrightarrow\infty.
\]
Thus the fixed real-valued coefficient $\gamma$ in~\eqref{eq:sharpness-gamma-data} also prevents a uniform estimate over these terminal data.

\section{Application to stochastic HJB equations}\label{sec:hjb}
We apply the linear Schauder theory to uniformly super-parabolic stochastic HJB equations with control-independent diffusion and a compact control set. We obtain spatial $C_x^{2+\alpha}$ estimates uniform in  $(\omega, t)$  under spatial H\"older assumptions. These equations arise from dynamic programming for control problems with random coefficients~\cite{Peng1992,MollaQiu2026}.

The probability space, still denoted by $(\Omega,\mathcal F,\mathbb P)$, carries two independent Brownian motions $W$ and $\overline W$ of dimensions $m_1$ and $m_2$, respectively. Let $(\mathcal F_s)$ be the augmented natural filtration of $W$, and $\mathcal P$ its predictable $\sigma$-algebra. For each initial time $t\in[0,T]$, the controller's filtration $\mathbb G^t=(\mathcal G_s^t)_{t\le s\le T}$ is the usual augmentation of
\[
\mathcal G_s^t:=\mathcal F_s\vee
\sigma(\overline W_r-\overline W_t:t\le r\le s),
\qquad t\le s\le T.
\]
In particular, $\mathcal G_t^t=\mathcal F_t$ up to completion.

Let $\mathcal A$ be a nonempty compact metric control space. Define the class of admissible controls by
\[
\mathcal U_t:=
\left\{\theta:\Omega\times[t,T]\to\mathcal A:
\theta\text{ is }\mathbb G^t\text{-predictable}\right\}.
\]
For $(t,x)\in[0,T]\times\mathbb T^n$ and $\theta\in\mathcal U_t$, the state process $X=X^{t,x;\theta}$ on $\mathbb T^n$ satisfies
\[
\left\{
\begin{aligned}
dX_s&=b^{\theta_s}(s,X_s)\,ds+\sigma(s,X_s)\,dW_s
+\overline\sigma(s,X_s)\,d\overline W_s,\qquad t\le s\le T,\\
X_t&=x.
\end{aligned}
\right.
\]
With running cost $f^\theta$ and terminal cost $g$, the control problem is to minimize the cost functional $J$, with value function $V$, defined by
\[
\begin{aligned}
J(t,x;\theta)&:=\E\left[
\int_t^T f^{\theta_s}(s,X_s^{t,x;\theta})\,ds
+g(X_T^{t,x;\theta})\,\middle|\,\mathcal F_t\right],\\
V(t,x)&:=\operatorname*{ess\,inf}_{\theta\in\mathcal U_t}J(t,x;\theta).
\end{aligned}
\]
Hence $J(t,x;\theta)$ and $V(t,x)$ are $\mathcal F_t$-measurable, and the value function is adapted to the filtration generated by $W$.

The drift $b^\theta$ and running cost $f^\theta$ are jointly $\mathcal P\otimes\mathcal B(\mathbb T^n)\otimes\mathcal B(\mathcal A)$-measurable and continuous in $\theta\in\mathcal A$. The diffusion coefficients $(\sigma,\overline\sigma):\Omega\times[0,T]\times\mathbb T^n\to\R^{n\times m_1}\times\R^{n\times m_2}$ are predictable and independent of the control, and the terminal cost $g$ is $\mathcal F_T\otimes\mathcal B(\mathbb T^n)$-measurable. The control enters only through the drift and running cost, so the Hamiltonian is
\begin{equation}\label{eq:hamiltonian}
H(t,x,p):=\inf_{\theta\in\mathcal A}
\{b^\theta(t,x)\cdot p+f^\theta(t,x)\},
\qquad p\in\R^n.
\end{equation}
When $V$ is sufficiently regular, dynamic programming and the It\^o--Wentzell formula give the semilinear BSPDE~\cite{Peng1992}
\begin{equation}\label{eq:hjb}
-du=[a:D^2u+\sigma:Dv+H(t,x,Du)]dt-v\cdot dW_t,
\qquad u(T)=g,
\end{equation}
where $u=V$ and
\[
a=\tfrac12(\sigma\sigma^\top+\overline\sigma\overline\sigma^\top),
\qquad
A=a-\tfrac12\sigma\sigma^\top=\tfrac12\overline\sigma\overline\sigma^\top.
\]
Uniform nondegeneracy of the independent noise gives $A\ge\lambda I$.

The analytic result below uses the following spatial H\"older assumptions; no time continuity is required.
\begin{assumption}\label{ass:hjb}
\begin{enumerate}
\renewcommand{\labelenumi}{(\roman{enumi})}
\item The coefficients $a,\sigma$ are $\mathcal P\otimes\mathcal B(\mathbb T^n)$-measurable, with $a=a^\top$. There are deterministic constants $K<\infty$ and $\lambda>0$ such that
\[
\|a\|_{L^\infty_{\omega,t}C_x^\alpha}
+\|\sigma\|_{L^\infty_{\omega,t}C_x^{1+\alpha}}\le K,
\qquad a-\tfrac12\sigma\sigma^\top\ge\lambda I.
\]
\item The control space $\mathcal A$ is a nonempty compact metric space. The coefficients $b^\theta,f^\theta$ are jointly $\mathcal P\otimes\mathcal B(\mathbb T^n)\otimes\mathcal B(\mathcal A)$-measurable and continuous in $\theta$. There is a deterministic constant $L<\infty$ such that
\begin{equation}\label{eq:hjb-coeff}
\left\{
\begin{aligned}
&\sup_{\theta\in\mathcal A}\|b^\theta\|_{L^\infty_{\omega,t}C_x^\alpha}\le L,\\
&\sup_{\theta\in\mathcal A}\|f^\theta\|_{L^\infty_{\omega,t}C_x^\alpha}<\infty.
\end{aligned}
\right.
\end{equation}
\item The terminal datum $g$ is $\mathcal F_T\otimes\mathcal B(\mathbb T^n)$-measurable and belongs to $L^\infty_\omega C_x^{2+\alpha}$.
\end{enumerate}
\end{assumption}

Sobolev strong and Schauder solutions are understood as in Definition~\ref{def:energy}, with the drift replaced by that of~\eqref{eq:hjb}. The theorem below gives Schauder regularity under Assumption~\ref{ass:hjb}, which requires only $C_x^\alpha$ regularity of the drift and running cost, uniformly over the controls. Boundedness of the drift makes $H$ globally Lipschitz in $p$. The estimate requires neither differentiability of $H$ nor spatial regularity of a minimizing feedback.

\begin{theorem}\label{thm:hjb}
Under Assumption~\ref{ass:hjb}, equation~\eqref{eq:hjb} has a unique Sobolev strong solution $(u,v)$, which satisfies
\begin{equation}\label{eq:hjb-estimate}
\|u\|_{L^\infty_{\omega,t}C_x^{2+\alpha}}
\le C\Bigl(\|g\|_{L^\infty_\omega C_x^{2+\alpha}}
+\sup_{\theta\in\mathcal A}\|f^\theta\|_{L^\infty_{\omega,t}C_x^\alpha}\Bigr),
\end{equation}
where $C=C(\lambda,K,L)$. In particular, $(u,v)$ is a Schauder solution.
\end{theorem}

We first record the Hamiltonian bounds needed for Sobolev solvability and spatial regularity. Suppressing $\omega$ in the notation, define $\tilde H_t:H_x^1\to L_x^2$ by
\[
\tilde H_t(w)(x):=H(t,x,Dw(x)).
\]
Continuity in $\theta$ allows the infimum in~\eqref{eq:hamiltonian} to be taken over a fixed countable dense subset of $\mathcal A$, so $H$ is predictable. For families $(A_\theta)$ and $(B_\theta)$ with finite infima, we have
\[
\left|\inf_\theta A_\theta-\inf_\theta B_\theta\right|
\le\sup_\theta|A_\theta-B_\theta|.
\]
Applying this inequality at a fixed spatial point and using~\eqref{eq:hamiltonian}, we obtain
\[
\left\{
\begin{aligned}
&|H(t,x,p)-H(t,x,q)|\le L|p-q|,\\
&|H(t,x,0)|\le\sup_{\theta\in\mathcal A}|f^\theta(t,x)|.
\end{aligned}
\right.
\]
Consequently,
\begin{equation}\label{eq:hjb-sobolev-lipschitz}
\left\{
\begin{aligned}
&\|\tilde{H}_t(w)-\tilde{H}_t(z)\|_{L_x^2}
\le L\|D(w-z)\|_{L_x^2},\\
&\|\tilde{H}_t(0)\|_{L_x^2}
\le C\sup_{\theta\in\mathcal A}\|f^\theta\|_{L^\infty_{\omega,t}C_x^\alpha}.
\end{aligned}
\right.
\end{equation}
These bounds give the following Sobolev result for~\eqref{eq:hjb}.

\begin{lemma}\label{lem:semilinear-sobolev}
Under Assumption~\ref{ass:hjb}(i), (ii), for every $\mathcal F_T$-measurable $g\in L^2_\omega H_x^1$, equation~\eqref{eq:hjb} has a unique Sobolev strong solution $(u,v)$. For every deterministic $t\in[0,T]$,
\begin{equation}\label{eq:semilinear-unconditional}
\begin{aligned}
&\E\sup_{t\le s\le T}\|u_s\|_{H_x^1}^2
+\E\int_t^T\bigl(\|u_s\|_{H_x^2}^2+\|v_s\|_{H_x^1}^2\bigr)\,ds\\
&\le C\E\left[\|g\|_{H_x^1}^2
+\int_t^T\|H(s,\cdot,0)\|_{L_x^2}^2\,ds\right],
\end{aligned}
\end{equation}
where $C=C(\lambda,K,L)$ is independent of $t$.
\end{lemma}

\begin{proof}
For every $\varepsilon>0$, Sobolev interpolation and~\eqref{eq:hjb-sobolev-lipschitz} give
\[
\|\tilde H_t(w)-\tilde H_t(z)\|_{L_x^2}
\le\varepsilon\|w-z\|_{H_x^2}
+C_{\varepsilon,L}\|w-z\|_{L_x^2},
\qquad w,z\in H_x^2.
\]
Thus $\tilde H$ satisfies~\cite[Assumption~5.4]{DuQiuTang2012}, with $\tilde H_s(0)=H(s,\cdot,0)$. The proof of~\cite[Theorem~5.5]{DuQiuTang2012}, using the periodic linear estimates of Lemma~\ref{lem:periodic-existence}, gives the result.
\end{proof}

For $B\in\mathcal F_t$, multiplying $\tilde H$ and $g$ by $\mathbf1_B$ preserves the Lipschitz bound in~\eqref{eq:hjb-sobolev-lipschitz}. The conditioning and stopping-time argument for~\eqref{eq:conditional-sobolev} therefore also gives the conditional version of~\eqref{eq:semilinear-unconditional}.
\begin{corollary}\label{cor:semilinear-conditional}
Under the hypotheses of Lemma~\ref{lem:semilinear-sobolev}, the solution $(u,v)$ of~\eqref{eq:hjb} satisfies, for every stopping time $\tau\le T$,
\begin{equation}\label{eq:semilinear-conditional}
\begin{aligned}
&\|u_\tau\|_{H_x^1}^2
+\E_\tau\int_\tau^T\bigl(\|u_s\|_{H_x^2}^2+\|v_s\|_{H_x^1}^2\bigr)\,ds\\
&\le C\E_\tau\left[\|g\|_{H_x^1}^2
+\int_\tau^T\|H(s,\cdot,0)\|_{L_x^2}^2\,ds\right],
\end{aligned}
\end{equation}
where $C=C(\lambda,K,L)$.
\end{corollary}

For the approximation argument, use the periodic mollifier $\varrho_\epsilon$ from Section~\ref{sec:closure} and set
\[
\tilde H^\epsilon_t(w):=\varrho_\epsilon*\tilde H_t(w),\qquad \epsilon>0.
\]
Convolution preserves measurability. Its contraction property on $L_x^2$ and~\eqref{eq:hjb-sobolev-lipschitz} give, for $w,z\in H_x^1$,
\begin{equation}\label{eq:semilinear-lipschitz}
\begin{aligned}
\|\tilde H^\epsilon_t(w)-\tilde H^\epsilon_t(z)\|_{L_x^2}
&\le\|\tilde H_t(w)-\tilde H_t(z)\|_{L_x^2}\\
&\le L\|D(w-z)\|_{L_x^2}
\le L\|w-z\|_{H_x^1}.
\end{aligned}
\end{equation}
Similarly,
\[
\|\tilde H^\epsilon_t(0)\|_{L_x^2}\le\|H(t,\cdot,0)\|_{L_x^2},
\qquad \tilde H^\epsilon(0)\in L^2_{\omega,t,x}.
\]
Consider the regularized equation
\begin{equation}\label{eq:hjb-regularized}
-du^\epsilon=[a:D^2u^\epsilon+\sigma:Dv^\epsilon
+\tilde H^\epsilon_t(u^\epsilon)]dt-v^\epsilon\cdot dW_t,
\qquad u^\epsilon(T)=g.
\end{equation}
The proofs of Lemma~\ref{lem:semilinear-sobolev} and Corollary~\ref{cor:semilinear-conditional} therefore apply to~\eqref{eq:hjb-regularized}, giving a unique Sobolev strong solution $(u^\epsilon,v^\epsilon)$ and the corresponding estimates with constants independent of $\epsilon$.
\begin{corollary}\label{cor:semilinear-stability}
Under the hypotheses of Lemma~\ref{lem:semilinear-sobolev}, let $(u,v)$ solve~\eqref{eq:hjb} and $(u^\epsilon,v^\epsilon)$ solve~\eqref{eq:hjb-regularized}, with the same terminal datum $g$. Then
\begin{equation}\label{eq:semilinear-stability}
\begin{aligned}
&\|u^\epsilon-u\|_{L^2_\omega C_tH_x^1}
+\|u^\epsilon-u\|_{L^2_{\omega,t}H_x^2}
+\|v^\epsilon-v\|_{L^2_{\omega,t}H_x^1}\\
&\le C\|\tilde H^\epsilon(u)-\tilde H(u)\|_{L^2_{\omega,t,x}},
\end{aligned}
\end{equation}
where $C=C(\lambda,K,L)$ is independent of $\epsilon$.
\end{corollary}

\begin{proof}
The difference equation for $(u^\epsilon-u,v^\epsilon-v)$ has zero terminal datum and nonlinear term
\[
w\longmapsto\tilde H^\epsilon_s(u_s+w)-\tilde H_s(u_s).
\]
By~\eqref{eq:semilinear-lipschitz}, this map has Lipschitz constant at most $L$, and its value at zero is the forcing term $\tilde H^\epsilon(u)-\tilde H(u)$. The proof of Lemma~\ref{lem:semilinear-sobolev}, applied to this difference equation at $t=0$, gives~\eqref{eq:semilinear-stability}.
\end{proof}

For the spatial H\"older bound, the inequality for infima applied at two spatial points gives, for $w\in C_x^{1+\alpha}$,
\[
\begin{aligned}
|H(t,x,Dw(x))-H(t,y,Dw(y))|
\le{}&L|Dw(x)-Dw(y)|\\
&+\Bigl(L\|Dw\|_\infty+\sup_{\theta\in\mathcal A}\|f^\theta\|_{L^\infty_{\omega,t}C_x^\alpha}\Bigr)\rho(x-y)^\alpha.
\end{aligned}
\]
Together with the bound on $\|\tilde H_t(w)\|_\infty$, this yields
\begin{equation}\label{eq:hjb-composition}
\|\tilde{H}_t(w)\|_{C_x^\alpha}
\le C\Bigl(\sup_{\theta\in\mathcal A}\|f^\theta\|_{L^\infty_{\omega,t}C_x^\alpha}+L\|w\|_{C_x^{1+\alpha}}\Bigr).
\end{equation}

\begin{proof}[Proof of Theorem~\ref{thm:hjb}]
Lemma~\ref{lem:semilinear-sobolev} gives a unique Sobolev strong solution $(u,v)$ of~\eqref{eq:hjb}. Set
\[
D_H:=\|g\|_{L^\infty_\omega C_x^{2+\alpha}}+\sup_{\theta\in\mathcal A}\|f^\theta\|_{L^\infty_{\omega,t}C_x^\alpha}.
\]

\smallskip\noindent\emph{Step 1: Uniform Sobolev bound.}
For the solutions $(u^\epsilon,v^\epsilon)$ of~\eqref{eq:hjb-regularized}, the conditional estimate from Corollary~\ref{cor:semilinear-conditional} and the bound on $\tilde H^\epsilon(0)$ above give, for every stopping time $\tau$,
\[
\|u^\epsilon_\tau\|_{H_x^1}^2\le C D_H^2.
\]
It follows that
\begin{equation}\label{eq:hjb-uniform-h1}
\sup_{\epsilon>0}\|u^\epsilon\|_{L^\infty_{\omega,t}H_x^1}
\le C D_H.
\end{equation}

\smallskip\noindent\emph{Step 2: Uniform Schauder estimate.}
For fixed $\epsilon>0$, convolution with $\varrho_\epsilon$ is bounded from $L_x^2$ to $C_x^\alpha$. Together with~\eqref{eq:hjb-sobolev-lipschitz} and~\eqref{eq:hjb-uniform-h1}, this gives
\[
\|\tilde{H}^\epsilon(u^\epsilon)\|_{L^\infty_{\omega,t}C_x^\alpha}
\le C_\epsilon\left(\sup_{\theta\in\mathcal A}\|f^\theta\|_{L^\infty_{\omega,t}C_x^\alpha}
+L\|u^\epsilon\|_{L^\infty_{\omega,t}H_x^1}\right)<\infty.
\]
Thus Theorem~\ref{thm:linear}, applied with forcing $\tilde H^\epsilon(u^\epsilon)$ and $b=c=\gamma=0$, gives $u^\epsilon\in L^\infty_{\omega,t}C_x^{2+\alpha}$.

To obtain a bound independent of $\epsilon$, use the contraction property of convolution with $\varrho_\epsilon$ on $C_x^\alpha$ and~\eqref{eq:hjb-composition}:
\[
\begin{aligned}
\|\tilde{H}^\epsilon(u^\epsilon)\|_{L^\infty_{\omega,t}C_x^\alpha}
&\le\|\tilde{H}(u^\epsilon)\|_{L^\infty_{\omega,t}C_x^\alpha}\\
&\le C\left(\sup_{\theta\in\mathcal A}\|f^\theta\|_{L^\infty_{\omega,t}C_x^\alpha}
+L\|u^\epsilon\|_{L^\infty_{\omega,t}C_x^{1+\alpha}}\right).
\end{aligned}
\]
Theorem~\ref{thm:linear} therefore gives
\[
\|u^\epsilon\|_{L^\infty_{\omega,t}C_x^{2+\alpha}}
\le C\left(D_H+L\|u^\epsilon\|_{L^\infty_{\omega,t}C_x^{1+\alpha}}\right).
\]
By~\eqref{eq:holder-sobolev-interpolation}, for every $\delta>0$,
\[
\|u^\epsilon\|_{L^\infty_{\omega,t}C_x^{1+\alpha}}
\le\delta\|u^\epsilon\|_{L^\infty_{\omega,t}C_x^{2+\alpha}}
+C_\delta\|u^\epsilon\|_{L^\infty_{\omega,t}H_x^1}.
\]
Choose $\delta$ sufficiently small to absorb the $C_x^{2+\alpha}$ term and apply~\eqref{eq:hjb-uniform-h1}. This yields, uniformly in $\epsilon$,
\begin{equation}\label{eq:hjb-uniform-regularization}
\|u^\epsilon\|_{L^\infty_{\omega,t}C_x^{2+\alpha}}\le C D_H.
\end{equation}

\smallskip\noindent\emph{Step 3: Passage to the limit.}
By~\eqref{eq:hjb-sobolev-lipschitz}, the forcing $\tilde H(u)$ belongs to $L^2_{\omega,t,x}$. Spatial convolution therefore gives
\[
\|\tilde{H}^\epsilon(u)-\tilde{H}(u)\|_{L^2_{\omega,t,x}}
=\|\varrho_\epsilon*\tilde{H}(u)-\tilde{H}(u)\|_{L^2_{\omega,t,x}}
\longrightarrow0.
\]
Corollary~\ref{cor:semilinear-stability} now yields
\[
u^\epsilon\longrightarrow u\quad\text{in }L^2_{\omega,t}H_x^2.
\]
The lower semicontinuity of the spatial H\"older norm established in Lemma~\ref{lem:energy-limit}, together with~\eqref{eq:hjb-uniform-regularization}, gives
\[
\|u\|_{L^\infty_{\omega,t}C_x^{2+\alpha}}
\le\liminf_{\epsilon\downarrow0}\|u^\epsilon\|_{L^\infty_{\omega,t}C_x^{2+\alpha}}
\le C D_H,
\]
which proves~\eqref{eq:hjb-estimate}.
\end{proof}

\appendix
\section{Frequency and conditional energy estimates}\label{app:windows-frequency}
We first prove the Littlewood--Paley kernel bounds and the Bernstein and interpolation inequalities. We then establish convolution and scale-comparison estimates for the conditional energy, including its extension to measurable random centers. The final two subsections prove the localized heat-kernel bounds and the commutator and lower-order estimates used in the perturbation argument.

Throughout Appendix~\ref{app:windows-frequency}, we use the notation of Section~\ref{subsec:frequency-strategy}. For the fixed smooth cutoff $\chi$,
\[
r_j=2^{-j},\qquad
\varphi(\xi)=\chi(\xi)-\chi(2\xi),\qquad
m_0=\chi,\quad m_j(\xi)=\varphi(r_j\xi)\quad(j\ge1).
\]
The function $\psi_j$ is the convolution kernel of $\Delta_j$ on $\mathbb T^n$, and $\rho(y)=\operatorname{dist}_{\mathbb T^n}(y,0)$.

\subsection{Littlewood--Paley kernel estimates}\label{app:block-kernels}
We recall the standard decay estimates for the Euclidean kernels
\[
k_j(z):=\frac1{(2\pi)^n}\int_{\R^n}e^{iz\cdot\xi}m_j(\xi)\,d\xi,
\qquad z\in\R^n.
\]
Since $\chi,\varphi\in C_c^\infty(\R^n)$, their inverse Fourier transforms are Schwartz functions. Scaling therefore gives, for every multi-index $\beta$ and $N\ge0$,
\[
|D^\beta k_j(z)|
\le C_{\beta,N}r_j^{-n-|\beta|}(1+|z|/r_j)^{-N},
\qquad z\in\R^n,\quad j\ge0,
\]
where $r_0=1$. We prove Lemma~\ref{lem:block-kernels} by periodizing these Euclidean kernels.

\begin{proof}[Proof of Lemma~\ref{lem:block-kernels}]
By the Poisson summation formula~\cite[Lemma~6]{GubinelliPerkowski2015},
\[
\psi_j(y)=\sum_{q\in\mathbb Z^n}k_j(y+2\pi q).
\]
Fix $N>n+2$. For $0\le\nu\le2$, unfolding the integral over the torus, using $\rho(y)\le|y+2\pi q|$, and scaling $z=r_jw$ give
\[
\begin{aligned}
\int_{\mathbb T^n}\rho(y)^\nu|\psi_j(y)|\,dy
&\le\int_{\R^n}|z|^\nu|k_j(z)|\,dz\\
&\le Cr_j^\nu\int_{\R^n}|w|^\nu(1+|w|)^{-N}\,dw
\le C_\nu r_j^\nu.
\end{aligned}
\]
The last integral is finite since $N>n+\nu$, and $\nu=0$ gives the $L^1$ bound.

Finally, the Fourier coefficient of $\psi_j$ at frequency zero gives
\[
\int_{\mathbb T^n}\psi_j(y)\,dy=m_j(0)
=\begin{cases}
1,&j=0,\\
0,&j\ge1,
\end{cases}
\]
since $\chi(0)=1$ and $\varphi(0)=0$.
\end{proof}

\subsection{Bernstein and H\"older interpolation inequalities}\label{app:frequency}
Fix $C_*\ge2$ and define the periodic convolution kernel of the multiplier $\chi(r_k\xi/C_*)$ by
\begin{equation}\label{eq:reproducing-kernel}
\Psi_k(y):=\frac1{(2\pi)^n}\sum_{\xi\in\mathbb Z^n}
\chi(r_k\xi/C_*)e^{i\xi\cdot y},
\qquad y\in\mathbb T^n,\quad k\ge0.
\end{equation}
The multiplier $\chi(r_k\xi/C_*)$ equals one for $|\xi|\le C_*r_k^{-1}$. Thus every $h_k$ with Fourier support in $\{|\xi|\le C_*r_k^{-1}\}$ satisfies
\begin{equation}\label{eq:reproducing-kernel-identity}
h_k=\Psi_k*h_k,\qquad
D^\ell h_k=(D^\ell\Psi_k)*h_k,\quad \ell=1,2.
\end{equation}

\begin{proof}[Proof of Lemma~\ref{lem:bernstein}]
Define the Euclidean kernels
\[
\Theta(z):=\frac1{(2\pi)^n}\int_{\R^n}e^{iz\cdot\xi}\chi(\xi/C_*)\,d\xi,
\qquad
\Theta_k(z):=r_k^{-n}\Theta(z/r_k).
\]
Then $\Theta_k$ is the inverse Fourier transform of $\chi(r_k\xi/C_*)$. Since $\chi(\cdot/C_*)\in C_c^\infty(\R^n)$, $\Theta$ is a Schwartz function. Scaling gives, for $\ell=0,1,2$,
\[
|D^\ell\Theta_k(z)|
\le Cr_k^{-n-\ell}(1+|z|^2/r_k^2)^{-(n+1)/2}
\le Cr_k^{-\ell}w_{r_k}(z),
\qquad z\in\R^n,
\]
where $w_{r_k}$ is the Euclidean weight defined in~\eqref{eq:old-2-3} and $C=C(C_*)$ is independent of $k$.

Poisson summation and the Schwartz decay of $\Theta_k$ give
\[
D^\ell\Psi_k(y)=\sum_{q\in\mathbb Z^n}D^\ell\Theta_k(y+2\pi q),
\qquad\ell=0,1,2.
\]
Hence,
\begin{equation}\label{eq:reproducing-kernel-bounds}
|D^\ell\Psi_k(y)|\le Cr_k^{-\ell}\kappa_{r_k}(y) \qquad\text{and}\qquad \|D^\ell\Psi_k\|_{L^1(\mathbb T^n)}\le Cr_k^{-\ell}.
\end{equation}
The reproducing identity~\eqref{eq:reproducing-kernel-identity} and Cauchy--Schwarz on $\mathbb T^n$ yield
\[
\begin{aligned}
|D^\ell h_k(x)|^2
&\le \|D^\ell\Psi_k\|_1
\int_{\mathbb T^n}|D^\ell\Psi_k(z)|\,|h_k(x-z)|^2\,dz\\
&\le Cr_k^{-2\ell}(\kappa_{r_k}*|h_k|^2)(x).
\end{aligned}
\]
This proves (i). Since $\kappa_{r_k}*|h_k|^2\le\|h_k\|_\infty^2$, (ii) follows.
\end{proof}

\begin{proof}[Proof of Lemma~\ref{lem:holder-interpolation}]
Write $h=\sum_{j\ge0}h_j$, where $h_j:=\Delta_jh$, and split the sum at an integer $J\ge1$. Bernstein's inequality and the H\"older--Besov characterization give
\[
\|D^2h_j\|_\infty\le Cr_j^{-2}\|h_j\|_\infty,
\qquad
\|h_j\|_{C^{1+\alpha}}\le Cr_j^{-1-\alpha}\|h_j\|_\infty.
\]
For the second estimate, we use the finite overlap of the block supports and the uniform $L^1$ bounds for their kernels.

For the high frequencies, Lemma~\ref{lem:frequency} with $\beta=2+\alpha$ gives
\[
\|h_j\|_\infty\le Cr_j^{2+\alpha}\|h\|_{C^{2+\alpha}}.
\]
Consequently,
\[
\sum_{j\ge J}\bigl(\|D^2h_j\|_\infty+\|h_j\|_{C^{1+\alpha}}\bigr)\le Cr_J^\alpha\|h\|_{C^{2+\alpha}}.
\]
These summable bounds also justify convergence of the block series to $h$ in $C^2$ and $C^{1+\alpha}$.

For the finitely many low-frequency blocks, Lemma~\ref{lem:block-kernels} gives $\|h_j\|_\infty\le C\|h\|_\infty$. Hence
\[
\sum_{j<J}\bigl(\|D^2h_j\|_\infty+\|h_j\|_{C^{1+\alpha}}\bigr)
\le C_J\|h\|_\infty.
\]
Alternatively,~\eqref{eq:finite-fourier-block} gives $\|h_j\|_\infty\le C_j\|h\|_{L_x^2}$, and therefore
\[
\sum_{j<J}\|h_j\|_{C^{1+\alpha}}
\le C_J\|h\|_{L_x^2}
\le C_J\|h\|_{H_x^1}.
\]
Given $\delta>0$, choose $J=J(\delta)$ so that $Cr_J^\alpha\le\delta$. Combining these bounds proves~\eqref{eq:old-2-16} and~\eqref{eq:holder-sobolev-interpolation}.
\end{proof}

\subsection{Conditional energy: convolution and scale comparison}\label{app:weights}
$\mathcal E_r(h)$ is invariant under spatial translations of $h$, since translating $h$ merely shifts the window center over which the supremum is taken. Minkowski's inequality then gives the following convolution estimate.

\begin{lemma}\label{lem:calculus}
Let $\mathcal E_r$ be defined by~\eqref{eq:old-2-6}. For every $0<r\le1$ and deterministic $L\in L^1(\mathbb T^n)$,
\begin{equation}\label{eq:old-2-17}
\mathcal E_r(L*h)
\le
\|L\|_{L^1}\mathcal E_r(h).
\end{equation}
\end{lemma}
\begin{proof}
Fix $\tau,x$ and set $\zeta=(\tau+r^2)\wedge T$. Minkowski's inequality in the weighted spatial space gives
\[
\|(L*h_s)(x+\cdot)\|_{\kappa_r}
\le\int|L(z)|\,\|h_s(x-z+\cdot)\|_{\kappa_r}\,dz,
\]
where $\kappa_r$ is the periodic weight defined in~\eqref{eq:periodic-window-weight}. Cauchy--Schwarz and conditional Fubini give
\[
\begin{aligned}
\E_\tau\int_\tau^\zeta\|(L*h_s)(x+\cdot)\|_{\kappa_r}^2\,ds
&\le\E_\tau\int_\tau^\zeta
\left(\int|L(z)|\,\|h_s(x-z+\cdot)\|_{\kappa_r}\,dz\right)^2ds\\
&\le\|L\|_1\int|L(z)|\,\E_\tau\int_\tau^\zeta
\|h_s(x-z+\cdot)\|_{\kappa_r}^2\,ds\,dz\\
&\le\|L\|_1^2\mathcal E_r(h)^2.
\end{aligned}
\]
Taking the supremum over the centers and stopping times proves~\eqref{eq:old-2-17}.
\end{proof}

For $\beta>n$, define the Euclidean kernels
\begin{equation}\label{eq:euclidean-dyadic-kernels}
\eta_{k,\beta}(y):=2^{kn}(1+2^k|y|)^{-\beta},
\qquad k\in\mathbb N_0,\quad y\in\R^n.
\end{equation}

\begin{lemma}[Diening--H\"ast\"o--Roudenko {\cite[Lemma~A.3]{DieningHastoRoudenko2009}}]\label{lem:euclidean-weight-convolution}
Let $\beta>n$. There exists a constant $C=C(n,\beta)\ge1$ such that, for all $k,\ell\in\mathbb N_0$ and $y\in\R^n$,
\[
C^{-1}\eta_{\min\{k,\ell\},\beta}(y)
\le(\eta_{k,\beta}*\eta_{\ell,\beta})(y)
\le C\eta_{\min\{k,\ell\},\beta}(y),
\]
where convolution is taken on $\R^n$.
\end{lemma}

Comparing $w_r$ with these dyadic kernels and periodizing gives the following estimate.
\begin{lemma}\label{lem:window-auxiliary}
There is a constant $C=C(n)\ge1$ such that, for all $r,R\in(0,1]$ and $y\in\mathbb T^n$,
\begin{equation}\label{eq:old-2-11}
C^{-1}\kappa_{\max\{r,R\}}(y)
\le(\kappa_r*\kappa_R)(y)
\le C\kappa_{\max\{r,R\}}(y).
\end{equation}
\end{lemma}
\begin{proof}
Choose $k,\ell\in\mathbb N_0$ so that $r\in(2^{-k-1},2^{-k}]$ and $R\in(2^{-\ell-1},2^{-\ell}]$.
By direct computation, we have
\[
w_r\asymp\eta_{k,n+1}\qquad\text{and}\qquad w_R\asymp\eta_{\ell,n+1},
\]
where $w_r,w_R$ are defined in~\eqref{eq:old-2-3} and $\eta_{k,\beta}$ in~\eqref{eq:euclidean-dyadic-kernels}. Lemma~\ref{lem:euclidean-weight-convolution} with $\beta=n+1$ therefore yields
\[
w_r*w_R\asymp\eta_{\min\{k,\ell\},n+1}
\asymp w_{\max\{r,R\}}
\quad\text{on }\R^n.
\]
By Tonelli's theorem,
\[
\begin{aligned}
(\kappa_r*\kappa_R)(y)
&=\int_{\mathbb T^n}\sum_{p,q\in\mathbb Z^n}
w_r(y-z+2\pi p)w_R(z+2\pi q)\,dz\\
&=\sum_{\ell\in\mathbb Z^n}\int_{\R^n}
w_r(y-z+2\pi\ell)w_R(z)\,dz\\
&=\sum_{\ell\in\mathbb Z^n}(w_r*w_R)(y+2\pi\ell)
\asymp\kappa_{\max\{r,R\}}(y).
\end{aligned}
\]
This proves~\eqref{eq:old-2-11}.
\end{proof}

For $0<r\le R\le1$, averaging the weighted energy at scale $r$ with $\kappa_R$ gives the weighted energy at scale $R$, up to uniform constants. Indeed, Lemma~\ref{lem:window-auxiliary} and Fubini's theorem give, for every $h\in L^2(\mathbb T^n)$ and $x\in\mathbb T^n$,
\[
\begin{aligned}
\int_{\mathbb T^n}\kappa_R(z)\|h(x+z+\cdot)\|_{\kappa_r}^2\,dz
&=\int_{\mathbb T^n}(\kappa_R*\kappa_r)(y)|h(x+y)|^2\,dy\\
&\asymp\|h(x+\cdot)\|_{\kappa_R}^2.
\end{aligned}
\]
This comparison is used in Lemma~\ref{lem:window-scales} to control the energy at scale $R$ by an average of energies at scale $r$.

\begin{proof}[Proof of Lemma~\ref{lem:window-scales}]
For $\ell=1,2$, apply Lemma~\ref{lem:calculus} to the convolution identity~\eqref{eq:reproducing-kernel-identity}, with kernel $L=D^\ell\Psi_k$ and localization scale $r_j$. The kernel bound~\eqref{eq:reproducing-kernel-bounds} gives
\[
\begin{aligned}
\mathcal E_{r_j}(D^\ell h_k)
&=\mathcal E_{r_j}\bigl((D^\ell\Psi_k)*h_k\bigr)\\
&\le\|D^\ell\Psi_k\|_{L^1(\mathbb T^n)}\mathcal E_{r_j}(h_k)
\le Cr_k^{-\ell}\mathcal E_{r_j}(h_k),
\end{aligned}
\]
which proves~\eqref{eq:old-2-18}.

For the scale comparison, put $r=r_j$, $R=r_k$, and $\zeta_r=(\tau+r^2)\wedge T$.
If $r\le R$, the local Bernstein bound~\eqref{eq:old-2-13} gives
\[
|h_k(s,y)|^2\le C(\kappa_R*|h_k(s)|^2)(y),
\qquad y\in\mathbb T^n.
\]
Combining this with the convolution comparison~\eqref{eq:old-2-11} gives
\[
\|h_k(s)\|_{\kappa_{r,x}}^2
\le C\int(\kappa_r*\kappa_R)(x-y)|h_k(s,y)|^2dy
\le C\|h_k(s)\|_{\kappa_{R,x}}^2.
\]
The inclusion $[\tau,\zeta_r]\subset[\tau,(\tau+R^2)\wedge T]$ then gives
$\mathcal E_r(h_k)\le C\mathcal E_R(h_k)$.

If $r>R$, use $\kappa_r\le C\kappa_r*\kappa_R$ from~\eqref{eq:old-2-11} to obtain
\begin{equation}\label{eq:old-2-21}
\|h_k(s)\|_{\kappa_{r,x}}^2
\le C\int\kappa_r(z)\|h_k(s)\|_{\kappa_{R,x+z}}^2\,dz.
\end{equation}
Set $N:=\lceil r^2/R^2\rceil$, so that $N\le2(r/R)^2$, and define the stopping times
\[
\tau_\nu:=(\tau+\nu R^2)\wedge\zeta_r,
\qquad \nu=0,\ldots,N.
\]
Then $\tau_0=\tau$ and $\tau_N=\zeta_r$, and each interval $[\tau_\nu,\tau_{\nu+1}]$ has length at most $R^2$.
For every fixed $z$, the definition of $\mathcal E_R$ gives
\[
\E_{\tau_\nu}\int_{\tau_\nu}^{\tau_{\nu+1}}
\|h_k(s)\|_{\kappa_{R,x+z}}^2ds\le\mathcal E_R(h_k)^2.
\]
Conditioning at $\tau$, summing the intervals, and integrating~\eqref{eq:old-2-21} in $z$ yields
\[
\E_\tau\int_\tau^{\zeta_r}\|h_k(s)\|_{\kappa_{r,x}}^2ds
\le CN\mathcal E_R(h_k)^2.
\]
Taking square roots proves~\eqref{eq:old-2-19}.
\end{proof}

\subsection{Conditional energy at random centers}\label{app:window-calculus}
\begin{proof}[Proof of Lemma~\ref{lem:random-center-energy}]
Fix a stopping time $\tau\le T$ and an $\mathcal F_\tau$-measurable center $z$. Put $\zeta=(\tau+r^2)\wedge T$ and write
\[
J_{\tau,x}:=\int_\tau^\zeta\|h_s\|_{\kappa_{r,x}}^2\,ds,
\qquad x\in\mathbb T^n.
\]
By the definition of $\mathcal E_r(h)$, $\E_\tau J_{\tau,x}\le\mathcal E_r(h)^2$ for $(\mathbb P\otimes dx)$-almost every $(\omega,x)$. Fubini's theorem therefore gives a deterministic dense set $G_\tau\subset\mathbb T^n$ such that this estimate holds $\mathbb P$-a.s.\ for each fixed $x\in G_\tau$.

Choose $x_*\in G_\tau$ and set $\Omega_\tau:=\{J_{\tau,x_*}<\infty\}$. Taking expectations gives
\[
\E J_{\tau,x_*}\le\mathcal E_r(h)^2<\infty,
\]
so $\mathbb P(\Omega_\tau)=1$. Since $c_r:=\min_{y\in\mathbb T^n}\kappa_r(y)>0$, on $\Omega_\tau$ we have $\int_\tau^\zeta\|h_s\|_{L_x^2}^2\,ds\le c_r^{-1}J_{\tau,x_*}<\infty$ and
\[
|J_{\tau,x}-J_{\tau,x'}|
\le c_r^{-1}J_{\tau,x_*}
\|\kappa_r(\cdot-x)-\kappa_r(\cdot-x')\|_\infty
\longrightarrow0
\qquad\text{as }x'\to x,
\]
by uniform continuity of $\kappa_r$. Thus $x\mapsto J_{\tau,x}(\omega)$ is continuous for every $\omega\in\Omega_\tau$.

Since $G_\tau$ is dense and the torus is compact, for each $N\ge1$ we can choose finitely many points $x_{N,1},\ldots,x_{N,m_N}\in G_\tau$ whose balls of radius $1/N$ cover the torus. Define an $\mathcal F_\tau$-measurable partition $(A_{N,i})_{i=1}^{m_N}$ of the sample space $\Omega$ by
\[
\begin{aligned}
A_{N,i}:={}&\{\omega\in\Omega:\rho(z(\omega)-x_{N,i})<1/N\}\\
&\cap\bigcap_{\ell=1}^{i-1}\{\omega\in\Omega:\rho(z(\omega)-x_{N,\ell})\ge1/N\},
\qquad 1\le i\le m_N.
\end{aligned}
\]
Thus $A_{N,i}$ is the event that $i$ is the first index whose ball contains $z$. Since $z$ is $\mathcal F_\tau$-measurable, each $A_{N,i}$ belongs to $\mathcal F_\tau$. The covering property ensures that these events are disjoint and their union is $\Omega$.

We can therefore define
\[
z_N(\omega):=x_{N,i}\qquad\text{if }\omega\in A_{N,i}.
\]
This random center is finite-valued and $\mathcal F_\tau$-measurable, and satisfies $\rho(z_N-z)<1/N$. On $A_{N,i}$, we have $J_{\tau,z_N}=J_{\tau,x_{N,i}}$. Taking conditional expectation and using $A_{N,i}\in\mathcal F_\tau$ gives
\[
\E_\tau J_{\tau,z_N}
=\sum_{i=1}^{m_N}\mathbf1_{A_{N,i}}\E_\tau J_{\tau,x_{N,i}}
\le\mathcal E_r(h)^2
\qquad\text{a.s.}
\]
Here we used the fixed-center estimate for $x_{N,i}\in G_\tau$.

For every $\omega\in\Omega_\tau$, continuity of $x\mapsto J_{\tau,x}(\omega)$ gives $J_{\tau,z_N}(\omega)\to J_{\tau,z}(\omega)$. Since $\mathbb P(\Omega_\tau)=1$ and these random variables are nonnegative, conditional Fatou's lemma yields
\[
\E_\tau J_{\tau,z}
\le\liminf_{N\to\infty}\E_\tau J_{\tau,z_N}
\le\mathcal E_r(h)^2
\qquad\text{a.s.}
\]
This proves the lemma.
\end{proof}

\subsection{Localized heat-kernel estimates}\label{app:heat-kernel-proof}
\begin{proof}[Proof of Lemma~\ref{lem:heat-kernel-decay}]
We estimate the Fourier series directly by finite differences and summation by parts on $\mathbb Z^n$; see~\cite[Section~2]{RuzhanskyTurunen2007} for this discrete calculus.
Fix $t,s,j$, put $r:=r_j$ and $\eta:=(s-t)/r^2$, and let $\mu$ be a multi-index with $|\mu|=\ell\le1$. Define
\[
a_{j,\mu}(\xi):=(i\xi)^\mu\widetilde m_j(\xi)e^{-\xi^\top Q_{t,s}\xi},
\qquad \xi\in\R^n.
\]
On the support of $\widetilde m_j$,
\[
\frac14\le r|\xi|\le4,\qquad
\xi^\top Q_{t,s}\xi\ge\frac{\lambda\eta}{16},\qquad
\|Q_{t,s}\|\le\Lambda r^2\eta.
\]
For every multi-index $\gamma$, differentiation of the cutoff and the Gaussian factor gives, on this support,
\[
\begin{aligned}
\bigl|\partial_\xi^\gamma[(i\xi)^\mu\widetilde m_j(\xi)]\bigr|
&\le C_\gamma r^{|\gamma|-\ell},\\
\bigl|\partial_\xi^\gamma e^{-\xi^\top Q_{t,s}\xi}\bigr|
&\le C_\gamma r^{|\gamma|}(1+\eta)^{|\gamma|}e^{-\lambda\eta/16}.
\end{aligned}
\]
Indeed, each derivative of the cutoff contributes a factor $r$, while differentiating the Gaussian produces factors $Q_{t,s}\xi$ and $Q_{t,s}$, bounded by $Cr\eta$ and $Cr^2\eta$, respectively. The product rule therefore gives, for every multi-index $\beta$ and $\xi\in\R^n$,
\begin{equation}\label{eq:heat-multiplier-derivatives}
\begin{aligned}
|\partial_\xi^\beta a_{j,\mu}(\xi)|
&\le\sum_{\gamma\le\beta}\binom{\beta}{\gamma}
\bigl|\partial_\xi^\gamma[(i\xi)^\mu\widetilde m_j(\xi)]\bigr|
\bigl|\partial_\xi^{\beta-\gamma}e^{-\xi^\top Q_{t,s}\xi}\bigr|\\
&\le C_\beta r^{|\beta|-\ell}(1+\eta)^{|\beta|}e^{-\lambda\eta/16}\\
&\le C_\beta r^{|\beta|-\ell}e^{-c\eta},
\qquad c:=\lambda/32.
\end{aligned}
\end{equation}
The last inequality uses $(1+\eta)^{|\beta|}e^{-\lambda\eta/32}\le C_\beta$ for $\eta\ge0$.

For the $p$th coordinate vector $e_p$, set $\delta_p a(\xi):=a(\xi+e_p)-a(\xi)$. Repeated use of the fundamental theorem of calculus yields
\[
\delta_p^q a_{j,\mu}(\xi)
=\int_{[0,1]^q}\partial_{\xi_p}^q a_{j,\mu}
\bigl(\xi+(\theta_1+\cdots+\theta_q)e_p\bigr)\,d\theta,
\qquad q\ge1.
\]
The difference vanishes for $|\xi|>4/r+q$, so at most $C_qr^{-n}$ lattice points contribute. By~\eqref{eq:heat-multiplier-derivatives},
\[
\sum_{\xi\in\mathbb Z^n}|\delta_p^q a_{j,\mu}(\xi)|
\le C_qr^{q-n-\ell}e^{-c\eta},\qquad q\ge0,
\]
where $q=0$ denotes the multiplier itself. Choose the representative $y=(y_1,\ldots,y_n)\in[-\pi,\pi]^n$, where $y_p$ denotes its $p$th coordinate. Shifting indices in the finite Fourier sum gives
\[
(e^{-iy_p}-1)^qD^\mu K_j(t,s,y)
=\frac1{(2\pi)^n}\sum_{\xi\in\mathbb Z^n}
\delta_p^q a_{j,\mu}(\xi)e^{i\xi\cdot y}.
\]
Since $\rho(y)=|y|$, for $y\ne0$ we can choose $p$ such that $|y_p|\ge\rho(y)/\sqrt n$. Hence
\[
|e^{-iy_p}-1|=2|\sin(y_p/2)|\ge\frac{2\rho(y)}{\pi\sqrt n}.
\]
Combining the estimates with $q=0$ and $q\ge1$ gives
\begin{equation}\label{eq:periodic-heat-kernel-decay}
|D^\mu K_j(t,s,y)|
\le C_qr^{-n-\ell}e^{-c\eta}
\min\{1,(r/\rho(y))^q\}
\le C_qr^{-n-\ell}e^{-c\eta}(1+\rho(y)/r)^{-q}.
\end{equation}
At $y=0$, the final bound follows directly from $q=0$.

Finally, the nearest lattice translate in the definition of $\kappa_r$ gives
\[
\kappa_r(y)\ge c_nr^{-n}(1+\rho(y)^2/r^2)^{-(n+1)/2}
\ge c_nr^{-n}(1+\rho(y)/r)^{-(n+1)}.
\]
Take $q=n+3$ in~\eqref{eq:periodic-heat-kernel-decay}. For $0\le\nu\le2$,
\[
\rho(y)^\nu|D^\mu K_j(t,s,y)|
\le Cr^{\nu-\ell}e^{-c\eta}
r^{-n}(1+\rho(y)/r)^{-(n+1)}
\le Cr^{\nu-\ell}e^{-c\eta}\kappa_r(y).
\]
Summing over $|\mu|=\ell$ proves~\eqref{eq:old-4-9}.
\end{proof}

\subsection{Commutator and lower-order estimates}\label{app:commutator-proofs}
We first establish conditional energy estimates for products with bounded coefficients and commutators with spatially Lipschitz coefficients, using the convention $D^0v_k:=v_k$.

\begin{lemma}\label{lem:conditional-products}
Let $h\in L^\infty_{\omega,t,x}$ and $v\in L^2_{\omega,t}H_x^1$ be predictable. For all $j,k\ge0$ and $d=0,1$, we have
\begin{equation}\label{eq:conditional-product-blocks}
\mathcal E_{r_j}\bigl(\Delta_j(hD^dv_k)\bigr)+\mathcal E_{r_j}\bigl(hD^d\Delta_jv_k\bigr)\le C\|h\|_{L^\infty_{\omega,t,x}}r_k^{-d}\mathcal E_{r_j}(v_k).
\end{equation}
If, in addition, $h\in L^\infty_{\omega,t}W^{1,\infty}(\mathbb T^n)$, then
\begin{equation}\label{eq:conditional-lipschitz-commutator}
\mathcal E_{r_j}\bigl([\Delta_j,h]Dv_k\bigr)
\le C\frac{r_j}{r_k}\|Dh\|_{L^\infty_{\omega,t,x}}\mathcal E_{r_j}(v_k).
\end{equation}
The constant $C$ is independent of $j,k$.
\end{lemma}
\begin{proof}
To simplify the notation, we omit $(\omega,t)$ in the following spatial estimates. The identities
\[
\Delta_jq=\psi_j*q,
\qquad D^d\Delta_jv_k=\psi_j*D^dv_k
\]
give, for $d=0,1$,
\[
\begin{aligned}
|\Delta_j(hD^dv_k)|
&\le \|h\|_{L^\infty_{\omega,t,x}}\bigl(|\psi_j|*|D^dv_k|\bigr),\\
|hD^d\Delta_jv_k|
&\le \|h\|_{L^\infty_{\omega,t,x}}\bigl(|\psi_j|*|D^dv_k|\bigr).
\end{aligned}
\]
If $h\in L^\infty_{\omega,t}W^{1,\infty}(\mathbb T^n)$, then
$|h(x-z)-h(x)|\le\|Dh\|_{L^\infty_{\omega,t,x}}\rho(z)$, so
\[
\begin{aligned}
|[\Delta_j,h]Dv_k(x)|
&=\left|\int_{\mathbb T^n}\psi_j(z)[h(x-z)-h(x)]Dv_k(x-z)\,dz\right|\\
&\le\|Dh\|_{L^\infty_{\omega,t,x}}\bigl[(\rho|\psi_j|)*|Dv_k|\bigr](x).
\end{aligned}
\]

By Lemma~\ref{lem:block-kernels}, with $\nu=0$ and $\nu=1$, respectively,
\[
\|\psi_j\|_{L^1}\le C,
\qquad
\|\rho|\psi_j|\|_{L^1}
=\int_{\mathbb T^n}\rho(z)|\psi_j(z)|\,dz\le Cr_j.
\]
Applying Lemma~\ref{lem:calculus} with kernel $|\psi_j|$ to the two product bounds and adding yields
\[
\begin{aligned}
&\mathcal E_{r_j}\bigl(\Delta_j(hD^dv_k)\bigr)
+\mathcal E_{r_j}\bigl(hD^d\Delta_jv_k\bigr)\\
&\qquad\le 2\|h\|_{L^\infty_{\omega,t,x}}\|\psi_j\|_{L^1}
\mathcal E_{r_j}(D^dv_k)
\le C\|h\|_{L^\infty_{\omega,t,x}}\mathcal E_{r_j}(D^dv_k).
\end{aligned}
\]
For the commutator, applying Lemma~\ref{lem:calculus} with kernel $\rho|\psi_j|$ gives
\[
\begin{aligned}
\mathcal E_{r_j}\bigl([\Delta_j,h]Dv_k\bigr)
&\le\|Dh\|_{L^\infty_{\omega,t,x}}\|\rho|\psi_j|\|_{L^1}
\mathcal E_{r_j}(Dv_k)\\
&\le Cr_j\|Dh\|_{L^\infty_{\omega,t,x}}\mathcal E_{r_j}(Dv_k).
\end{aligned}
\]

Finally, $v_k=\Delta_kv$ has Fourier support in $\{|\xi|\le2r_k^{-1}\}$. The differentiation estimate~\eqref{eq:old-2-18} in Lemma~\ref{lem:window-scales}, applied at scale $r_j$, gives
\[
\mathcal E_{r_j}(Dv_k)\le Cr_k^{-1}\mathcal E_{r_j}(v_k).
\]
Substituting this bound into the preceding product and commutator estimates proves~\eqref{eq:conditional-product-blocks} for $d=1$ and~\eqref{eq:conditional-lipschitz-commutator}, while the case $d=0$ follows directly from $D^0v_k=v_k$.
\end{proof}

For each integer $k\ge0$, define
\[
S_kh:=
\begin{cases}
\displaystyle\sum_{\ell=0}^{k-4}\Delta_\ell h,&k\ge4,\\
0,&0\le k\le3.
\end{cases}
\]
For each coefficient $h$, write
\[
h=S_kh+R_kh,
\qquad R_kh:=h-S_kh.
\]
The low-frequency factor $S_kh$ ensures that $\Delta_j((S_kh)D^dv_k)=0$ for $d=0,1$ whenever $|j-k|\ge4$. The remainder $R_kh$ is estimated using the spatial regularity of $h$.

For $k\ge4$, the multiplier of $S_k$ is
\[
\sum_{\ell=0}^{k-4}m_\ell(\xi)=\chi(2^{-(k-4)}\xi).
\]
The corresponding convolution kernel has an $L^1$ bound independent of $k$. Hence, for all $k\ge0$,
\begin{equation}\label{eq:coefficient-low-cutoff}
\|S_kh\|_{L^\infty_{\omega,t,x}}\le C\|h\|_{L^\infty_{\omega,t,x}},
\qquad
\|DS_kh\|_{L^\infty_{\omega,t,x}}\le C\|Dh\|_{L^\infty_{\omega,t,x}}.
\end{equation}
Applying the H\"older--Besov characterization in Lemma~\ref{lem:frequency} with $\beta=1+\alpha$, we obtain
\begin{equation}\label{eq:coefficient-tail-bound}
\|R_kh\|_{L^\infty_{\omega,t,x}}
\le\sum_{\ell\ge\max\{0,k-3\}}\|\Delta_\ell h\|_{L^\infty_{\omega,t,x}}
\le Cr_k^{1+\alpha}\|h\|_{L^\infty_{\omega,t}C_x^{1+\alpha}}.
\end{equation}
Choose a fixed integer $c_0\ge4$ large enough for the Fourier support conditions in Lemmas~\ref{lem:low-coefficient-commutator} and~\ref{lem:coefficient-remainder}. For the low-frequency part $S_kh$ of the coefficient, only $|k-j|\le c_0$ contributes; hence the ratio $r_j/r_k$ in~\eqref{eq:conditional-lipschitz-commutator} is bounded.

\begin{lemma}\label{lem:low-coefficient-commutator}
Let $h\in L^\infty_{\omega,t}C_x^1$ and $v\in L^2_{\omega,t}H_x^1$ be predictable, with $\|h\|_{L^\infty_{\omega,t}C_x^1}\le K$. For all $j,k\ge0$ with $|k-j|>c_0$, we have
\[
\Delta_j\bigl((S_kh)D^dv_k\bigr)=0\quad(d=0,1),
\qquad [\Delta_j,S_kh]Dv_k=0.
\]
For all $j,k\ge0$ with $|k-j|\le c_0$, the commutator satisfies
\begin{equation}\label{eq:old-2-35}
\mathcal E_{r_j}\bigl([\Delta_j,S_kh]Dv_k\bigr)
\le C\mathcal E_{r_j}(v_k),
\end{equation}
where $C=C(K)$ is independent of $j,k$.
\end{lemma}
\begin{proof}
If $k\le3$, then $S_k=0$. If $k>3$,
\[
\begin{aligned}
\operatorname{supp}\widehat{S_kh}
&\subset\{\xi\in\mathbb Z^n:|\xi|\le2^{k-3}\},\\
\operatorname{supp}\widehat v_k
&\subset\{\xi\in\mathbb Z^n:2^{k-1}\le|\xi|\le2^{k+1}\}.
\end{aligned}
\]
Thus $(S_kh)D^dv_k$ has Fourier support in $2^{k-2}\le|\xi|\le2^{k+2}$ for $d=0,1$. When $|k-j|>c_0$, both $\Delta_j((S_kh)D^dv_k)$ and $\Delta_jv_k$ vanish, so $[\Delta_j,S_kh]Dv_k=0$.

Apply~\eqref{eq:conditional-lipschitz-commutator} with coefficient $S_kh$ and use~\eqref{eq:coefficient-low-cutoff}. Since $r_j/r_k=2^{k-j}\le2^{c_0}$ for the remaining indices, this proves~\eqref{eq:old-2-35}.
\end{proof}

For the remainder, high coefficient frequencies can contribute to $\Delta_j(R_kh\,D^dv_k)$ even when $k\le j-c_0$. The next lemma estimates these terms.

\begin{lemma}\label{lem:coefficient-remainder}
Let $h\in L^\infty_{\omega,t}C_x^{1+\alpha}$ and $v\in L^2_{\omega,t}H_x^1$ be predictable, with $\|h\|_{L^\infty_{\omega,t}C_x^{1+\alpha}}\le K$. For all $j,k\ge0$ and $d=0,1$, we have
\begin{equation}\label{eq:old-2-36}
\mathcal E_{r_j}\bigl(\Delta_j(R_kh\,D^dv_k)\bigr)
\le Cc_{jk}^{(d)}\mathcal E_{r_j}(v_k),
\end{equation}
where
\[
c_{jk}^{(d)}:=
\begin{cases}
r_j^{1+\alpha}r_k^{-d},&k\le j,\\
r_k^{1+\alpha-d},&k>j.
\end{cases}
\]
Moreover, $R_kh\,D^d\Delta_jv_k=0$ when $|k-j|>c_0$, while for $|k-j|\le c_0$,
\begin{equation}\label{eq:old-2-37}
\mathcal E_{r_j}\bigl(R_kh\,D^d\Delta_jv_k\bigr)
\le Cr_k^{1+\alpha-d}\mathcal E_{r_j}(v_k).
\end{equation}
The constant $C=C(K)$ is independent of $j,k$.
\end{lemma}

\begin{proof}[Proof of Lemma~\ref{lem:coefficient-remainder}]
Suppose first that $k\le j-c_0$. Expand the coefficient remainder as
\[
R_kh=\sum_{\ell\ge\max\{0,k-3\}}\Delta_\ell h.
\]
If $\Delta_j((\Delta_\ell h)D^dv_k)\ne0$, there are frequencies $\eta$ of $\Delta_\ell h$ and $\zeta$ of $v_k$ such that
\[
|\zeta|\le2^{k+1}\le2^{j-3},
\qquad 2^{j-1}\le|\eta+\zeta|\le2^{j+1},
\]
where we used $c_0\ge4$. The triangle inequality gives $2^{j-2}\le|\eta|\le2^{j+2}$. Since $j\ge4$, this excludes $\Delta_0h$, whose frequencies satisfy $|\eta|\le2$. For $\ell\ge1$, comparison with $2^{\ell-1}\le|\eta|\le2^{\ell+1}$ yields
\[
2^{\ell-1}\le2^{j+2},\qquad 2^{j-2}\le2^{\ell+1}
\quad\Longrightarrow\quad j-3\le\ell\le j+3.
\]
Thus $\Delta_j((\Delta_\ell h)D^dv_k)=0$ for $|\ell-j|>3$. Summing over the remaining coefficient blocks yields
\[
\begin{aligned}
\mathcal E_{r_j}\bigl(\Delta_j(R_kh\,D^dv_k)\bigr)
&\le Cr_k^{-d}\sum_{|\ell-j|\le3}\|\Delta_\ell h\|_{L^\infty_{\omega,t,x}}\mathcal E_{r_j}(v_k)\\
&\le Cr_j^{1+\alpha}r_k^{-d}\mathcal E_{r_j}(v_k).
\end{aligned}
\]
The first inequality uses~\eqref{eq:conditional-product-blocks}; the second uses $\|\Delta_\ell h\|_{L^\infty_{\omega,t,x}}\le Cr_\ell^{1+\alpha}$ and $r_\ell\asymp r_j$ for $|\ell-j|\le3$.

If $k>j-c_0$, apply the product estimate directly to the whole remainder:
\[
\begin{aligned}
\mathcal E_{r_j}\bigl(\Delta_j(R_kh\,D^dv_k)\bigr)
&\le C\|R_kh\|_{L^\infty_{\omega,t,x}}r_k^{-d}\mathcal E_{r_j}(v_k)\\
&\le Cr_k^{1+\alpha-d}\mathcal E_{r_j}(v_k),
\end{aligned}
\]
where the last inequality is~\eqref{eq:coefficient-tail-bound}. If $j-c_0<k\le j$, then
\[
r_k^{1+\alpha-d}
=\left(\frac{r_k}{r_j}\right)^{1+\alpha}r_j^{1+\alpha}r_k^{-d}
\le2^{c_0(1+\alpha)}r_j^{1+\alpha}r_k^{-d}.
\]
Together with the estimate for $k\le j-c_0$, this proves the first case of~\eqref{eq:old-2-36} for every $k\le j$. For $k>j$, the remainder bound above gives the second case directly.

For the product with $\Delta_j$ acting on $v_k$, the same estimates give
\[
\begin{aligned}
\mathcal E_{r_j}\bigl(R_kh\,D^d\Delta_jv_k\bigr)
&\le C\|R_kh\|_{L^\infty_{\omega,t,x}}r_k^{-d}\mathcal E_{r_j}(v_k)\\
&\le Cr_k^{1+\alpha-d}\mathcal E_{r_j}(v_k),
\end{aligned}
\]
which is~\eqref{eq:old-2-37}. Moreover, disjoint block supports imply
\[
R_kh\,D^d\Delta_jv_k
=R_kh\,D^d\Delta_j\Delta_kv=0,\qquad |k-j|>c_0,
\]
where $\Delta_j$ acts on $v_k$ alone, whereas in~\eqref{eq:old-2-36} it acts on the product $R_kh\,D^dv_k$.
\end{proof}

To prove Proposition~\ref{prop:commutators}, we apply the preceding lemmas with $(h,d)=(\sigma,1)$ for $\sigma:Dv$ and $(h,d)=(\gamma,0)$ for $\gamma\cdot v$. In particular, splitting $\sigma=S_k\sigma+R_k\sigma$ in each block of $v=\sum_kv_k$ gives the following decomposition of the commutator $R_j^\sigma$ defined in~\eqref{eq:old-2-27}:
\begin{equation}\label{eq:old-2-34}
\begin{aligned}
[\Delta_j,\sigma]Dv
={}&\sum_k[\Delta_j,S_k\sigma]Dv_k\\
&+\sum_k\Delta_j\bigl[(R_k\sigma):Dv_k\bigr]
-\sum_k(R_k\sigma):D\Delta_jv_k.
\end{aligned}
\end{equation}

\begin{proof}[Proof of Proposition~\ref{prop:commutators}]
For the terms involving $u$, the kernel formula and H\"older product estimate give
\[
\left\{
\begin{aligned}
&\|R_j^a\|_{L^\infty_{\omega,t,x}}
\le Cr_j^\alpha\|D^2u\|_{L^\infty_{\omega,t,x}},\\
&\|\Delta_j(b\cdot Du+cu)\|_{L^\infty_{\omega,t,x}}
\le Cr_j^\alpha\|u\|_{L^\infty_{\omega,t}C_x^{1+\alpha}}.
\end{aligned}
\right.
\]
Indeed,
\[
R_j^a(x)=\int_{\mathbb T^n}\psi_j(z)[a(x-z)-a(x)]:D^2u(x-z)\,dz,
\]
so the first bound follows from the $\alpha$-moment of $\psi_j$ in Lemma~\ref{lem:block-kernels}. Applying~\eqref{eq:old-2-16} in Lemma~\ref{lem:holder-interpolation} with parameter $\delta/2$ and taking essential suprema in $(\omega,t)$ gives
\[
\|D^2u\|_{L^\infty_{\omega,t,x}}+\|u\|_{L^\infty_{\omega,t}C_x^{1+\alpha}}
\le\delta\|u\|_{L^\infty_{\omega,t}C_x^{2+\alpha}}
+C_\delta\|u\|_{L^\infty_{\omega,t,x}}.
\]
Combining this with the two bounds above and the triangle inequality proves~\eqref{eq:perturbation-u}.

For the martingale commutator, the scale comparison of Lemma~\ref{lem:window-scales} reads
\begin{equation}\label{eq:old-2-28}
\mathcal E_{r_j}(v_k)\le C
\begin{cases}
e_k,&k\le j,\\
\dfrac{r_j}{r_k}e_k,&k>j.
\end{cases}
\end{equation}
For $|k-j|\le c_0$, the ratio $r_j/r_k$ is bounded by $2^{c_0}$, so $\mathcal E_{r_j}(v_k)\le Ce_k$.

We claim that
\begin{equation}\label{eq:old-2-30}
\mathcal E_{r_j}(R_j^\sigma)
\le C\biggl[\sum_{|k-j|\le c_0}e_k
+r_j^{1+\alpha}\sum_{k\le j}r_k^{-1}e_k
+r_j\sum_{k>j}r_k^{\alpha-1}e_k\biggr].
\end{equation}
It suffices to prove the estimates for $v^N:=\sum_{k=0}^Nv_k$ uniformly in $N$, since $v^N\to v$ in $L^2_{\omega,t}H_x^1$ and conditional Fatou gives the passage to the limit. By the scale comparison~\eqref{eq:old-2-28}, to prove~\eqref{eq:old-2-30} it suffices to establish the following estimate at scale $r_j$:
\begin{equation}\label{eq:commutator-fixed-window}
\mathcal E_{r_j}(R_j^\sigma)\le C\biggl[\sum_{|k-j|\le c_0}\mathcal E_{r_j}(v_k)+r_j^{1+\alpha}\sum_{k\le j}r_k^{-1}\mathcal E_{r_j}(v_k)+\sum_{k>j}r_k^\alpha\mathcal E_{r_j}(v_k)\biggr].
\end{equation}
We estimate the three sums in~\eqref{eq:old-2-34} separately. Lemma~\ref{lem:low-coefficient-commutator}, with $h=\sigma$, gives
\[
\sum_{k\ge0}\mathcal E_{r_j}\bigl([\Delta_j,S_k\sigma]Dv_k\bigr)
\le C\sum_{|k-j|\le c_0}\mathcal E_{r_j}(v_k).
\]
For the second sum, use~\eqref{eq:old-2-36} with $h=\sigma$ and $d=1$:
\[
\sum_{k\ge0}\mathcal E_{r_j}\bigl(\Delta_j[(R_k\sigma):Dv_k]\bigr)
\le Cr_j^{1+\alpha}\sum_{k\le j}r_k^{-1}\mathcal E_{r_j}(v_k)
+C\sum_{k>j}r_k^\alpha\mathcal E_{r_j}(v_k).
\]
For the last sum,~\eqref{eq:old-2-37} and $r_k^\alpha\le1$ give
\[
\sum_{k\ge0}\mathcal E_{r_j}\bigl((R_k\sigma):D\Delta_jv_k\bigr)
\le C\sum_{|k-j|\le c_0}r_k^\alpha\mathcal E_{r_j}(v_k)
\le C\sum_{|k-j|\le c_0}\mathcal E_{r_j}(v_k).
\]
Adding these bounds and using the triangle inequality in~\eqref{eq:old-2-34} gives~\eqref{eq:commutator-fixed-window} for finite sums of blocks.

For $\Delta_j(\gamma\cdot v)$, use
\[
\Delta_j(\gamma\cdot v)
=\sum_k\Delta_j((S_k\gamma)\cdot v_k)
+\sum_k\Delta_j((R_k\gamma)\cdot v_k).
\]
The first sum contains only $|k-j|\le c_0$, and its $\mathcal E_{r_j}$ norm is bounded by $C\sum_{|k-j|\le c_0}e_k$. For the second sum, Lemma~\ref{lem:coefficient-remainder} with $h=\gamma$ and $d=0$, followed by the scale comparison~\eqref{eq:old-2-28}, gives
\[
\begin{aligned}
\sum_{k\ge0}\mathcal E_{r_j}\bigl(\Delta_j((R_k\gamma)\cdot v_k)\bigr)
&\le Cr_j^{1+\alpha}\sum_{k\le j}\mathcal E_{r_j}(v_k)
+C\sum_{k>j}r_k^{1+\alpha}\mathcal E_{r_j}(v_k)\\
&\le Cr_j^{1+\alpha}\sum_{k\le j}e_k
+Cr_j\sum_{k>j}r_k^\alpha e_k.
\end{aligned}
\]
Since $r_k\le1$, the right-hand side is bounded by the last two terms on the right-hand side of~\eqref{eq:old-2-30}.

We now prove~\eqref{eq:perturbation-v-energy} and~\eqref{eq:perturbation-v-small}. Substituting $e_k\le r_k^{2+\alpha}\|v\|_{\mathcal E}$ into the sums over $|k-j|\le c_0$ and $k>j$ in~\eqref{eq:old-2-30} gives
\[
\sum_{|k-j|\le c_0}e_k\le Cr_j^{2+\alpha}\|v\|_{\mathcal E},
\qquad r_j\sum_{k>j}r_k^{\alpha-1}e_k\le Cr_j^{2+2\alpha}\|v\|_{\mathcal E}.
\]
Thus
\begin{equation}\label{eq:compact-commutator-sum}
r_j^{-1-\alpha}\left[e_j+\mathcal E_{r_j}(R_j^\sigma)
+\mathcal E_{r_j}(\Delta_j(\gamma\cdot v))\right]
\le C\left[r_j\|v\|_{\mathcal E}+\sum_{k\le j}r_k^{-1}e_k\right].
\end{equation}
The bound $\sum_{k\le j}r_k^{-1}e_k\le C\|v\|_{\mathcal E}$ proves~\eqref{eq:perturbation-v-energy}.

To prove~\eqref{eq:perturbation-v-small}, fix $N$ and use~\eqref{eq:finite-fourier-block} and the definition~\eqref{eq:global-conditional-energy} of $\mathcal B_1(v)$ for $k\le N$:
\[
\sum_{k\le j}r_k^{-1}e_k
\le C_N\mathcal B_1(v)+\|v\|_{\mathcal E}\sum_{k>N}r_k^{1+\alpha}
\le C_N\mathcal B_1(v)+Cr_N^{1+\alpha}\|v\|_{\mathcal E}.
\]
Choose $N$ large and then $J>N$ large so that $C(r_J+r_N^{1+\alpha})\le\varepsilon$. Equation~\eqref{eq:compact-commutator-sum} gives the desired bound for $j\ge J$. For $j<J$, the finite Fourier bounds imply
\[
\|v_j(s)\|_\infty+\|R_j^\sigma(s)\|_\infty
+\|\Delta_j(\gamma\cdot v)(s)\|_\infty\le C_J\|v_s\|_{H_x^1}.
\]
Indeed, $R_j^\sigma=\Delta_j(\sigma:Dv)-\sigma:D\Delta_jv$, and boundedness of $\sigma,\gamma$ controls each term. Conditional integration bounds these finitely many normalized blocks by $C_J\mathcal B_1(v)$. This proves~\eqref{eq:perturbation-v-small}.
\end{proof}

\section{Sobolev solvability and comparison}\label{app:sobolev}
We prove the Sobolev solvability, conditional estimates, and comparison principle used in Section~\ref{subsec:sobolev}.

\begin{proof}[Proof of Lemma~\ref{lem:periodic-existence}]
The principal coefficients satisfy~\cite[Section~5.1, Assumptions~5.1--5.3]{DuQiuTang2012} at Sobolev regularity index zero, while boundedness of $b,c,\gamma$ and Sobolev interpolation verify the generator condition in~\cite[Assumption~5.4]{DuQiuTang2012}. The whole-space proof of~\cite[Theorem~5.5]{DuQiuTang2012} carries over to $\mathbb T^n$ using Fourier series and a finite smooth partition of unity. Together with the standard $H_x^1$ time-continuity argument based on the Hilbert-space It\^o formula~\cite{Krylov2013}, this gives a unique Sobolev strong solution satisfying, for every deterministic $t\in[0,T]$,
\begin{equation}\label{eq:sobolev-trace-estimate}
\begin{aligned}
&\E\sup_{t\le s\le T}\|u_s\|_{H_x^1}^2
+\E\int_t^T\bigl(\|u_s\|_{H_x^2}^2+\|v_s\|_{H_x^1}^2\bigr)\,ds\\
&\le C\left(\E\|g\|_{H_x^1}^2
+\E\int_t^T\|f_s\|_{L_x^2}^2\,ds\right),
\end{aligned}
\end{equation}
where $C=C(\lambda,K)$ is independent of $t$. For more details, see~\cite[proof of Theorem~5.5]{DuQiuTang2012}.

We now prove the conditional estimates.
Put
\[
e_s:=\|u_s\|_{H_x^2}^2+\|v_s\|_{H_x^1}^2.
\]
\emph{Step 1: Deterministic starting times.}
Fix $t\in[0,T]$ and $B\in\mathcal F_t$. On $[t,T]$, linearity implies that $(\mathbf1_Bu,\mathbf1_Bv)$ solves the same equation with data $(\mathbf1_Bf,\mathbf1_Bg)$. Estimate~\eqref{eq:sobolev-trace-estimate} gives
\[
\E\left[\mathbf1_B\int_t^T e_s\,ds\right]
\le C\E\left[\mathbf1_B\left(\|g\|_{H_x^1}^2
+\int_t^T\|f_s\|_{L_x^2}^2\,ds\right)\right].
\]
Since this holds for every $B\in\mathcal F_t$, the defining property of conditional expectation yields
\[
\E_t\int_t^T e_s\,ds
\le C\E_t\left[\|g\|_{H_x^1}^2
+\int_t^T\|f_s\|_{L_x^2}^2\,ds\right].
\]

\smallskip\noindent\emph{Step 2: Stopping times.}
First let $\tau$ take finitely many values $t_1,\ldots,t_N$. For $B\in\mathcal F_\tau$, apply Step 1 to the events $B\cap\{\tau=t_i\}\in\mathcal F_{t_i}$ and sum over $i$. This gives
\begin{equation}\label{eq:simple-stopping-sobolev}
\E\left[\mathbf1_B\int_\tau^T e_s\,ds\right]
\le C\E\left[\mathbf1_B\left(\|g\|_{H_x^1}^2
+\int_\tau^T\|f_s\|_{L_x^2}^2\,ds\right)\right].
\end{equation}
For a general stopping time $\tau$, choose finite-valued stopping times $\tau_k\downarrow\tau$. If $B\in\mathcal F_\tau$, then $B\in\mathcal F_{\tau_k}$ for every $k$, so~\eqref{eq:simple-stopping-sobolev} applies with $\tau_k$ in place of $\tau$. Monotone convergence yields~\eqref{eq:simple-stopping-sobolev} at $\tau$, for every $B\in\mathcal F_\tau$. This proves~\eqref{eq:conditional-sobolev}.

\smallskip\noindent\emph{Step 3: Uniformly bounded data norms.}
If $f\in L^\infty_\omega L_t^2L_x^2$ and $g\in L^\infty_\omega H_x^1$, then, for every stopping time $\tau$,
\[
\|g\|_{H_x^1}^2+\int_\tau^T\|f_s\|_{L_x^2}^2\,ds
\le\|g\|_{L^\infty_\omega H_x^1}^2
+\|f\|_{L^\infty_\omega L_t^2L_x^2}^2
\quad\text{a.s.}
\]
Substituting into~\eqref{eq:conditional-sobolev} gives a deterministic bound for the conditional integral of $\|v_s\|_{H_x^1}^2$, uniform in $\tau$. This is~\eqref{eq:old-2-23}.
\end{proof}

\begin{proof}[Proof of Lemma~\ref{lem:periodic-sobolev}]
Set $w:=\bar u-u$. The approximation and comparison argument of~\cite[Theorem~6.3]{DuQiuTang2012} applies on the torus. Since the difference equation has nonnegative forcing $\bar f-f$ and terminal datum $\bar g-g$, it gives
\[
\mathbb P\bigl(\|w_t^-\|_{L_x^2}=0\bigr)=1
\qquad\text{for each fixed }t\in[0,T],
\]
where $h^-:=\max\{-h,0\}$.

To show that the comparison holds almost surely for every time, let $\mathcal T:=[0,T]\cap\mathbb Q$ and define
\[
\Omega_0:=\{w\in C([0,T];L_x^2)\}
\cap\bigcap_{q\in\mathcal T}\{\|w_q^-\|_{L_x^2}=0\}.
\]
This event has probability one. For $\omega\in\Omega_0$ and $t\in[0,T]$, choose $q_j\in\mathcal T$ with $q_j\to t$. Then
\[
\|w_t^-\|_{L_x^2}
=\|w_t^--w_{q_j}^-\|_{L_x^2}
\le\|w_t-w_{q_j}\|_{L_x^2}\longrightarrow0.
\]
Thus $u_t\le\bar u_t$ almost everywhere in space for every $t\in[0,T]$, on the same event $\Omega_0$.
\end{proof}

\section*{Declaration of generative AI and AI-assisted technologies in the manuscript preparation process}
During the preparation of this work, the authors used OpenAI's Codex for assistance with language editing, manuscript organization, and the presentation and checking of intermediate mathematical arguments. The authors reviewed and edited the output as needed and take full responsibility for the content of the published article.

\bibliographystyle{siam}
\bibliography{ref}

@book{GilbargTrudinger2001,
  author    = {Gilbarg, David and Trudinger, Neil S.},
  title     = {Elliptic Partial Differential Equations of Second Order},
  edition   = {2},
  series    = {Classics in Mathematics},
  publisher = {Springer},
  address   = {Berlin},
  year      = {2001},
  note      = {Reprint of the 1998 edition},
  doi       = {10.1007/978-3-642-61798-0},
  isbn      = {978-3-540-41160-4},
  url       = {https://link.springer.com/book/10.1007/978-3-642-61798-0},
}

@article{Bensoussan1983PartialObservation,
  author    = {Bensoussan, Alain},
  title     = {Maximum principle and dynamic programming approaches of the optimal control of partially observed diffusions},
  journal   = {Stochastics},
  volume    = {9},
  number    = {3},
  pages     = {169--222},
  year      = {1983},
  doi       = {10.1080/17442508308833253},
}

@article{Bensoussan1983Distributed,
  author    = {Bensoussan, Alain},
  title     = {Stochastic maximum principle for distributed parameter systems},
  journal   = {J. Franklin Inst.},
  volume    = {315},
  number    = {5--6},
  pages     = {387--406},
  year      = {1983},
  doi       = {10.1016/0016-0032(83)90059-5},
}

@article{Bayer2022,
  author    = {Bayer, Christian and Qiu, Jinniao and Yao, Yao},
  title     = {Pricing options under rough volatility with backward {SPDE}s},
  journal   = {SIAM J. Financial Math.},
  volume    = {13},
  number    = {1},
  pages     = {179--212},
  year      = {2022},
  doi       = {10.1137/20M1357639},
}

@article{BenderDokuchaev2016Swing,
  author    = {Bender, Christian and Dokuchaev, Nikolai},
  title     = {A first-order {BSPDE} for swing option pricing},
  journal   = {Math. Finance},
  volume    = {26},
  number    = {3},
  pages     = {461--491},
  year      = {2016},
  doi       = {10.1111/mafi.12067},
}

@article{TangZhang2009,
  author    = {Tang, Shanjian and Zhang, Xu},
  title     = {Null controllability for forward and backward stochastic parabolic equations},
  journal   = {SIAM J. Control Optim.},
  volume    = {48},
  number    = {4},
  pages     = {2191--2216},
  year      = {2009},
  doi       = {10.1137/050641508},
}

@article{QiuWei2014,
  author    = {Qiu, Jinniao and Wei, Wenning},
  title     = {On the quasi-linear reflected backward stochastic partial differential equations},
  journal   = {J. Funct. Anal.},
  volume    = {267},
  number    = {10},
  pages     = {3598--3656},
  year      = {2014},
  doi       = {10.1016/j.jfa.2014.08.023},
}

@article{YangTang2013,
  author    = {Yang, Zhou and Tang, Shanjian},
  title     = {{Dynkin} game of stochastic differential equations with random coefficients and associated backward stochastic partial differential variational inequality},
  journal   = {SIAM J. Control Optim.},
  volume    = {51},
  number    = {1},
  pages     = {64--95},
  year      = {2013},
  doi       = {10.1137/110850980},
}

@article{Du2019,
  author    = {Du, Kai and Liu, Jiakun},
  title     = {On the {C}auchy problem for stochastic parabolic equations in {H}\"older spaces},
  journal   = {Trans. Amer. Math. Soc.},
  volume    = {371},
  number    = {4},
  pages     = {2643--2664},
  year      = {2019},
  doi       = {10.1090/tran/7533},
}

@article{Du2010,
  author    = {Du, Kai and Meng, Qingxin},
  title     = {A revisit to {$W_2^n$}-theory of super-parabolic backward stochastic partial differential equations in {$\mathbb R^d$}},
  journal   = {Stochastic Process. Appl.},
  volume    = {120},
  number    = {10},
  pages     = {1996--2015},
  year      = {2010},
  doi       = {10.1016/j.spa.2010.06.001},
}

@article{DuQiuTang2012,
  author    = {Du, Kai and Qiu, Jinniao and Tang, Shanjian},
  title     = {{$L^p$} theory for super-parabolic backward stochastic partial differential equations in the whole space},
  journal   = {Appl. Math. Optim.},
  volume    = {65},
  number    = {2},
  pages     = {175--219},
  year      = {2012},
  doi       = {10.1007/s00245-011-9154-9},
}

@article{Du2012Dirichlet,
  author    = {Du, Kai and Tang, Shanjian},
  title     = {Strong solution of backward stochastic partial differential equations in {$C^2$} domains},
  journal   = {Probab. Theory Related Fields},
  volume    = {154},
  number    = {1--2},
  pages     = {255--285},
  year      = {2012},
  doi       = {10.1007/s00440-011-0369-0},
}

@article{Hu2002,
  author    = {Hu, Ying and Ma, Jin and Yong, Jiongmin},
  title     = {On semi-linear degenerate backward stochastic partial differential equations},
  journal   = {Probab. Theory Related Fields},
  volume    = {123},
  number    = {3},
  pages     = {381--411},
  year      = {2002},
  doi       = {10.1007/s004400100193},
}

@article{Ma1997,
  author    = {Ma, Jin and Yong, Jiongmin},
  title     = {Adapted solution of a degenerate backward {SPDE}, with applications},
  journal   = {Stochastic Process. Appl.},
  volume    = {70},
  number    = {1},
  pages     = {59--84},
  year      = {1997},
  doi       = {10.1016/S0304-4149(97)00057-4},
}

@article{Ma1999,
  author    = {Ma, Jin and Yong, Jiongmin},
  title     = {On linear, degenerate backward stochastic partial differential equations},
  journal   = {Probab. Theory Related Fields},
  volume    = {113},
  number    = {2},
  pages     = {135--170},
  year      = {1999},
  doi       = {10.1007/s004400050205},
}

@article{Peng1992,
  author    = {Peng, Shige},
  title     = {Stochastic {H}amilton-{J}acobi-{B}ellman equations},
  journal   = {SIAM J. Control Optim.},
  volume    = {30},
  number    = {2},
  pages     = {284--304},
  year      = {1992},
  doi       = {10.1137/0330018},
}

@article{tang2016cauchy,
  author    = {Tang, Shanjian and Wei, Wenning},
  title     = {On the {C}auchy problem for backward stochastic partial differential equations in {H}\"{o}lder spaces},
  journal   = {Ann. Probab.},
  volume    = {44},
  number    = {1},
  pages     = {360--398},
  year      = {2016},
  doi       = {10.1214/14-AOP976},
}

@article{Ye2025,
  author    = {Ye, Yuyang and Li, Yunzhang and Tang, Shanjian},
  title     = {Fractional backward stochastic partial differential equations with applications to stochastic optimal control of partially observed systems driven by {L\'{e}vy} processes},
  journal   = {SIAM J. Control Optim.},
  volume    = {63},
  number    = {2},
  pages     = {1314--1347},
  year      = {2025},
  doi       = {10.1137/24M1692666},
}

@article{Zhao2022,
  author    = {Zhao, Guohuan},
  title     = {{SDEs} with random and irregular coefficients},
  journal   = {Rev. Mat. Iberoam.},
  volume    = {38},
  number    = {3},
  pages     = {947--980},
  year      = {2022},
  doi       = {10.4171/RMI/1313},
}

@article{Zhou1993,
  author    = {Zhou, Xun Yu},
  title     = {On the necessary conditions of optimal controls for stochastic partial differential equations},
  journal   = {SIAM J. Control Optim.},
  volume    = {31},
  number    = {6},
  pages     = {1462--1478},
  year      = {1993},
  doi       = {10.1137/0331068},
}

@article{AraiSuzuki2021,
  author    = {Arai, Takuji and Suzuki, Ryoichi},
  title     = {A {Clark--Ocone} type formula via {It\^o} calculus and its application to finance},
  journal   = {J. Stoch. Anal.},
  volume    = {2},
  number    = {4},
  year      = {2021},
  doi       = {10.31390/josa.2.4.05},
  note      = {Article 5},
}

@book{CDLL2019,
  author    = {Cardaliaguet, Pierre and Delarue, Fran\c{c}ois and Lasry, Jean-Michel and Lions, Pierre-Louis},
  title     = {The Master Equation and the Convergence Problem in Mean Field Games},
  series    = {Annals of Mathematics Studies},
  volume    = {201},
  publisher = {Princeton University Press},
  address   = {Princeton, NJ},
  year      = {2019},
  doi       = {10.23943/princeton/9780691190716.001.0001},
}

@misc{CardaliaguetMailletYan2025,
  author    = {Cardaliaguet, Pierre and Maillet, Rapha{\"e}l and Yan, Wenbin},
  title     = {On the long-time behavior of mean field game systems with a common noise},
  year      = {2025},
  howpublished = {arXiv:2509.17443},
  eprint    = {2509.17443},
  archivePrefix = {arXiv},
  url       = {https://arxiv.org/abs/2509.17443},
}

@article{LittlewoodPaley1931,
  author    = {Littlewood, J. E. and Paley, R. E. A. C.},
  title     = {Theorems on {Fourier} series and power series},
  journal   = {J. London Math. Soc.},
  volume    = {6},
  number    = {3},
  pages     = {230--233},
  year      = {1931},
  doi       = {10.1112/jlms/s1-6.3.230},
  url       = {https://academic.oup.com/jlms/article-abstract/s1-6/3/230/963501},
}

@article{GubinelliPerkowski2015,
  author    = {Gubinelli, Massimiliano and Perkowski, Nicolas},
  title     = {Lectures on singular stochastic {PDEs}},
  journal   = {Ensaios Mat.},
  volume    = {29},
  pages     = {1--89},
  year      = {2015},
  doi       = {10.21711/217504322015/em291},
}

@article{GubinelliImkellerPerkowski2015,
  author    = {Gubinelli, Massimiliano and Imkeller, Peter and Perkowski, Nicolas},
  title     = {Paracontrolled distributions and singular {PDEs}},
  journal   = {Forum Math. Pi},
  volume    = {3},
  pages     = {e6},
  year      = {2015},
  doi       = {10.1017/fmp.2015.2},
}

@article{DieningHastoRoudenko2009,
  author    = {Diening, Lars and H\"ast\"o, Peter and Roudenko, Svetlana},
  title     = {Function spaces of variable smoothness and integrability},
  journal   = {J. Funct. Anal.},
  volume    = {256},
  number    = {6},
  pages     = {1731--1768},
  year      = {2009},
  doi       = {10.1016/j.jfa.2009.01.017},
}

@article{Du2016,
  author    = {Du, Kai and Liu, Jiakun},
  title     = {A {Schauder} estimate for stochastic {PDEs}},
  journal   = {C. R. Math. Acad. Sci. Paris},
  volume    = {354},
  number    = {4},
  pages     = {371--375},
  year      = {2016},
  doi       = {10.1016/j.crma.2016.01.010},
}

@article{Krylov2013,
  author    = {Krylov, N. V.},
  title     = {A relatively short proof of {It\^o}'s formula for {SPDEs} and its applications},
  journal   = {Stoch. Partial Differ. Equ. Anal. Comput.},
  volume    = {1},
  number    = {1},
  pages     = {152--174},
  year      = {2013},
  doi       = {10.1007/s40072-013-0003-5},
}

@misc{MollaQiu2026,
  author    = {Molla, H. U. and Qiu, Jinniao},
  title     = {A policy iteration scheme for semilinear stochastic {Hamilton--Jacobi--Bellman} equations with exponential convergence},
  year      = {2026},
  howpublished = {arXiv:2607.29024},
  eprint    = {2607.29024},
  archivePrefix = {arXiv},
  url       = {https://arxiv.org/abs/2607.29024},
}

@misc{DelarueMou2025,
  author    = {Delarue, Fran\c{c}ois and Mou, Chenchen},
  title     = {Major-minor mean field games: common noise helps},
  year      = {2025},
  howpublished = {arXiv:2501.02627},
  eprint    = {2501.02627},
  archivePrefix = {arXiv},
  url       = {https://arxiv.org/abs/2501.02627},
}

@article{Lorenzi2000,
  author    = {Lorenzi, Luca},
  title     = {Optimal {Schauder} estimates for parabolic problems with data measurable with respect to time},
  journal   = {SIAM J. Math. Anal.},
  volume    = {32},
  number    = {3},
  pages     = {588--615},
  year      = {2000},
  doi       = {10.1137/S0036141098342842},
}

@incollection{RuzhanskyTurunen2007,
  author    = {Ruzhansky, Michael and Turunen, Ville},
  title     = {On the {Fourier} analysis of operators on the torus},
  editor    = {Toft, Joachim and Wong, M. W. and Zhu, Hongmei},
  booktitle = {Modern Trends in Pseudo-Differential Operators},
  series    = {Operator Theory: Advances and Applications},
  volume    = {172},
  pages     = {87--105},
  publisher = {Birkh{\"a}user},
  address   = {Basel},
  year      = {2007},
  doi       = {10.1007/978-3-7643-8116-5_5},
}
\end{document}